# INTRODUCTION TO n-ADAPTIVE FUZZY MODELS TO ANALYZE PUBLIC OPINION ON AIDS

**W. B. Vasantha Kandasamy**
**Florentin Smarandache**

Translation of interviews
from Tamil by

**Meena Kandasamy**

**February 2006**

# INTRODUCTION TO n-ADAPTIVE FUZZY MODELS TO ANALYZE PUBLIC OPINION ON AIDS

**W. B. Vasantha Kandasamy**
e-mail: vasantha@iitm.ac.in
web: http://mat.iitm.ac.in/~wbv

**Florentin Smarandache**
e-mail: smarand@unm.edu

Translation of interviews
from Tamil by

**Meena Kandasamy**

**February 2006**

## CONTENTS

## Preface

"AIDS is not simply a physical malady, it is also an artifact of social and sexual transgression, violated taboo, fractured identity—political and personal projections. Its key words are primarily the property of the powerful. *AIDS: Keywords* – is my attempt to identify and contest some of the assumptions underlying our current 'knowledge'. In this effort I am joined by many AIDS activists including people living with AIDS—Acquired Immuno Deficiency Syndrome.

"A syndrome is a pattern of symptoms pointing to a "morbid state" which may or may not be caused by infectious agents; a disease, on the other hand is, "any deviation from or interruption of the normal structure or function of any part, organ or system of the body that is manifested by a characteristic set of symptoms or signs and whose etiology, pathology and prognosis may be known or unknown". In other words, a syndrome points to or signifies the underlying disease process(es), a disease on the other hand is constituted in and by those process(es). The syndrome AIDS, in other words, cannot be communicated nor can the opportunistic infections that constitute the syndrome be readily communicated to those with healthy immune systems. Misunderstanding AIDS for a disease is one of our cultures profoundest confusions of a signifier for a sign. We keep pushing the signifying chain toward that ultimate sign—our collective mortality. Diseases, we are taught are communicable. When AIDS is identified as disease many consequences follows -— not the least of which is wide spread public terror about "Catching" AIDS from people in public places or during casual contact. The communicability of AIDS was fixed in the public and medical memory by easy public health reports. AIDS was repeatedly compared to hepatitis B-virus since both appeared to be blood borne and sexually transmitted."

Jan Zita Grover writes so in his article on *AIDS: Keywords* in "The State of the Language" (1990) ed. Christopher Ricks and Leonard Michaels, Faber and Faber, London. In the same book, in his article "*Speaking in the Shadow of AIDS*", Wayne Koestenbaun says, "From the start, AIDS posed vocabulary difficulties. Even when AIDS is disclosed it means more than simply a sickness. The word AIDS implies sociological categories — a population, an inclination, a failure, a primary affiliation, that is assumed to pave the way to sickness. The word AIDS thus becomes an accusation and an act of housekeeping, it clarifies the caste system of the culture it assaults.

"When speaking of AIDS, words cease to offer documentation; words become the makers of fear and hesitation. It is not simple to interpret contradictory data to phrase each physical sensation or to calculate the nearness of death. Though AIDS, as we know it, is inseparable from the verbal and prejudicial structures that mantle it, it also contains elements of a more archaic, essential dialect, sounds of pain and phrases of regret. Because a person with AIDS facing isolation or eviction may hide the disease and because the syndrome first achieved notoriety among gay men, AIDS—as guilty secret—has been linked, like homosexuality itself, to unspeakability.

"In the era of Silence = Death, sly or shy signs cease to be revolutionary, even if they promise to save lives. Though AIDS remains closely associated with silence what is finally, silent about it? Michel Foucault, who died of AIDS, gave us a theoretical framework to explain how private language and life are webbed by public forces. We may think we control words, as if they were neutral transparently instrumental but language is larger than our efforts to wrangle it into reflecting our desires. Think of all the words spent in the name of AIDS—preventing it, fearing it, spreading our ignorance around like music. AIDS may present us with an excuse for declaring our alienation from public language and yet the suffering and sorrow we have learned to call "AIDS" only shows how our bodies are veined by coursing powers that seems too remote but is nearby…

Because language is multiple the word 'AIDS' always means more and differently than we intend.

"There is no possibility of an independently elected language every word is stained by community AIDS may seem to have little to do with subjects and verbs; but it has everything to do with the acts that depend on language—renting an apartment, hiring a lawyer, voting… The syndrome may not have a will or a mission but it has been imprinted with the directives and prejudices of the societies where it flourishes; so AIDS like hunger, like wars is a affair of money and its transport, a matter of unjust circulation."

Michael Callen, in the article on "*AIDS: The Linguistic Battlefield*" writes, "The language of AIDS depends upon many presumptions, piled one upon another like a house of cards when you happen to doubt a particular pre-assumption often the entire superstructure collapses making communication difficult if not impossible. Three commonly used terms, which betray no evidence of the enormously complex; presumptions encoded within them are the cause the AIDS virus and HIV disease! … Houses of HIV-infected children have been burned down. There are people currently in jail and others under trial for being HIV antibody positive and having sex without telling their partners. In Cuba, HIV antibody positive individuals are quarantined for life." We have mainly given these short excerpts of three persons mainly to see how we can make some changes in our society so that the stigma that the term HIV/AIDS holds is at least lessened, if not made to entirely disappear.

A public debate on AIDS is virtually absent because it is viewed as something that only affects 'others'—people who deviate from the existing 'moral' codes of society such as gays, sex workers, intravenous drug users and so on. Its virulent spread among the entire population (including the heterosexual community) has become a major cause for concern.

Persons with AIDS have faced isolation and discrimination and AIDS itself has been constantly viewed with hostility and fear. Panic and an unhealthy silence have been the two extreme reactions of public hysteria over AIDS. Any person who is affected with AIDS is immediately labeled 'characterless.' Widespread advertisements have to a great extent increased the awareness levels. But social acceptance of the disease has not been similarly forthcoming. A recent survey estimated that 57% of Indian women are aware of AIDS. This clearly shows that there is a vast rural population of women who are simply unaware of the existence of AIDS or the information of how it spreads. Thus it becomes imperative for NGOs and governmental agencies to extend awareness programmes to rural areas. Thus we have to build a bridge between public and patients (who are in fact and in actuality not divided entities but merely consider themselves so). We have interviewed 101 people from all walks of life and we have carefully analyzed their feelings using 2 adaptive fuzzy model.

There are many fuzzy models like Fuzzy matrices, Fuzzy Cognitive Maps, Fuzzy relational Maps, Fuzzy Associative Memories, Bidirectional Associative memories and so on. But almost all these models can give only one sided solution like hidden pattern or a resultant output vector dependent on the input vector depending in the problem at hand. So for the first time we have defined a n-adaptive fuzzy model which can view or analyze the problem in n ways ($n \geq 2$). Though we have defined these n- adaptive fuzzy models theorectically we are not in a position to get a n-adaptive fuzzy model for $n > 2$ for practical real world problems.

We have used the 2-adaptive fuzzy model having the two fuzzy models, fuzzy matrices model and BAMs viz. model to analyze the views of public about HIV/ AIDS disease, patient and the awareness program. This book has five chapters and 6 appendices. The first chapter just recalls the definition of four fuzzy models used in this book and gives illustration of some of them. Chapter two introduces the new n-adaptive fuzzy models. Chapter three uses for the first time 2 adaptive fuzzy models to study psychological and sociological problems about HIV/AIDS. Chapter four gives an outline of the interviews. Chapter five gives the suggestions and conclusion based on our study. Of the 6 appendices four of them are C-program made to make the working of the fuzzy model simple.

We thank Dr.K.Kandasamy for his intellectual support and unstinted cooperation that made this book possible.

W.B.VASANTHA KANDASAMY
FLORENTIN SMARANDACHE
18-02-2006

## Chapter One

# SOME BASIC FUZZY MODELS WITH ILLUSTRATIONS

This chapter has four sections. In the first section we just recall the notion of fuzzy matrices or Combined Effect Time Dependent data (CETD) matrix and illustrate it with examples. In section two the definition and illustration of the fuzzy Relational Maps is given. Third section gives the basic concept of Bidirectional Associative memories (BAM) models and illustrates it with an example. In the forth section we describe the Fuzzy Associative Memories (FAM) model. As all these four models are used by us in chapter II and III. We have just recalled them in order to make this book a self contained one.

## 1.1 Fuzzy matrices and their applications

In this section we just recall the definition of fuzzy matrices and illustrate with a model. Matrix theory has become a very simple but an effective tool in the analysis of collected raw data, be it from a transportation department or from migrant labourers with HIV/AIDS or public opinion about HIV/AIDS. Fuzzy matrix or CETD matrix model is the one, which helps to analyze the raw data. The analysis is carried out in five stages. In the first stage the collected raw data, which is time dependent is made into a Initial Raw Data (IRD) matrix. By taking along the rows the time period (it may be age group, hours of the day depending on the problem under analysis).

The columns of the initial Raw data matrix correspond to the attributes or concepts carefully chosen by the expert.

The IRD matrix may not in general be uniform for the length of the interval of time period may not be the same. In the second stage the IRD matrix is transformed into an Average Time Dependent data (ATD) matrix ($a_{ij}$).

So, in the second stage, in order to obtain an unbiased uniform effect on each and every data so collected, we transform this initial matrix into an Average Time-dependent Data (ATD) matrix ($a_{ij}$). To make the calculations simpler, in the third stage, we use the simple average techniques and convert the above Average Time-dependent Data (ATD) matrix with entries $e_{ij}$ where $e_{ij} \in \{-1, 0, 1\}$. We name this matrix as the Refined Time-dependent Data (RTD) matrix.

$e_{ij}$ are calculated using the formula if

$$a_{ij} \leq (\mu_j - \alpha * \sigma_j) \text{ then } e_{ij} = -1 \text{ else if}$$
$$a_{ij} \in (\mu_j - \alpha * \sigma_j, \mu_j + \alpha * \sigma_j) \text{ then } e_{ij} = 0 \text{ else if}$$
$$a_{ij} \geq (\mu_j + \alpha * \sigma_j) \text{ then } e_{ij} = 1$$

where $a_{ij}$ are entries of the ATD matrix and $\mu_j$ and $\sigma_j$ are the mean and the standard deviation of the $j^{th}$ column of the ATD matrix respectively. $\alpha$ is a parameter taken from the interval [0, 1]. The value of the $e_{ij}$ corresponding to each entry is determined in a special way. At the next stage using the RTD matrix we get the Combined Effect Time Dependent data matrix (CETD matrix), which gives the cumulative effect of all these entries. In the final stage we obtain the row sums of the CETD matrix. A program in C language is written which easily estimates all these five stages. The C program shall give details of the mathematical working and method of calculation. Now we give the description of the problem and the proposed solution to the problem.

So now we will work with by redefining the age group to six intervals and draw conclusions from the model.

Let

$A_1$ – CSWs
$A_2$ – Other women
$A_3$ – Smoking
$A_4$ – Alcohol
$A_5$ – Bad Company
$A_6$ – Quacks.

The initial raw data matrix of the 55 migrant labourers affected with HIV/AIDS relating their age group and the diseases are given by the following matrix.

Initial Raw Data Matrix

| | $A_1$ | $A_2$ | $A_3$ | $A_4$ | $A_5$ | $A_6$ |
|---|---|---|---|---|---|---|
| 20-23 | 3 | 2 | 2 | 3 | 3 | 2 |
| 24-30 | 20 | 6 | 15 | 16 | 14 | 10 |
| 31-34 | 15 | 4 | 13 | 14 | 9 | 6 |
| 35-37 | 8 | 2 | 6 | 6 | 3 | 2 |
| 38-40 | 6 | 1 | 4 | 4 | 2 | 1 |
| 41-47 | 8 | 1 | 5 | 6 | 3 | 2 |

ATD Matrix

| | $A_1$ | $A_2$ | $A_3$ | $A_4$ | $A_5$ | $A_6$ |
|---|---|---|---|---|---|---|
| 20-23 | 0.75 | 0.5 | 0.5 | 0.75 | 0.75 | 0.5 |
| 24-30 | 3.33 | 1 | 2.5 | 2.67 | 2.33 | 1.67 |
| 31-34 | 3.75 | 1 | 3.25 | 3.5 | 2.25 | 1.5 |
| 35-37 | 2.67 | 0.67 | 2 | 2 | 1 | 0.67 |
| 38-40 | 2 | 0.33 | 1.33 | 1.33 | 0.67 | 0.28 |
| 41-47 | 1.14 | 0.14 | 0.71 | 0.86 | 0.48 | 0.28 |

Average and Standard Deviation of the data

| | | | | | | |
|---|---|---|---|---|---|---|
| Average | 2.27 | 0.61 | 1.72 | 1.85 | 1.24 | 0.83 |
| Standard Deviation | 1.09 | 0.32 | 0.97 | 0.99 | 0.76 | 0.55 |

The RTD matrix for $\alpha = 0.5$

$$\begin{bmatrix} -1 & 0 & -1 & -1 & -1 & -1 \\ 1 & 1 & 1 & 1 & 1 & 1 \\ 1 & 1 & 1 & 1 & 1 & 1 \\ 0 & 0 & 0 & 0 & 0 & 0 \\ 0 & -1 & 0 & -1 & -1 & -1 \\ -1 & -1 & -1 & -1 & -1 & -1 \end{bmatrix}$$

Row sum of the RTD or fuzzy matrix is as follows

$$\begin{array}{r} 20-23 \\ 24-30 \\ 31-34 \\ 35-37 \\ 38-40 \\ 41-47 \end{array}\begin{bmatrix} -5 \\ 6 \\ 6 \\ 0 \\ -4 \\ -6 \end{bmatrix}$$

Graph of the row sum of the RTD matrix is given for is $\alpha = 0.5$

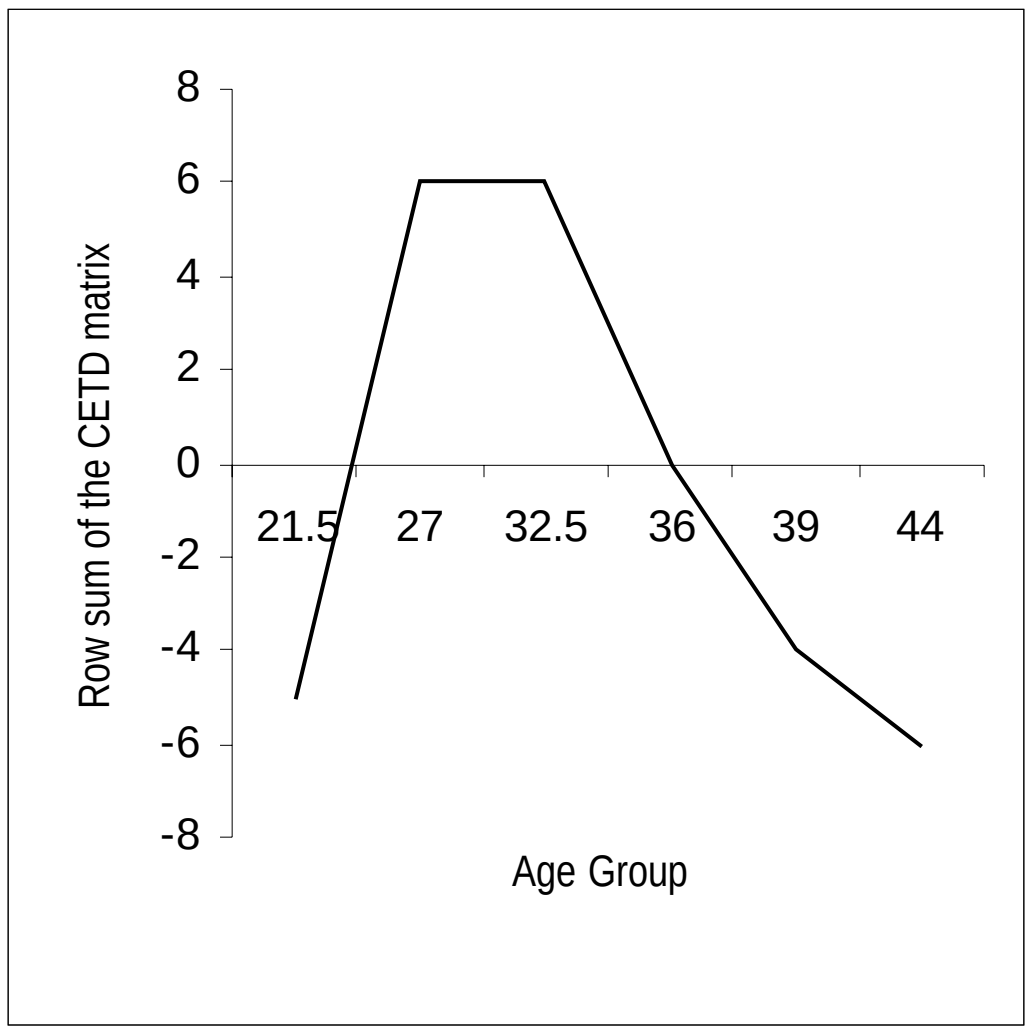

**Graph 1.1**

The maximum vulnerable age group being 24-30 and 31-34.

The RTD matrix for the parameter $\alpha = 0.2$ is given below:

$$\begin{bmatrix} -1 & -1 & -1 & -1 & -1 & -1 \\ 1 & 1 & 1 & 1 & 1 & 1 \\ 1 & 1 & 1 & 1 & 1 & 1 \\ 1 & 0 & 1 & 0 & -1 & -1 \\ 0 & -1 & -1 & -1 & -1 & -1 \\ -1 & -1 & -1 & -1 & -1 & -1 \end{bmatrix}$$

The row sum of the RTD matrix

$$\begin{matrix} 20-23 \\ 24-30 \\ 31-34 \\ 35-37 \\ 38-40 \\ 41-47 \end{matrix}\begin{bmatrix} -6 \\ 6 \\ 6 \\ 0 \\ -5 \\ -6 \end{bmatrix}$$

We see that there are nil or very negligible HIV/AIDS patients in the age group 20-23, as its value is –6.

The graph for the parameter α = 0.2 of the RTD matrix is given below:

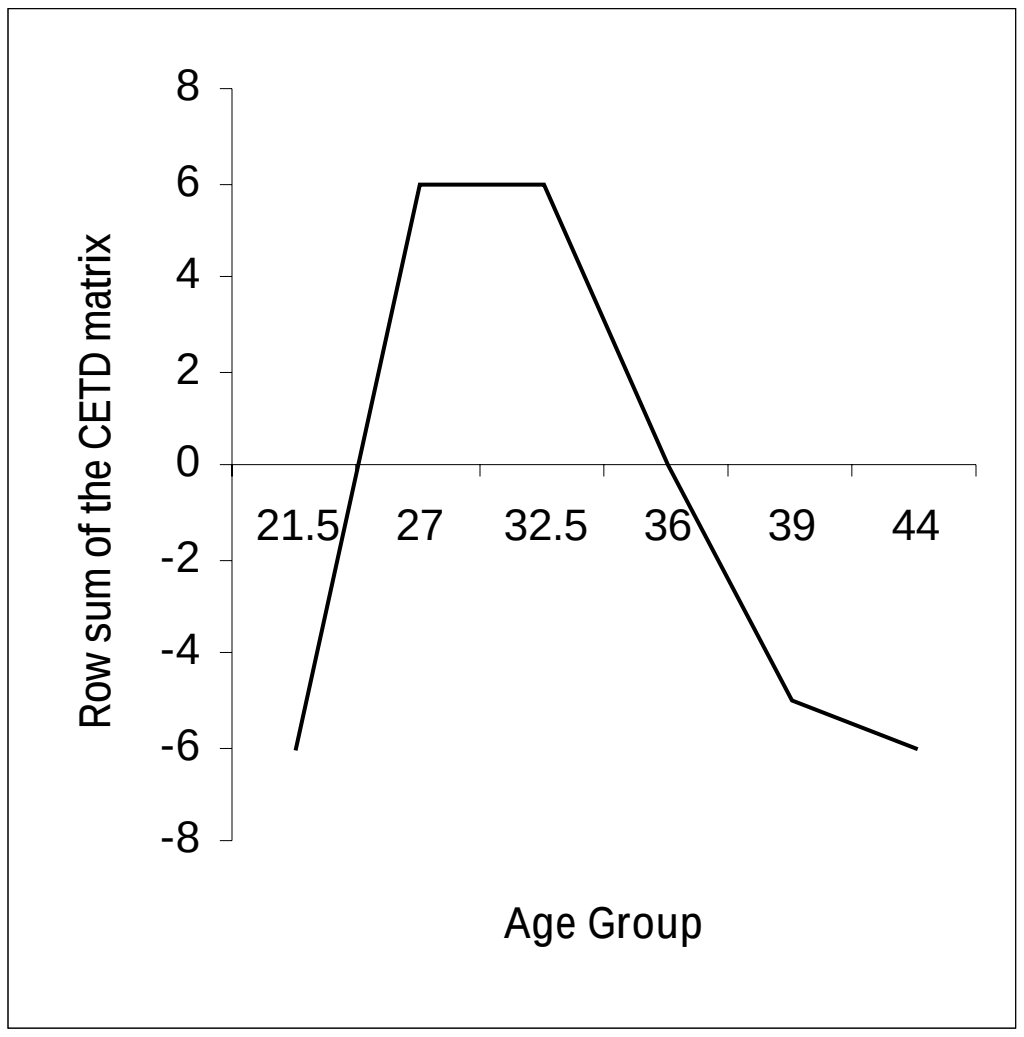

**Graph 1.2**

For α = 1, the RTD matrix is as follows

$$\begin{bmatrix} -1 & 0 & -1 & -1 & 0 & 0 \\ 0 & 1 & 0 & 0 & 1 & 1 \\ 1 & 1 & 1 & 1 & 1 & 1 \\ 0 & 0 & 0 & 0 & 0 & 0 \\ 0 & 0 & 0 & 0 & 0 & 0 \\ -1 & -1 & -1 & -1 & -1 & -1 \end{bmatrix}$$

Row sum of RTD matrix

$$\begin{matrix} 20-23 \\ 24-30 \\ 31-34 \\ 35-37 \\ 38-40 \\ 41-47 \end{matrix} \begin{bmatrix} -3 \\ 3 \\ 6 \\ 0 \\ 0 \\ -6 \end{bmatrix}$$

Graph of the RTD matrix for α = 1

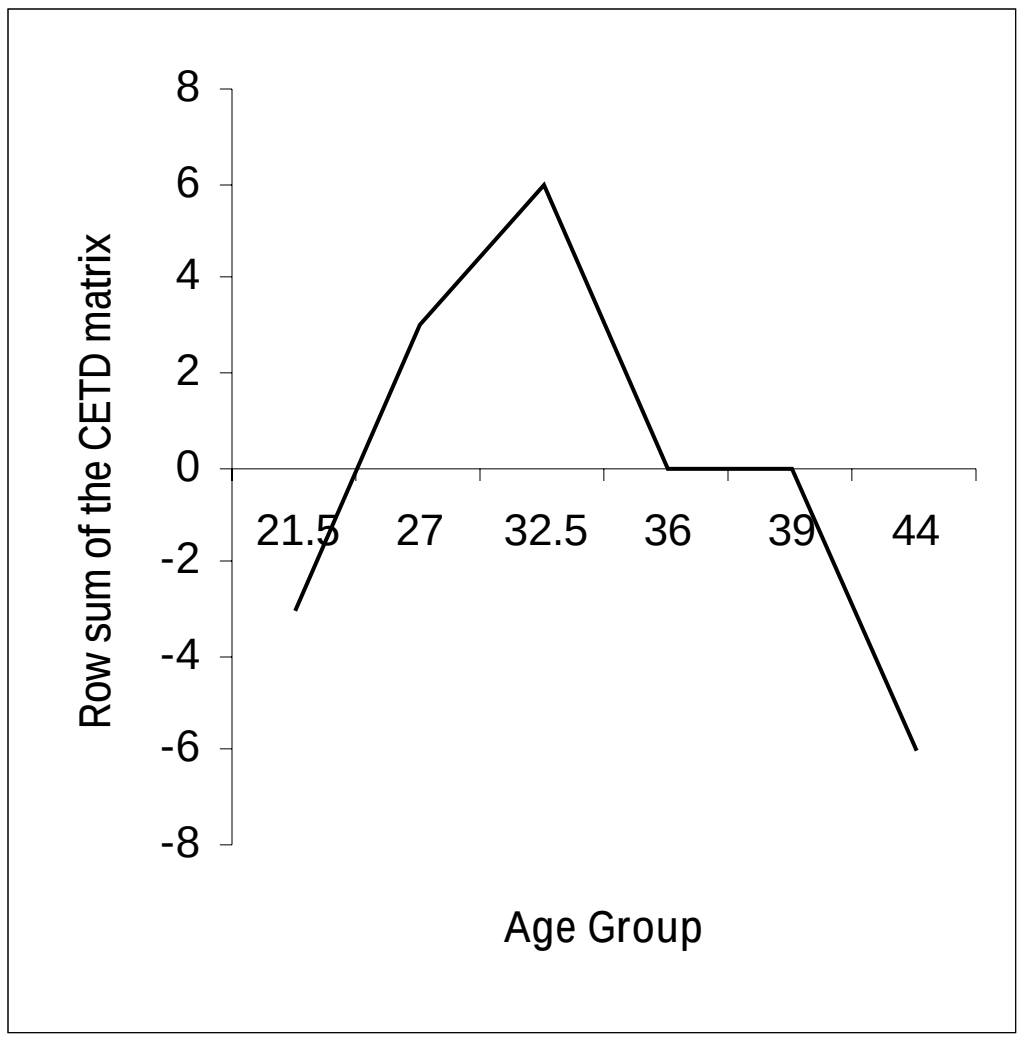

**Graph 1.3**

CETD matrix

$$\begin{bmatrix} -3 & -1 & -3 & -3 & -2 & -2 \\ 2 & 3 & 2 & 2 & 3 & 3 \\ 3 & 3 & 3 & 3 & 3 & 3 \\ 1 & 0 & 1 & 0 & -1 & -1 \\ 0 & -2 & -1 & -2 & -2 & -2 \\ -3 & -3 & -3 & -3 & -3 & -3 \end{bmatrix}$$

Row sum of CETD Matrix

$$\begin{matrix} 20-23 \\ 24-30 \\ 31-34 \\ 35-37 \\ 38-40 \\ 41-47 \end{matrix}\begin{bmatrix} -14 \\ 15 \\ 18 \\ 0 \\ -9 \\ -18 \end{bmatrix}$$

Graph of CETD Matrix

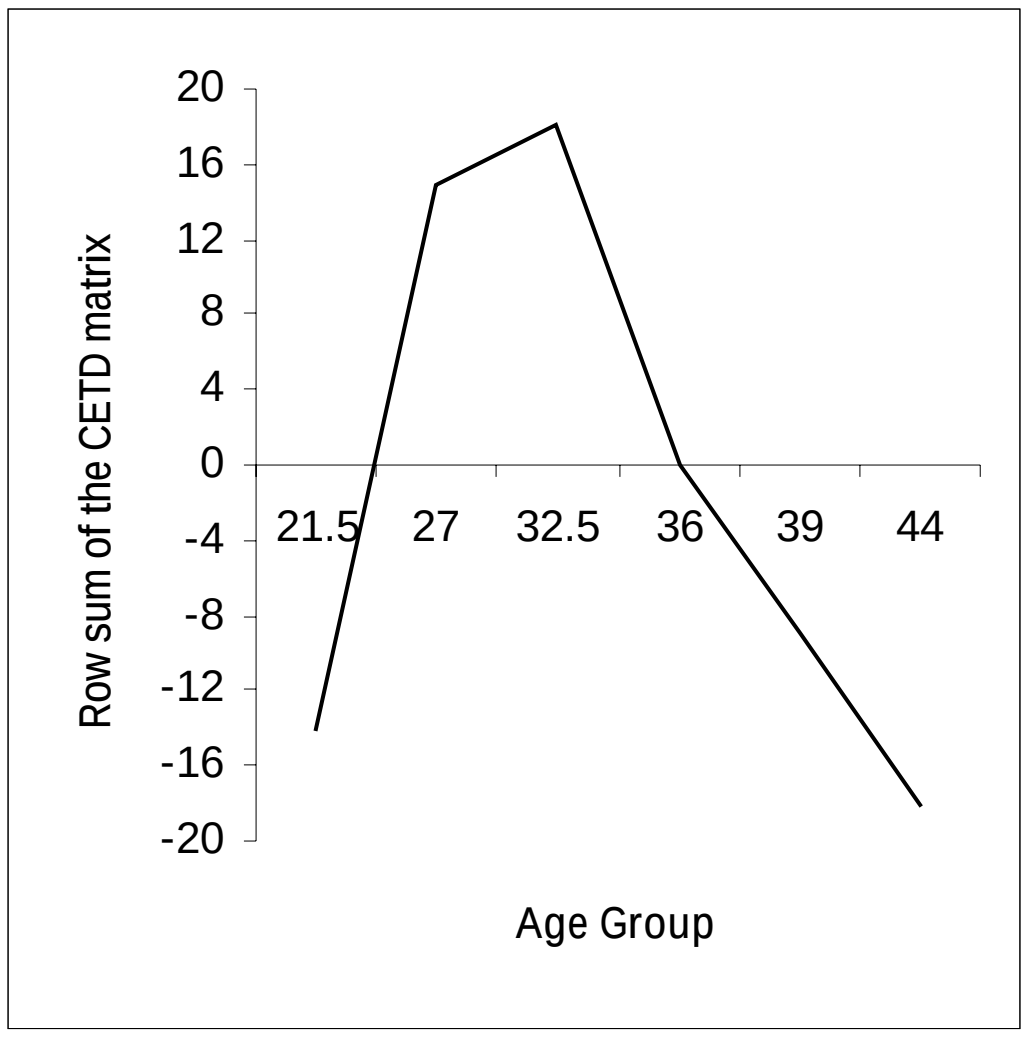

**Graph 1.4**

This graph shows that the age group 21.5 to 39.5 are most vulnerable to HIV/AIDS among migrant labourers the peak falling in the interval 31-34. The C- program given in Appendix 3 can be used to make calculations simple.

## 1.2 Definition and illustration of Fuzzy Relational Maps (FRMs)

In this section, we introduce the notion of Fuzzy Relational Maps (FRMs); they are constructed analogous to FCMs described and discussed in the earlier sections. In FCMs we promote the correlations between causal associations among concurrently active units. But in FRMs we divide the very causal associations into two disjoint units, for example, the relation between a teacher and a student or relation between an employee or employer or a relation between doctor and patient and so on. Thus for us to define a FRM we need a domain space and a range space which are disjoint in the sense of concepts. We further assume no intermediate relation exists within the domain elements or node and the range space elements. The number of elements in the range space need not in general be equal to the number of elements in the domain space.

Thus throughout this section we assume the elements of the domain space are taken from the real vector space of dimension n and that of the range space are real vectors from the vector space of dimension m (m in general need not be equal to n). We denote by R the set of nodes $R_1,\ldots, R_m$ of the range space, where $R = \{(x_1,\ldots, x_m) \mid x_j = 0 \text{ or } 1 \}$ for $j = 1, 2,\ldots, m$. If $x_i = 1$ it means that the node $R_i$ is in the on state and if $x_i = 0$ it means that the node $R_i$ is in the off state. Similarly D denotes the nodes $D_1, D_2,\ldots, D_n$ of the domain space where $D = \{(x_1,\ldots, x_n) \mid x_j = 0 \text{ or } 1\}$ for $i = 1, 2,\ldots, n$. If $x_i = 1$ it means that the node $D_i$ is in the on state and if $x_i = 0$ it means that the node $D_i$ is in the off state.

Now we proceed on to define a FRM.

**DEFINITION 1.2.1:** *A FRM is a directed graph or a map from D to R with concepts like policies or events etc, as nodes and*

*causalities as edges. It represents causal relations between spaces D and R .*

*Let $D_i$ and $R_j$ denote that the two nodes of an FRM. The directed edge from $D_i$ to $R_j$ denotes the causality of $D_i$ on $R_j$ called relations. Every edge in the FRM is weighted with a number in the set {0, ±1}. Let $e_{ij}$ be the weight of the edge $D_iR_j$, $e_{ij} \in \{0, \pm 1\}$. The weight of the edge $D_i$ $R_j$ is positive if increase in $D_i$ implies increase in $R_j$ or decrease in $D_i$ implies decrease in $R_j$ ie causality of $D_i$ on $R_j$ is 1. If $e_{ij} = 0$, then $D_i$ does not have any effect on $R_j$ . We do not discuss the cases when increase in $D_i$ implies decrease in $R_j$ or decrease in $D_i$ implies increase in $R_j$ .*

**DEFINITION 1.2.2:** *When the nodes of the FRM are fuzzy sets then they are called fuzzy nodes. FRMs with edge weights {0, ±1} are called simple FRMs.*

**DEFINITION 1.2.3:** *Let $D_1, \ldots, D_n$ be the nodes of the domain space D of an FRM and $R_1, \ldots, R_m$ be the nodes of the range space R of an FRM. Let the matrix E be defined as $E = (e_{ij})$ where $e_{ij}$ is the weight of the directed edge $D_iR_j$ (or $R_jD_i$), E is called the relational matrix of the FRM.*

***Note:*** It is pertinent to mention here that unlike the FCMs the FRMs can be a rectangular matrix with rows corresponding to the domain space and columns corresponding to the range space. This is one of the marked difference between FRMs and FCMs.

**DEFINITION 1.2.4:** *Let $D_1, ..., D_n$ and $R_1,\ldots, R_m$ denote the nodes of the FRM. Let $A = (a_1,\ldots,a_n)$, $a_i \in \{0, +1\}$. A is called the instantaneous state vector of the domain space and it denotes the on-off position of the nodes at any instant. Similarly let $B = (b_1,\ldots, b_m)$ $b_i \in \{0, +1\}$. B is called instantaneous state vector of the range space and it denotes the on-off position of the nodes at any instant $a_i = 0$ if $a_i$ is off and $a_i = 1$ if $a_i$ is on for i = 1, 2,…, n Similarly, $b_i = 0$ if $b_i$ is off and $b_i = 1$ if $b_i$ is on, for i= 1, 2,…, m.*

**DEFINITION 1.2.5:** *Let $D_1, \ldots, D_n$ and $R_1,\ldots, R_m$ be the nodes of an FRM. Let $D_iR_j$ (or $R_j$ $D_i$) be the edges of an FRM, $j = 1, 2,\ldots, m$ and $i= 1, 2,\ldots, n$. Let the edges form a directed cycle. An FRM is said to be a cycle if it posses a directed cycle. An FRM is said to be acyclic if it does not posses any directed cycle.*

**DEFINITION 1.2.6:** *An FRM with cycles is said to be an FRM with feedback.*

**DEFINITION 1.2.7:** *When there is a feedback in the FRM, i.e. when the causal relations flow through a cycle in a revolutionary manner, the FRM is called a dynamical system.*

**DEFINITION 1.2.8:** *Let $D_i$ $R_j$ (or $R_j$ $D_i$), $1 \leq j \leq m$, $1 \leq i \leq n$. When $R_i$ (or $D_j$) is switched on and if causality flows through edges of the cycle and if it again causes $R_i$ (or$D_j$), we say that the dynamical system goes round and round. This is true for any node $R_j$ (or $D_i$ ) for $1 \leq i \leq n$, (or $1 \leq j \leq m$). The equilibrium state of this dynamical system is called the hidden pattern.*

**DEFINITION 1.2.9:** *If the equilibrium state of a dynamical system is a unique state vector, then it is called a fixed point. Consider an FRM with $R_1, R_2,\ldots, R_m$ and $D_1, D_2,\ldots, D_n$ as nodes. For example, let us start the dynamical system by switching on $R_1$ (or $D_1$). Let us assume that the FRM settles down with $R_1$ and $R_m$ (or $D_1$ and $D_n$) on, i.e. the state vector remains as (1, 0, …, 0, 1) in R (or 1, 0, 0, … , 0, 1) in D), This state vector is called the fixed point.*

**DEFINITION 1.2.10:** *If the FRM settles down with a state vector repeating in the form*

$$A_1 \to A_2 \to A_3 \to \ldots \to A_i \to A_1 \text{ (or } B_1 \to B_2 \to \ldots \to B_i \to B_1)$$

*then this equilibrium is called a limit cycle.*

Now we give the methods of determining the hidden pattern.

Let $R_1, R_2,\ldots, R_m$ and $D_1, D_2,\ldots, D_n$ be the nodes of a FRM with feedback. Let E be the relational matrix. Let us find a hidden pattern when $D_1$ is switched on i.e. when an input is given as vector $A_1 = (1, 0, \ldots, 0)$ in $D_1$, the data should pass through the relational matrix E. This is done by multiplying $A_1$ with the relational matrix E. Let $A_1E = (r_1, r_2,\ldots, r_m)$, after thresholding and updating the resultant vector we get $A_1 E \in R$. Now let $B = A_1E$ we pass on B into $E^T$ and obtain $BE^T$. We update and threshold the vector $BE^T$ so that $BE^T \in D$. This procedure is repeated till we get a limit cycle or a fixed point.

**DEFINITION 1.2.11:** *Finite number of FRMs can be combined together to produce the joint effect of all the FRMs. Let $E_1,\ldots, E_p$ be the relational matrices of the FRMs with nodes $R_1, R_2,\ldots, R_m$ and $D_1, D_2,\ldots, D_n$, then the Combined FRM (CFRM) is represented by the relational matrix $E = E_1+\ldots+ E_p$.*

Now we give a simple illustration of a FRM, for more about FRMs please refer [130-132, 143, 144].

***Example 1.2.1:*** Let us consider the relationship between the teacher and the student. Suppose we take the domain space as the concepts belonging to the teacher say $D_1,\ldots, D_5$ and the range space denote the concepts belonging to the student say $R_1$, $R_2$ and $R_3$.

We describe the nodes $D_1,\ldots, D_5$ and $R_1$ , $R_2$ and $R_3$ as follows:
Domain Space

$D_1$ – Teaching is good
$D_2$ – Teaching is poor
$D_3$ – Teaching is mediocre
$D_4$ – Teacher is kind
$D_5$ – Teacher is harsh [or rude]

(We can have more concepts like teacher is non-reactive, unconcerned etc.)
Range Space

$R_1$ – Good Student
$R_2$ – Bad Student

$R_3$ – Average Student

The relational directed graph of the teacher-student model is given in Figure 1.2.1.

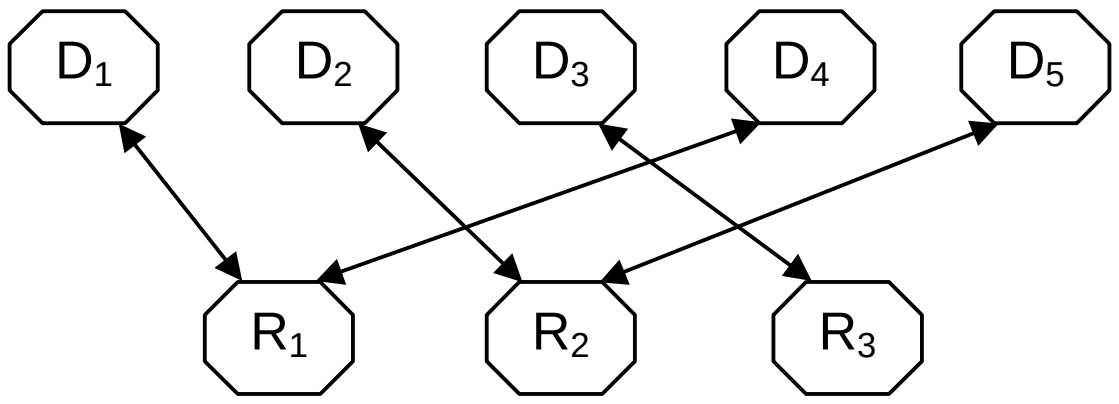

**FIGURE: 1.2.1**

The relational matrix E got from the above map is

$$E = \begin{bmatrix} 1 & 0 & 0 \\ 0 & 1 & 0 \\ 0 & 0 & 1 \\ 1 & 0 & 0 \\ 0 & 1 & 0 \end{bmatrix}$$

If A = (1 0 0 0 0) is passed on in the relational matrix E, the instantaneous vector, AE = (1 0 0) implies that the student is a good student . Now let AE = B, $BE^T$ = (1 0 0 1 0) which implies that the teaching is good and he / she is a kind teacher. Let $BE^T$ = $A_1$, $A_1E$ = (2 0 0) after thresholding we get $A_1E$ = (1 0 0) which implies that the student is good, so on and so forth.

The C-program given in appendix 4 and 5 will help in finding the hidden pattern easily.

## 1.3 Some Basic concepts of BAM with illustration

In this section BAM model [64] is just described to make the book a self contained one. This section has two subsections. Section 1.3.1 describes the BAM model and section 1.3.2 gives an illustration of the model.

We live in a world of marvelous complexity and variety, a world where events never repeat exactly. Even though events are never exactly the same they are also not completely different. There is a threat of continuity, similarity and predictability that allows us to generalize often correctly from past experience to future events. Neural networks or neuro-computing or brain like computation is based on the wistful hope that we can reproduce at least some of the flexibility and power of the human brain by artificial brains. Neural networks consists of many simple computing elements generally simple non linear summing junctions connected together by connections of varying strength a gross abstraction of the brain which consists of very large number of far more complex neurons connected together with far more complex and far more structured couplings, neural networks architecture cover a wide range.

In one sense every computer is a neural net, because we can view traditional digital logic as constructed from inter connected McCullouch-Pitts neurons. McCullouch-Pitts neurons were proposed in 1943 as models of biological neurons and arranged in networks for a specific purpose of computing logic functions. The problems where artificial neural networks have the most promise are those with a real-world flavor: medical research, signal processing, sociological problems etc.

Neural networks helps to solve these problems with natural mechanisms of generalizations. To over-simplify, suppose we represent an object in a network as a pattern of activation of several units. If a unit or two responds incorrectly the overall pattern stays pretty much of the same, and the network still respond correctly to stimuli when neural networks operate similar inputs naturally produce similar outputs. Most real world perceptual problems have this structure of input-output continuity. The prototype model provides a model for human categorization with a great deal of psychological support. The computational strategy leads to some curious human psychology. For instance in United States people imagine a prototype bird that looks somewhat like a sparrow or a robin. So they learn to judge ‘penguins’ or “ostriches” as “bad” birds because these birds do not resemble the prototype bird even though they are birds "Badness" shows up in a number of ways.

Neural networks naturally develop this kind of category structure. The problems that neural networks solved well and solved poorly were those where human showed comparable strengths and weaknesses in their cognitive computations". For this reason until quite recently most of the study of neural networks has been carried out by psychologists and cognitive scientists who sought models of human cognitive function. Neural networks deal with uncertainty as humans do, not by deliberate design but as a by product of their parallel distributed structure. Because general statements about both human psychology and the structure of the world embed so deeply in both neural networks and fuzzy systems, it is very appropriate to study the psychological effects of HIV/AIDS patients and their influence on public using this theory. Like social customs these assumptions are obvious only if you grew up with them. Both neural networks and fuzzy systems break with the historical tradition prominent in western thought and we can precisely and unambiguously characterize the world, divide into two categories and then manipulate these descriptions according to precise and formal rules. Huang Po, a Buddhist teacher of the ninth century observed that "To make use of your minds to think conceptually is to leave the substance and attach yourself to form", and "from discrimination between this and that a host of demons blazes forth".

### 1.3.1 Some Basic Concepts of BAM

Now in this section we go forth to describe the mathematical structure of the Bidirectional Associative Memories (BAM) model [64]. Neural networks recognize ill defined problems without an explicit set of rules. Neurons behave like functions, neurons transduce an unbounded input activation $x(t)$ at time $t$ into a bounded output signal $S(x(t))$ i.e. Neuronal activations change with time.

Artificial neural networks consists of numerous simple processing units or neurons which can be trained to estimate sampled functions when we do not know the form of the functions. A group of neurons form a field. Neural networks contain many field of neurons. In our text $F_x$ will denote a neuron field, which contains n neurons, and $F_y$ denotes a neuron field, which contains p neurons. The neuronal dynamical system is described by a system of first order differential equations that govern the time-evolution of the neuronal activations or which can be called also as membrane potentials.

$$\dot{x}_i = g_i(X, Y, ...)$$
$$\dot{y}_j = h_j(X, Y, ...)$$

where $\dot{x}_i$ and $\dot{y}_j$ denote respectively the activation time function of the $i^{th}$ neuron in $F_X$ and the $j^{th}$ neuron in $F_Y$. The over dot denotes time differentiation, $g_i$ and $h_j$ are some functions of X, Y, ... where $X(t) = (x_1(t), ..., x_n(t))$ and $Y(t) = (y_1(t), ..., y_p(t))$ define the state of the neuronal dynamical system at time t.

The passive decay model is the simplest activation model, where in the absence of the external stimuli, the activation decays in its resting value

$$\dot{x}_i = x_i$$
and $$\dot{y}_j = y_j$$

The passive decay rate $A_i > 0$ scales the rate of passive decay to the membranes resting potentials $\dot{x}_i = -A_i x_i$. The default rate is $A_i = 1$, i.e. $\dot{x}_i = -A_i x_i$. The membrane time constant $C_i > 0$ scales the time variables of the activation dynamical system. The default time constant is $C_i = 1$. Thus $C_i \dot{x}_i = -A_i x_i$.

The membrane resting potential $P_i$ is defined as the activation value to which the membrane potential equilibrates in the absence of external inputs. The resting potential is an additive constant and its default value is zero. It need not be positive.

$$P_i = C_i \dot{x}_i + A_i x_i$$

$$I_i = \dot{x}_i + x_i$$

is called the external input of the system. Neurons do not compute alone. Neurons modify their state activations with external input and with feed back from one another. Now, how do we transfer all these actions of neurons activated by inputs their resting potential etc. mathematically. We do this using what are called synaptic connection matrices.

Let us suppose that the field $F_X$ with n neurons is synaptically connected to the field $F_Y$ of p neurons. Let $m_{ij}$ be a synapse where the axon from the $i^{th}$ neuron in $F_X$ terminates. $M_{ij}$ can be positive, negative or zero. The synaptic matrix M is a n by p matrix of real numbers whose entries are the synaptic efficacies $m_{ij}$.

The matrix M describes the forward projections from the neuronal field $F_X$ to the neuronal field $F_Y$. Similarly a p by n synaptic matrix N describes the backward projections from $F_Y$ to $F_X$. Unidirectional networks occur when a neuron field synaptically intra connects to itself. The matrix M be a n by n square matrix. A Bidirectional network occur if $M = N^T$ and $N = M^T$. To describe this synaptic connection matrix more simply, suppose the n neurons in the field $F_X$ synaptically connect to the p-neurons in field $F_Y$. Imagine an axon from the $i^{th}$ neuron in $F_X$ that terminates in a synapse $m_{ij}$, that about the $j^{th}$ neuron in $F_Y$. We assume that the real number $m_{ij}$ summarizes the synapse and that $m_{ij}$ changes so slowly relative to activation fluctuations that is constant.

Thus we assume no learning if $m_{ij} = 0$ for all t. The synaptic value $m_{ij}$ might represent the average rate of release of a neuro-transmitter such as norepinephrine. So, as a rate, $m_{ij}$ can be positive, negative or zero.

When the activation dynamics of the neuronal fields $F_X$ and $F_Y$ lead to the overall stable behaviour the bidirectional networks are called as Bidirectional Associative Memories (BAM). As in the analysis of the HIV/AIDS patients relative to the migrancy we state that the BAM model studied presently and predicting the future after a span of 5 or 10 years may not be the same.

For the system would have reached stability and after the lapse of this time period the activation neurons under investigations and which are going to measure the model would be entirely different.

Thus from now onwards more than the uneducated poor the educated rich and the middle class will be the victims of HIV/AIDS. So for this study presently carried out can only give how migration has affected the life style of poor labourer and had led them to be victims of HIV/AIDS.

Further not only a Bidirectional network leads to BAM also a unidirectional network defines a BAM if M is symmetric i.e. $M = M^T$. We in our analysis mainly use BAM which are bidirectional networks. However we may also use unidirectional BAM chiefly depending on the problems under investigations. We briefly describe the BAM model more technically and mathematically.

An additive activation model is defined by a system of n + p coupled first order differential equations that inter connects the fields $F_X$ and $F_Y$ through the constant synaptic matrices M and N.

$$x_i = -A_i x_i + \sum_{j=1}^{p} S_j(y_j) n_{ji} + I_i \qquad (1.3.1)$$

$$y_i = -A_j y_j + \sum_{i=1}^{n} S_i(x_i) m_{ij} + J_j \qquad (1.3.2)$$

$S_i(x_i)$ and $S_j(y_j)$ denote respectively the signal function of the $i^{th}$ neuron in the field $F_X$ and the signal function of the $j^{th}$ neuron in the field $F_Y$.

Discrete additive activation models correspond to neurons with threshold signal functions.

The neurons can assume only two values ON and OFF. ON represents the signal +1, OFF represents 0 or – 1 (– 1 when the representation is bipolar). Additive bivalent models describe asynchronous and stochastic behaviour.

At each moment each neuron can randomly decide whether to change state or whether to emit a new signal given its current activation. The Bidirectional Associative Memory or BAM is a non adaptive additive bivalent neural network. In neural literature the discrete version of the equation (1.3.1) and (1.3.2) are often referred to as BAMs.

A discrete additive BAM with threshold signal functions arbitrary thresholds inputs an arbitrary but a constant synaptic connection matrix M and discrete time steps K are defined by the equations

$$x_i^{k+1} = \sum_{j=1}^{p} S_j(y_j^k) m_{ij} + I_i \qquad (1.3.3)$$

$$y_j^{k+1} = \sum_{i=1}^{n} S_i\left(x_i^k\right) m_{ij} + J_i \qquad (1.3.4)$$

where $m_{ij} \in M$ and $S_i$ and $S_j$ are signal functions. They represent binary or bipolar threshold functions. For arbitrary real valued thresholds $U = (U_1, ..., U_n)$ for $F_X$ neurons and $V = (V_1, ..., V_P)$ for $F_Y$ neurons the threshold binary signal functions corresponds to

$$S_i(x_i^k) = \begin{cases} 1 \quad \text{if} \quad x_i^k > U_i \\ S_i(x_i^{k-1}) \text{ if } x_i^k = U_i \\ 0 \quad \text{if} \quad x_i^k < U_i \end{cases} \qquad (1.3.5)$$

and

$$S_j(x_j^k) = \begin{cases} 1 \quad \text{if} \quad y_j^k > V_j \\ S_j(y_j^{k-1}) \text{ if } y_j^k = V_j \\ 0 \quad \text{if} \quad y_j^k < V_j \end{cases} \qquad (1.3.6).$$

The bipolar version of these equations yield the signal value –1 when $x_i < U_i$ or when $y_j < V_j$. The bivalent signal functions allow us to model complex asynchronous state change patterns. At any moment different neurons can decide whether to compare their activation to their threshold. At each moment any of the 2n subsets of $F_X$ neurons or 2p subsets of the $F_Y$ neurons can decide to change state. Each neuron may randomly decide whether to check the threshold conditions in the equations (1.3.5) and (1.3.6). At each moment each neuron defines a random variable that can assume the value ON(+1) or OFF(0 or –1). The network is often assumed to be deterministic and state changes are synchronous i.e. an entire field of neurons is updated at a time. In case of simple asynchrony only one neuron makes a state change decision at a time. When the subsets represent the entire fields $F_X$ and $F_Y$ synchronous state change results.

In a real life problem the entries of the constant synaptic matrix M depends upon the investigator's feelings. The synaptic matrix is given a weightage according to their feelings. If $x \in F_X$ and $y \in F_Y$ the forward projections from $F_X$ to $F_Y$ is defined by the matrix M. $\{F(x_i, y_j)\} = (m_{ij}) = M$, $1 \le i \le n$, $1 \le j \le p$.

The backward projections is defined by the matrix $M^T$. $\{F(y_i, x_i)\} = (m_{ji}) = M^T$, $1 \le i \le n$, $1 \le j \le p$. It is not always true that the backward projections from $F_Y$ to $F_X$ is defined by the matrix $M^T$.

Now we just recollect the notion of bidirectional stability. All BAM state changes lead to fixed point stability. The property holds for synchronous as well as asynchronous state changes.

A BAM system ($F_X$, $F_Y$, M) is bidirectionally stable if all inputs converge to fixed point equilibria. Bidirectional stability is a dynamic equilibrium. The same signal information flows back and forth in a bidirectional fixed point. Let us suppose that A denotes a binary n-vector and B denotes a binary p-vector. Let A be the initial input to the BAM system. Then the BAM equilibrates to a bidirectional fixed point ($A_f$, $B_f$) as

$$
\begin{array}{lllll}
A & \to & M & \to & B \\
A' & \leftarrow & M^T & \leftarrow & B \\
A' & \to & M & \to & B' \\
A'' & \leftarrow & M^T & \leftarrow & B' \quad \text{etc.} \\
A_f & \to & M & \to & B_f \\
A_f & \leftarrow & M^T & \leftarrow & B_f \quad \text{etc.}
\end{array}
$$

where A', A'', ... and B', B'', ... represents intermediate or transient signal state vectors between respectively A and $A_f$ and B and $B_f$. The fixed point of a Bidirectional system is time dependent.

The fixed point for the initial input vectors can be attained at different times. Based on the synaptic matrix M which is developed by the investigators feelings the time at which bidirectional stability is attained also varies accordingly.

### 1.3.2 Use of BAM Model to Study the Cause of Vulnerability to HIV/AIDS and Factors for Migration

Now the object is to study the levels of knowledge and awareness relating to STD/HIV/AIDS existing among the migrant labourers in Tamil Nadu; and to understand the attitude, risk behavior and promiscuous sexual practice of migrant labourers. This study was mainly motivated from the data collected by us of the 60 HIV/AIDS infected persons who belonged to the category that comes to be defined as migrant labourers. Almost all of them were natives of (remote) villages and had migrated to the city, typically, “in search of jobs”, or because of caste and communal violence. We have noticed how, starting from small villages with hopes and dreams these people had set out to the city, only to succumb to various temptations and finally all their dreams turned into horrid nightmares.

Our research includes probing into areas like: patterns and history of migration work, vulnerabilities and risk exposure in an alien surrounding, ‘new’ sexual practices or attitudes and discrimination, effect of displacement, coping mechanism etc. We also study the new economic policies of liberalization and globalization and how this has affected people to lose their traditional livelihood and sources of local employment, forcing them into migration.

Our study has been conducted among this informal sector mainly because migrant labourers are more vulnerable to HIV/AIDS infection, when compared to the local population for reasons which include easy money, poverty, powerlessness, inaccessibility to health services, unstable life-style such that insecurity, in jobs, lack of skills, alienation from hometown, lack of community bondage. Moreover, migrant labourers are also not organized into trade unions, as a result of which, they are made victims of horrendous exploitation: they are paid less than minimum wages, they don't receive legal protection, they are unaware of worker's rights issues and essentially lack stability.

Their work periods are rarely permanent; they work as short-term unskilled/ semi-skilled contract labourers or as daily wagers; but they earn well for a day and spend it badly without any social binding or savings or investing on their family members or children. A vast majority around 65% of those interviewed were essentially also part of the 'mobile' population, which was wrapped not only in a single migration from native village to metropolitan city, but also involved in jobs like driving trucks, taxis, etc. which gave them increased mobility. We have also analyzed the patients' feelings about the outreach and intervention programs related to HIV/AIDS and we have sought to comprehend the patterns of marginalization that has increased the predisposition of migrants to HIV and other infectious diseases.

A linguistic questionnaire which was drafted and interviews were conducted for 60 HIV/AIDS patients from the hospitals was the main data used in this analysis. Then the questionnaire was transformed into a Bidirectional Associative Memory (BAM) model.

Our sample group consisted of HIV infected migrant labourers whose age group ranged between 20-58 and they were involved in a variety of deregulated labour such as transport or truck drivers, construction labourers, daily wagers or employed in hotels or eateries. We have also investigated the feminization of migration and how women were vulnerable to HIV/AIDS only because of their partners. Thus we have derived many notable conclusions and suggestion from our study of the socio-economic and psychological aspects of migrant labourers with reference to HIV/AIDS.

## DESCRIPTION OF THE PROBLEM

In view of the linguistic questionnaire we had analyzed the relation among

a. Causes for migrants' vulnerability to HIV/AIDS
b. Factors forcing migration
c. Role of the Government.

We take some subtitles for each of these three main titles.

For the sake of simplicity we are restricted to some major subtitles, which has primarily interested these experts. We use BAM model on the scale [-5, 5]. Here we mention that the analysis can be carried out on any other scale according to the whims and fancies of the investigator.

**A: CAUSES FOR MIGRANT LABOURERS VULNERABILITY TO HIV/AIDS**

$A_1$ - No awareness/ no education
$A_2$ - Social Status
$A_3$ - No social responsibility and social freedom
$A_4$ - Bad company and addictive habits
$A_5$ - Types of profession
$A_6$ - Cheap availability of CSWs.

**F: Factors forcing people for migration**

$F_1$ - Lack of labour opportunities in their hometown
$F_2$ - Poverty/seeking better status of life
$F_3$ - Mobilization of labour contractors
$F_4$ - Infertility of lands due to implementation of wrong research methodologies/failure of monsoon.

**G: Role of the Government**

$G_1$ - Alternate job if the harvest fails there by stopping migration
$G_2$ - Awareness clubs in rural areas about HIV/AIDS

| | | |
|---|---|---|
| $G_3$ | - | Construction of hospitals in rural areas with HIV/AIDS counseling cell/ compulsory HIV/AIDS test before marriage |
| $G_4$ | - | Failed to stop the misled agricultural techniques followed recently by farmers |
| $G_5$ | - | No foresight for the government and no precautionary actions taken from the past occurrences. |

### Experts opinion on the cause for vulnerability to HIV/AIDS and factors for migration

Taking the neuronal field $F_X$ as the attributes connected with the causes of vulnerability resulting in HIV/AIDS and the neuronal field $F_Y$ is taken as factors forcing migration.

The $6 \times 4$ matrix $M_1$ represents the forward synaptic projections from the neuronal field $F_X$ to the neuronal field $F_Y$.

The $4 \times 6$ matrix $M_1^T$ represents the backward synaptic projections $F_X$ to $F_Y$. Now, taking $A_1, A_2, \ldots, A_6$ along the rows and $F_1, \ldots, F_4$ along the columns we get the synaptic connection matrix $M_1$ which is modeled on the scale [–5, 5]

$$M_1 = \begin{bmatrix} 5 & 2 & 4 & 4 \\ 4 & 3 & 5 & 3 \\ -1 & -2 & 4 & 0 \\ 0 & 4 & 2 & 0 \\ 2 & 4 & 3 & 3 \\ 0 & 1 & 2 & 0 \end{bmatrix}$$

Let $X_K$ be the input vector given as (3 4 –1 –3 –2 1) at the $K^{th}$ time period. The initial vector is given such that illiteracy, lack of awareness, social status and cheap availability of CSWs have stronger impact over migration. We suppose that all neuronal state change decisions are synchronous.

The binary signal vector

$$S(X_K) = (1\ 1\ 0\ 0\ 0\ 1).$$

From the activation equation

$$S(X_K)M_1 = (9, 6, 11, 7) = Y_{K+1}.$$

From the activation equation

$$S(Y_{K+1}) = (1\ 1\ 1\ 1).$$

Now

$$S(Y_{K+1})M_1^T = (15, 15, 1, 6, 12, 3) = X_{K+2}.$$

From the activation equation,

$$S(X_{K+2}) = (1\ 1\ 1\ 1\ 1\ 1),$$
$$S(X_{K+2})M_1 = (10, 12, 20, 10) = Y_{K+3}.$$
$$S(Y_{K+3}) = (1\ 1\ 1\ 1).$$

Thus the binary pair {(1 1 1 1 1 1), (1 1 1 1)} represents a fixed point of the dynamical system. Equilibrium of the system has occurred at the time K + 2, when the starting time was K. Thus this fixed point suggests that illiteracy with unawareness, social status and cheap availability of CSW lead to the patients remaining or becoming socially free with no social responsibility, having all addictive habits and bad company which directly depends on the types of profession they choose.

On the other hand, all the factors of migration also come to on state. Suppose we take only the on state that the availability of CSWs at very cheap rates is in the on state. Say at the $K^{th}$ time we have

$$P_K = (0\ 0\ 0\ 0\ 0\ 4),$$
$$S(P_K) = (0\ 0\ 0\ 0\ 0\ 1),$$
$$S(P_K)M_1 = (0\ 1\ 2\ 0)$$

$$= Q_{K+1}.$$

$$S(Q_{K+1}) = (0\ 1\ 1\ 0)$$
$$S(Q_{K+1})M_1^T = (6\ 8\ 2\ 6\ 7\ 3)$$
$$= P_{K+2}$$

$$S(P_{K+2}) = (1\ 1\ 1\ 1\ 1\ 1)$$
$$S(P_{K+2})M_1 = (10\ 12\ 20\ 10)$$
$$= P_{K+3}$$

$$S(P_{K+3}) = (1\ 1\ 1\ 1).$$

Thus the binary pair {(1 1 1 1 1 1), (1 1 1 1)} represents a fixed point. Thus in the dynamical system given by the expert even if only the cheap availability of the CSW is in the on state, all the other states become on, i.e. they are unaware of the disease, their type of profession, they have bad company and addictive habits, they have no social responsibility and no social fear.

Thus one of the major causes for the spread of HIV/AIDS is the cheap availability of CSWs which is mathematically confirmed from our study. Several other states of vectors have been worked by us for deriving the conclusions.

Taking the neuronal field $F_X$ as the role of Government and the neuronal field $F_Y$ as the attributes connected with the causes of vulnerabilities resulting in HIV/ AIDS.

The 6 × 5 matrix $M_2$ represents the forward synaptic projections from the neuronal field $F_X$ to the neuronal field $F_Y$.

The 5 × 6 matrix $M_2^T$ represents the backward synaptic projection $F_X$ to $F_Y$. Now, taking $G_1, G_2, \ldots, G_5$ along the rows and $A_1, A_2, \ldots, A_6$ along the columns we get the synaptic connection matrix $M_2$ in the scale [–5,5] is as follows:

$$M_2 = \begin{bmatrix} 3 & 4 & -2 & 0 & -1 & 5 \\ 5 & 4 & 3 & -1 & 0 & 4 \\ 1 & 3 & 0 & 1 & 4 & 2 \\ 2 & 3 & -2 & -3 & 0 & 3 \\ 3 & 2 & 0 & 3 & 1 & 4 \end{bmatrix}$$

Let $X_K$ be the input vector (-3 4 -2 -1 3) at the $K^{th}$ instant. The initial vector is given such that Awareness clubs in the rural villages and the Government's inability in foreseeing the conflicts have a stronger impact over the vulnerability of HIV/AIDS. We suppose that all neuronal state change decisions are synchronous.

The binary signal vector

$$S(X_K) = (0\ 1\ 0\ 0\ 1).$$

From the activation equation

$$S(X_K)M_2 = (8\ 6\ 3\ 2\ 1\ 8) = Y_{K+1}.$$

Now,

$$S(Y_{K+1}) = (1\ 1\ 1\ 1\ 1),$$
$$S(Y_{K+1})\ M_2^T = (9\ 15\ 11\ 3\ 13) = X_{K+2}.$$

Now,

$$S(X_{K+2}) = (1\ 1\ 1\ 1\ 1),$$
$$S(X_{K+2})M_2 = (14\ 8\ -1\ 0\ 4\ 18) = Y_{K+3}.$$

Now,

$$S(Y_{K+3}) = (1\ 1\ 0\ 1\ 1\ 1),$$
$$S(Y_{K+3})\ M_2^T = (11,\ 12,\ 11,\ 5,\ 13) = X_{K+4}.$$

Thus

$$S(X_{K+4}) = (1\ 1\ 1\ 1\ 1) = X_{K+2}$$

and

$$S(Y_{K+5}) = Y_{K+3}.$$

The binary pair ((1 1 1 1 1), (1 1 0 1 1 1)) represents a fixed point of the BAM dynamical system. Equilibrium for the system occurs at the time K + 4, when the starting time was K.

This fixed point reveals that the other three conditions cannot be ignored and have its consequences in spreading HIV/ AIDS. Similarly by taking a vector $Y_K$ one can derive conclusions based upon the nature of $Y_K$.

We have given this model mainly for illustration. A C-program is given in appendix 6 of this book which makes the calculations simple.

## 1.4 Introduction to Fuzzy Associative Memories

In this section we just recall the FAM structure given by [62]. A fuzzy set is a map $\mu : X \to [0, 1]$ where X is any set called the domain and [0, 1] the range i.e., $\mu$ is thought of as a membership function i.e., to every element $x \in X$ $\mu$ assigns membership value in the interval [0, 1]. But very few try to visualize the geometry of fuzzy sets. It is not only of interest but is meaningful to see the geometry of fuzzy sets when we discuss fuzziness. Till date researchers over looked such visualization [60-69], instead they have interpreted fuzzy sets as generalized indicator or membership functions mappings $\mu$ from domain X to range [0, 1]. But functions are hard to visualize. Fuzzy theorist often picture membership functions as two-dimensional graphs with the domain X represented as a one-dimensional axis.

The geometry of fuzzy sets involves both domain X = $(x_1,\ldots, x_n)$ and the range [0, 1] of mappings $\mu : X \to [0, 1]$. The geometry of fuzzy sets aids us when we describe fuzziness, define fuzzy concepts and prove fuzzy theorems. Visualizing this geometry may by itself provide the most powerful argument for fuzziness.

An odd question reveals the geometry of fuzzy sets. What does the fuzzy power set $F(2^X)$, the set of all fuzzy subsets of X, look like? It looks like a cube, What does a fuzzy set look like? A fuzzy subsets equals the unit hyper cube $I^n = [0, 1]^n$. The fuzzy set is a point in the cube $I^n$. Vertices of the cube $I^n$ define a non-fuzzy set. Now with in the unit hyper cube $I^n = [0, 1]^n$ we are interested in a distance between points, which led to measures of size and fuzziness of a fuzzy set and more fundamentally to a measure. Thus within cube theory directly extends to the continuous case when the space X is a subset of $R^n$. The next step is to consider mappings between fuzzy cubes. This level of abstraction provides a surprising and fruitful alternative to the prepositional and predicate calculus reasoning techniques used in artificial intelligence (AI) expert systems. It allows us to reason with sets instead of propositions. The fuzzy set framework is numerical and multidimensional. The AI framework is symbolic and is one dimensional with usually only bivalent expert rules or propositions allowed. Both frameworks can encode structured knowledge in linguistic form. But the fuzzy approach translates the structured knowledge into a flexible numerical framework and processes it in a manner that resembles neural network processing. The numerical framework also allows us to adaptively infer and modify fuzzy systems perhaps with neural or statistical techniques directly from problem domain sample data.

Between cube theory is fuzzy-systems theory. A fuzzy set defines a point in a cube. A fuzzy system defines a mapping between cubes. A fuzzy system S maps fuzzy sets to fuzzy sets. Thus a fuzzy system S is a transformation S: $I^n \to I^P$. The n-dimensional unit hyper cube $I^n$ houses all the fuzzy subsets of the domain space or input universe of discourse X = $\{x_1, \ldots, x_n\}$. $I^p$ houses all the fuzzy subsets of the range space or output universe of discourse, Y = $\{y_1, \ldots, y_p\}$. X and Y can also denote subsets of $R^n$ and $R^p$. Then the fuzzy power sets F $(2^X)$ and F $(2^Y)$ replace $I^n$ and $I^p$.

In general a fuzzy system S maps families of fuzzy sets to families of fuzzy sets thus S: $I^{n_1} \times \ldots \times I^{n_r} \to I^{p_1} \times \ldots \times I^{p_s}$ Here too we can extend the definition of a fuzzy system to allow arbitrary products or arbitrary mathematical spaces to serve as the domain or range spaces of the fuzzy sets. We shall focus on fuzzy systems S: $I^n \to I^P$ that map balls of fuzzy sets in $I^n$ to balls of fuzzy set in $I^p$. These continuous fuzzy systems behave as associative memories. The map close inputs to close outputs. We shall refer to them as Fuzzy Associative Maps or FAMs.

The simplest FAM encodes the FAM rule or association ($A_i$, $B_i$), which associates the p-dimensional fuzzy set $B_i$ with the n-dimensional fuzzy set $A_i$. These minimal FAMs essentially map one ball in $I^n$ to one ball in $I^p$. They are comparable to simple neural networks. But we need not adaptively train the minimal FAMs. As discussed below, we can directly encode structured knowledge of the form, "If traffic is heavy in this direction then keep the stop light green longer" is a Hebbian-style FAM correlation matrix. In practice we sidestep this large numerical matrix with a virtual representation scheme. In the place of the matrix the user encodes the fuzzy set association (Heavy, longer) as a single linguistic entry in a FAM bank linguistic matrix. In general a FAM system F: $I^n \to I^b$ encodes the processes in parallel a FAM bank of m FAM rules ($A_1$, $B_1$), …, ($A_m$ $B_m$). Each input A to the FAM system activates each stored FAM rule to different degree. The minimal FAM that stores ($A_i$, $B_i$) maps input A to $B_i$' a partly activated version of $B_i$. The more A resembles $A_i$, the more $B_i$' resembles $B_i$. The corresponding output fuzzy set B combines these partially activated fuzzy sets $B_1^1, B_2^1, \ldots, B_m^1$. B equals a weighted average of the partially activated sets B = $w_1 B_1^1 + \ldots + w_n\ B_m^{1`}$ where $w_i$ reflects the credibility frequency or strength of fuzzy association ($A_i$, $B_i$). In practice we usually defuzzify the output waveform B to a single numerical value $y_j$ in Y by computing the fuzzy centroid of B with respect to the output universe of discourse Y.

More generally a FAM system encodes a bank of compound FAM rules that associate multiple output or consequent fuzzy sets $B_i$, …, $B_i^s$ with multiple input or antecedent fuzzy sets $A_i^1$, …, $A_i^r$. We can treat compound FAM rules as compound linguistic conditionals. This allows us to naturally and in many cases easily to obtain structural knowledge. We combine antecedent and consequent sets with logical conjunction, disjunction or negation. For instance, we could interpret the compound association ($A^1$, $A^2$, B), linguistically as the compound conditional "IF $X^1$ is $A^1$ AND $X^2$ is $A^2$, THEN Y is B" if the comma is the fuzzy association ($A^1$, $A^2$, B) denotes conjunction instead of, say, disjunction.

We specify in advance the numerical universe of discourse for fuzzy variables $X^1$, $X^2$ and Y. For each universe of discourse or fuzzy variable X, we specify an appropriate library of fuzzy set values $A_1^2$, …, $A_k^2$. Contiguous fuzzy sets in a library overlap. In principle a neural network can estimate these libraries of fuzzy sets. In practice this is usually unnecessary. The library sets represent a weighted though overlapping quantization of the input space X. They represent the fuzzy set values assumed by a fuzzy variable. A different library of fuzzy sets similarly quantizes the output space Y. Once we define the library of fuzzy sets we construct the FAM by choosing appropriate combinations of input and output fuzzy sets Adaptive techniques can make, assist or modify these choices.

An Adaptive FAM (AFAM) is a time varying FAM system. System parameters gradually change as the FAM system samples and processes data. Here we discuss how natural network algorithms can adaptively infer FAM rules from training data. In principle, learning can modify other FAM system components, such as the libraries of fuzzy sets or the FAM-rule weights $w_i$.

In the following section we propose and illustrate an unsupervised adaptive clustering scheme based on competitive learning to blindly generate and refine the bank of FAM rules. In some cases we can use supervised learning techniques if we have additional information to accurately generate error estimates. Thus Fuzzy Associative Memories (FAMs) are transformation. FAMs map fuzzy sets to fuzzy sets. They map unit cubes to unit cubes. In simplest case the FAM system consists of a single association. In general the FAM system consists of a bank of different FAM association. Each association corresponds to a different numerical FAM matrix or a different entry in a linguistic FAM-bank matrix. We do not combine these matrices as we combine or superimpose neural-network associative memory matrices. We store the matrices and access them in parallel. We begin with single association FAMs. We proceed on to adopt this model to the problem in chapter two.

**Chapter Two**

# ON A NEW CLASS OF n-ADAPTIVE FUZZY MODELS WITH ILLUSTRATIONS

In this chapter for the first time we introduce the new class of n-adaptive fuzzy models (n a positive integer and n ≥ 2) and illustrate it. These n-adaptive fuzzy models can analyze a problem using different models so one can get in one case the hidden pattern, in one case the maximum time period in some case output state vector for a given input vector and so on. So this new model has the capacity to analyze the problem in different angles, which can give multiple suggestions and solutions about the problem.

This chapter has two sections. Section one defines the new model gives illustrations and section two defines some special n-adaptive models and proposes some problems.

## 2.1 On a new class of n-adaptive fuzzy models with illustrations

In this section we proceed on to define the new class of n-adaptive fuzzy models and illustrate them with examples.

**DEFINITION 2.1.1:** *Let $C_n = \{M_1, M_2, \ldots, M_n\}$ be the collection of n distinct fuzzy models $M_1, \ldots, M_n$ ($n \geq 2$). We say the $C_n$ is a new class of n-adaptive fuzzy model if from any model $M_j$ we can go to another fuzzy model $M_i$ ($i \neq j$) by some simple transformation function or at time even with identity transformation (Only operations in Models $M_i$ and $M_j$ will be* *different). This is true for $1 \leq j \leq n$. Thus one can go from any $M_r$ to $M_s$, $r \neq s$, $1 \leq r, s \leq n$.*

It is interesting and important to make note of the following:

***Remark:*** If $C_n = \{M_1, M_2, \ldots, M_n\}$ be a new class of n-adaptive fuzzy models then for any $2 \leq m \leq n$ we have $C_m$ to be a new class of m-adaptive fuzzy model. Thus a 5-adaptive fuzzy model gives one a 4-adaptive fuzzy model, 3-adaptive fuzzy model and 2-adaptive fuzzy model. Thus given a n-adaptive fuzzy model we have (n – 2) distinct adaptive fuzzy models. Thus this helps us in reducing the number of models as per need. We in this section give a new class of 2-adaptive fuzzy model and illustrate it with a real world problem. We first describe the models in the new class of 2-adaptive fuzzy models, which we will be using in chapter III, the adaptiveness also is very much dependent on the problem under study. One cannot always say for all problems same $C_n$ will serve the purpose.

***Example 2.1.1:*** Let $C_2 = \{M_1, M_2\}$ where

$M_1$ is the fuzzy matrix model described in chapter 1 section 1.1 of the book.

$M_2$ is the Birectional Associative Memories (BAM) model described in section 1.3 of chapter 1.

Now we will describe how we go from one model to the other so that the pair in $C_2$ happens to be a 2-adaptive fuzzy model. Now the model $M_1$, the fuzzy matrix model is a five stage analysis we have the Initial Raw data matrix, ATD matrix, RTD matrix, CETD matrix and the row sum and $M_2$ the BAM model.

Now to convert model $M_1$ i.e., the fuzzy matrices model into a BAM model $M_2$. This is done by identical transformation. For any CETD matrix given in the model $M_1$ will have its entries from the interval [–t, t], t a positive integer that denotes directly the number of parameters used in forming the CETD matrix.

If t is the number of parameters used in forming the CETD matrix. The interval [–t, t] will form the scale of study for the BAM model. Now the CETD matrix itself can be taken as the associated matrix of BAM model $M_2$ in the interval [–t, t].

Thus the transformation in this case is the identical transformation. Thus we have transformed from the fuzzy matrix model to the BAM model.

Now we give a procedure of going from the model $M_2$ i.e., the BAM model to the model $M_1$ i.e., the fuzzy matrix model.

Now any BAM model will have its entries in the interval [–t, t].

Now if $P = (p_{ij})$ is any matrix associated with the BAM model its entries will be from the interval [–t, t]. Now divide each and every entry of P by t, so

$$\frac{P}{t} = \left(\frac{p_{ij}}{t}\right) = A = \left(a_{ij} = \frac{p_{ij}}{t}\right).$$

Now $A = (a_{ij})$ will be taken as the ATD matrix of the fuzzy matrix model $M_1$.

Using this ATD matrix we can construct the RTD matrix and the CETD matrix. Thus a BAM model $M_2$ can be easily converted into a fuzzy matrix model $M_1$. The conversion strongly depends on the expert and above all on the problem under investigation. Some problems may not yield to these transformation.

Now we illustrate this by the following example from the real world problem.

***Example 2.1.2:*** We take the 2-adaptive fuzzy model with the 2 distinct fuzzy models viz. {fuzzy matrices model ($M_1$) and BAM model ($M_2$)}.

The problem given in page 11 which uses RTD and CETD is already given in chapter 1 is taken. Now the 6 × 6 CETD matrix M page 15 is as follows:

$$M = \begin{bmatrix} -3 & -1 & -3 & -3 & -2 & -2 \\ 2 & 3 & 2 & 2 & 3 & 3 \\ 3 & 3 & 3 & 3 & 3 & 3 \\ 1 & 0 & 1 & 0 & -1 & -1 \\ 0 & -2 & -1 & -2 & -2 & -2 \\ -3 & -3 & -3 & -3 & -3 & -3 \end{bmatrix}$$

Now take M as the BAM ($M_2$) model. Now using the C-program given in appendix 5 one can work with the model and arrive at the conclusions.

We now proceed on to give yet another 2-adaptive model. $C'_2[M'_1, M'_2]$ be a 2 adaptive fuzzy model where $M_1$ is a FRM model and $M_2$ is a BAM model.

***Example 2.1.3:*** Take $C'_2[M'_1, M'_2]$ where $M'_1$ is the FRM model and $M'_2$ is the BAM model. Now we give the transformation of the FRM model $M_1$ to the BAM model $M_2$. Now first we decide on what scale the BAM model $M_2$ must be (say p). Accordingly we take p experts opinion and form the p combined FRM matrices say $P_1, P_2, \ldots, P_p$. Now $P = P_1 + \ldots + P_p$ the CFRM (matrix addition) is taken as the BAM model on the scale $[-p, p]$. Now we can work with the BAM model. Thus we transform the FRM model $M_1$ to the BAM model $M_2$ by the following function.

Now we give the transformation of the model $M_2$ (i.e., BAM model into the FRM model $M_1$. Suppose the BAM model $M_2$ is represented by the $M = (m_{ij})$ matrix on the scale $[-p\ p]$ then we divide

$$\frac{M}{p} = \left(\frac{m_{ij}}{p}\right)$$

now use a function f on

$$M = \left(\frac{m_{ij}}{p}\right) \text{ by defining } f\left(\frac{m_{ij}}{p}\right) \leq c - \in$$

then replace

$$\left(\frac{m_{ij}}{p}\right) \text{ by } -1 \text{ if } \left(\frac{m_{ij}}{p}\right) \geq c + \in$$

replace by +1 if

$$\left(\frac{m_{ij}}{p}\right) \in (c - \in,\ c + \in) \text{ put it as } 0,$$

($0 < \in < c$ is chosen, the choice of c and $\in$ depend on the expert and the problem under analysis).

The choice of the function is in the hands of the expert by studying the problem at hand. Now the entries of the new matrix

$$f(M) = \left\{ f\left(\frac{m_{ij}}{p}\right) \right\}$$

is a FRM matrix with entries from the set {–1, 0, 1}. Thus we have transformed the BAM model $M_2$ into the FRM model $M_1$. Now using the C-program given in appendices 4, 5 and 6 we can arrive at the solution.

We give yet another example of a 2-adaptive fuzzy model.

***Example 2.1.4:*** Let $C''_2$ [$M_1$, $M_2$] be a 2-adaptive fuzzy model where $M_1$ is a FRM model described in section 1.2 and $M_2$ be a FAM model described in section 1.4.

Let $W_1, W_2, \ldots, W_7$ be 7 attributes related to women given by expert as the nodes of the domain space and $R_1, R_2, \ldots, R_{10}$ be the nodes of the range space which gives the vulnerability of women becoming HIV/AIDS patients.

Now using the experts opinion we give the FRM directed graph, which is given in the following.

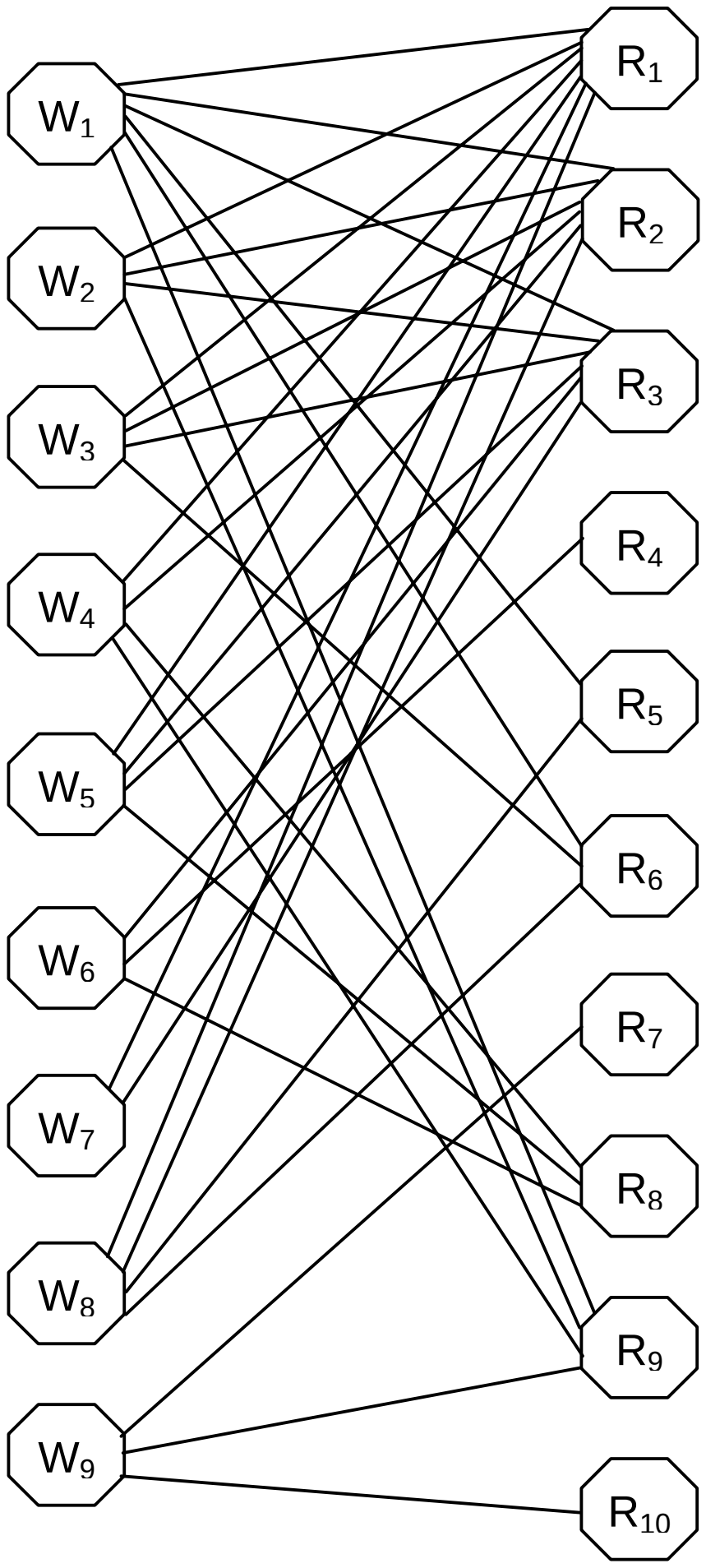

**FIGURE 1**

Using this directed graph we obtain the related connection matrix T which is a 9 × 10 matrix.

$$T = \begin{array}{c} \\ W_1 \\ W_2 \\ W_3 \\ W_4 \\ W_5 \\ W_6 \\ W_7 \\ W_8 \\ W_9 \end{array} \begin{array}{c} \begin{array}{cccccccccc} R_1 & R_2 & R_3 & R_4 & R_5 & R_6 & R_7 & R_8 & R_9 & R_{10} \end{array} \\ \left[\begin{array}{cccccccccc} 1 & 1 & 1 & 0 & 1 & 1 & 0 & 0 & 1 & 0 \\ 1 & 1 & 1 & 0 & 0 & 0 & 0 & 0 & 1 & 0 \\ 1 & 1 & 1 & 0 & 0 & 1 & 0 & 0 & 0 & 0 \\ 1 & 1 & 0 & 0 & 0 & 0 & 1 & 1 & 0 & 0 \\ 1 & 1 & 1 & 0 & 0 & 0 & 1 & 0 & 0 & 0 \\ 0 & 0 & 1 & 1 & 0 & 0 & 0 & 1 & 0 & 0 \\ 1 & 0 & 1 & 0 & 0 & 0 & 0 & 0 & 0 & 0 \\ 1 & 1 & 0 & 0 & 1 & 1 & 0 & 0 & 0 & 0 \\ 0 & 0 & 0 & 0 & 0 & 0 & 1 & 0 & 1 & 1 \end{array}\right] \end{array}$$

Now T is the dynamical system of our FRM model. We can study / analyze the effect of any state vector either from the domain space D or from the range space R.

Suppose one is interested to study the effect of $W_9$ on the dynamical system i.e. the rural poor uneducated women have no voice in the traditional set up alone is in the on state and all other nodes are in the off state to study the effect of this state vector. Let $X_1 = (0\ 0\ 0\ 0\ 0\ 0\ 0\ 0\ 1)$, i.e., $W_9$ alone is in the on state and all other nodes are in the off state. Its impact on the dynamical system T.

$$\begin{aligned} X_1T &= (0\,0\,0\,0\,0\,0\,1\,0\,1\,1) \\ Y_1T^{\,t} &= Y_1 \\ &\to (1\ 1\ 0\ 1\ 0\ 0\ 0\ 0\ 1) \\ &= X_2 \text{ (say)} \\ X_2T &= (3\ 3\ 2\ 0\ 1\ 1\ 1\ 1\ 4\ 1) \end{aligned}$$

after thresholding

$$Y_2 = (1\ 1\ 1\ 0\ 1\ 1\ 1\ 1\ 1\ 1\,).$$

$$Y_2T^t = (6\ 4\ 4\ 4\ 4\ 2\ 2\ 4\ 3)$$

after thresholding and updating we get

$$X_3 = (1\ 1\ 1\ 1\ 1\ 1\ 1\ 1\ 1).$$

Now we study the effect of $X_3$ on the dynamical system T.

$$X_3T = (7\ 6\ 6\ 1\ 2\ 3\ 1\ 3\ 4\ 1).$$

Now thresholding and updating the resultant vector

$$Y_3 = (1\ 1\ 1\ 1\ 1\ 1\ 1\ 1\ 1\ 1).$$

Thus the fixed point of the dynamical system is a pair {(1 1 1 1 1 1 1 1 1), (1 1 1 1 1 1 1 1 1 1)}.

Now we see when a single node "No voice for women basically in the traditional set up" alone is in the on state all the other nodes both in the domain and range space came to an state which shows that for the spread of HIV/AIDS among rural uneducated women in India is basically due to women voicelessness in the society.

They as a part of nations tradition cannot say no to sex to their HIV/AIDS infected husbands, for male denominated society demands it. So unless the tradition which is basically ruling the Indian society is broken there can be no solution to this problem. Poor rural uneducated women would continue to be the worst victims of HIV/AIDS. Thus the spread of HIV/AIDS in the nation above all in the rural India cannot be controlled.

Now we proceed on to study the effect of a state vector from the range space of the dynamical system T.

Let us consider the on state of the single state vector $R_5$ alone on the dynamical system i.e., "the role of media in inducing male ego" and all other nodes are assumed to be in the off state. Let Y = (0 0 0 0 1 0 0 0 0 0). Now we study the effect of Y on the dynamical system $T^t$.

$$YT^t = (1\ 0\ 0\ 0\ 0\ 0\ 0\ 1\ 0) = X.$$

The effect of X on the dynamical system T is given by

$$XT = (2\ 2\ 1\ 0\ 2\ 2\ 0\ 0\ 1\ 0).$$

After thresholding and updating we get the resultant as

$$Y_1 = (1\ 1\ 1\ 0\ 1\ 1\ 0\ 0\ 1\ 0).$$

Now we analyze the effect of $Y_1$ on $T^t$

$$Y_1 T^t = (6\ 4\ 4\ 3\ 3\ 1\ 2\ 3\ 1).$$

After thresholding and updating we get the resultant as

$$X_1 = (1\ 1\ 1\ 1\ 1\ 1\ 1\ 1\ 1).$$

Now the effect of $X_1$ on T is given by

$$X_1T = (7\ 6\ 7\ 1\ 2\ 3\ 1\ 3\ 4\ 1).$$

After updating and thresholding the resultant vector we get

$$Y_2 = (1\ 1\ 1\ 1\ 1\ 1\ 1\ 1\ 1\ 1\ ).$$

Thus we see the hidden pattern of the dynamical system is a fixed point given by the binary pair. {(1 1 1 1 1 1 1 1 1), (1 1 1 1 1 1 1 1 1 1)}. Thus when only the single fact that the role of media in inbreeding the male ego alone is in the on state, all the nodes both in the domain space and in the range space become on where by it is clearly seen that India is emotionally and culturally ruled by the cinema world. That is why the present youth, have mainly have their role models as actor and actresses, fully for getting themselves and the fact that it is only on the screen and not in real life. Unless politicians take up the issue for the women cause and see to it that movies which put down potent and individual women with a mission should be banned nothing concrete can be achieved. People should come forth to take up movies which give true importance to the role of women and their power. Unless such revolutionary steps are taken it may not be possible to put an end to spread of HIV/AIDS in India that too among poor rural in India.

Not only today's movie inbuilt in males their ego and superiority over women but also invariably in all most all of the movies heroes are portrayed with cigar in their hand / or in a pose of drinking alcohol so these bad habits in their real life. Thus the cinema seems to be one of the major cause which inbreeds male ego.

Using this FRM matrix we using a function transform the FRM connection matrix into a fuzzy vector matrix, i.e., a FAM model. At the outset we first wish to state that these women are mainly from the rural areas, they are economically poor and uneducated and majority of them are infected by their husbands. Second we use FAM because only this will indicate the gradations of the causes, which is the major cause for women being affected by HIV/AIDS followed in order of gradation the causes. Further among the fuzzy tools FAM model alone can give such gradations so we use them in this analysis. Another reason for using FAM is they can be used with the same FRM that is using the attributes of the FRM, FAMs can also be formulated. Already FAMs are very briefly described in section one of this chapter. We now give the sketch of the analysis of this problem with a view that any reader with a high school education will be in a position to follow it.

We just illustrate two experts opinion though we have used several experts' opinion for this analysis. Using the problems of women affected with HIV/AIDS along the rows and the causes of it along the column we obtain the related fuzzy vector matrix M. The entries are given by the experts using a transformation.

The gradations are given in the form of the fuzzy vector matrix M which is as follows:

$$\begin{array}{c} \\ W_1 \\ W_2 \\ W_3 \\ W_4 \\ W_5 \\ W_6 \\ W_7 \end{array}\begin{array}{c} \begin{array}{cccccccccc} R_1 & R_2 & R_3 & R_4 & R_5 & R_6 & R_7 & R_8 & R_9 & R_{10} \end{array} \\ \left[\begin{array}{cccccccccc} 0.9 & 0.8 & 0.7 & 0 & 0 & 0 & 0 & 0 & 0 & 0.7 \\ 0.5 & 0.8 & 0.6 & 0 & 0 & 0 & 0 & 0 & 0 & 0 \\ 0 & 0.3 & 0.6 & 0 & 0 & 0 & 0 & 0 & 0 & 0 \\ 0 & 0 & 0 & 0.6 & 0 & 0 & 0 & 0 & 0 & 0 \\ 0 & 0 & 0 & 0 & 0.9 & .6 & .7 & 0 & 0 & 0 \\ 0 & 0 & 0 & 0 & 0 & 0.7 & 0.5 & 0 & 0 & 0 \\ 0 & 0 & 0 & 0 & 0.6 & 0 & 0 & 0 & 0 & 0 \end{array}\right] \end{array}$$

Now using the expert's opinion we get the fit vectors. Suppose the fit vector B is given as B = (0 1 1 0 0 0 0 0 1 0). Using max-min in backward directional we get FAM described as

$$A = M \circ B$$

that is $a_i = \max \min (m_{ij}, b_j) \quad 1 \leq j \leq 10.$

Thus A = (0.8, 0.8, 0.6 0 0 0 0) since 0.8 is the largest value in the fit vector in A and it is associated with the two nodes vide

$W_1$ and $W_2$ child marriage / widower child marriage etc finds it first place also the vulnerability of rural uneducated women find the same states as that of $W_1$. Further the second place is given to $W_3$, disease untreated till it is chronic or they are in last stages, all other states are in the off state.

Suppose we consider the resultant vector

$$A = (0.8, 0.8, 0.6, 0\ 0\ 0\ 0)$$

$$A \text{ o } M = B. \text{ where}$$

$$B_j = \max (a_i, m_{ij}) \qquad 1 \leq j \leq 10.$$

A o M = (0.9, 0.8, 0.8, 0.6, 0.9, 0.7, 0.7, 0.0, 0.7) = $B_1$.

Thus we see the major cause being infected by HIV/AIDS is $R_1$ and $R_5$ having the maximum value. viz. female child a burden so the sooner they get married off the better relief economically; and women suffer no guilt and fear for life. The next value being given by $R_2$ and $R_3$, poverty and bad habits of men are the causes of women becoming HIV/AIDS victims.

The next value being $R_6$, $R_7$ and $R_{10}$ taking the value 0.7, women have not changed religion and developed faith in god after the disease, lost faith in god after the disease. Husbands hide their disease from family so wife becomes HIV/AIDS infected. Several conclusions can be derived from the analysis of similar form.

## 2.2 Some special n-adaptive models

In this section we define the notion of directed or ordered n adaptive model and semi n-adaptive model.

Now we define yet another new class of n-adaptive fuzzy model called directed (ordered) n-adaptive fuzzy model.

**DEFINITION 2.2.1:** *Let $C_n^d$ = {$M_1, \ldots, M_n$} where $M_1, M_2, \ldots, M_n$ denote a class of n distinct fuzzy models we say $C_n^d$ is a new class of directed or ordered n-adaptive fuzzy model if and only if we can have transformation from $M_1$ to $M_2$, $M_2$ to $M_3$, …, $M_{n-1}$ to $M_n$ i.e., we have transformation $M_i \rightarrow M_{i+1}$, i = 1, 2, …, n–1.*

That is given a directed or an ordered n-adaptive fuzzy model we cannot in general have a transformation from $M_r$ to $M_s$ if $s \neq r + 1$.

Now we proceed on to define the notion of semi directed r-adaptive fuzzy model $2 \leq r \leq n$.

**DEFINITION 2.2.2:** *Let $C_n^q = \{M_1, \ldots, M_n\}$ where $M_1, M_2, \ldots, M_n$ are distinct fuzzy models we say $C_n^q$ denote a class of semi directed r-adaptive fuzzy model ($2 \leq r < n$) if there exist r fuzzy models $\left\{M_{i_1}, \ldots, M_{i_r}\right\}$ in $C_n^q$, $1 \leq i_1\ i_2 \ldots i_r \leq n$ such that the set of r fuzzy models $\left\{M_{i_1}, \ldots, M_{i_r}\right\}$ forms a directed r-adaptive fuzzy model; that is we have a transformation i.e., $M_{i_t} \to M_{i_{t+1}}$. Clearly $M_{i_t}$ can be $M_s$ and $M_{i_{t+1}}$ can be $M_p$, we can have $p < s$.*

Now we propose some problems for the interested reader:

1. Give an example of a 4-adaptive fuzzy model.
2. Illustrate a 7 adaptive model using a real world problem.
3. What can be the approximate large number n such that we have a n-adaptive fuzzy model?

It is important to note that from the very definition of the new class of directed n-adaptive fuzzy model we see they are a particular case of the new class of n-adaptive fuzzy model, i.e., Every n-adaptive fuzzy model is a directed n- adaptive fuzzy model but a directed n-adaptive fuzzy model in general is not a n-adaptive fuzzy model.

We propose the following problems

1. Construct a directed 2-adaptive fuzzy model
2. Find an example of a directed 4 adaptive fuzzy model.
3. What is a maximum n such that we have a directed n-adaptive fuzzy model?

**Chapter Three**

# USE OF 2-ADAPTIVE FUZZY MODEL TO ANALYZE THE PUBLIC AWARENESS OF HIV/AIDS

In this chapter we using the 2-adaptive fuzzy model introduced in chapter two, analyze the data collected from the 101 public people living in and around Chennai about HIV/AIDS. The data was collected using a linguistic questionnaire enclosed in appendix 1. A brief table of statistics is given in appendix 2 for ready reference.

The 2-adaptive fuzzy model with the 2 models, fuzzy matrix model $M_1$, and BAM model $M_2$ are used. The proof how the new class, $C_2 = \{M_1, M_2\}$ is a 2-adaptive system is also proved in earlier chapter. Now using the 101 interviews we first use the fuzzy matrix model and convert the raw data into a ATD, a RTD and CETD matrices.

In this chapter we mathematically analyze the data collected from the 101 public people living in and around Chennai about HIV/ AIDS, and the awareness programs. The data pertains to rich and poor, educated and uneducated, labourers, students forming a heterogeneous collection.

This chapter has two sections, in the first section we use the CETD matrix to analyze that data. In section two we use the 2-adaptive fuzzy model to analyze the data. Conclusions are given based on this analysis. We have used C-program given in Appendix 3 and 4 to make easy and quick analysis.

### 3.1 Study of the psychological and social problems the public have about HIV/AIDS patients using CETD matrix

Here we give a brief description of the CETD matrix approach and a complete analysis of the data using CETD matrix.

This raw data collected from the public is converted into time-dependent matrices. By time dependent matrices we mean the matrices, which are dependent on the age of the public. So when we say time in this problem we only mean the age groups. After obtaining the time dependent matrices (age-dependent matrices) using the technique of usual mean/average and standard deviation we identify the age group in which they show aversion, fear, unconcern, social stigma etc. The identification of the age group when they show maximum of the negative attributes, can be used to prepare a well-tailored intervention and awareness programs.

The majority of such standards are based on their education, profession, social and economic status and above all their family background. The raw data under investigation is classified under the twelve broad heads: Averse, Positive / comforting attitude, Hostile, Indifferent, Concerned, Careless, Dislike or show open hatred, Unconcerned, Fear, Social Stigma, Shyness / Conservative and Sympathetic towards the patients, as these were given by many of the experts.

These twelve broad heads forms the columns of the matrices. The time periods are taken as the varying age groups.

The analysis is carried out in five stages. In the first stage we get a matrix representation of the collected raw data. Entries corresponding to the intersection of rows and columns are values corresponding to a live network. This initial $M \times N$ matrix is not uniform i.e. the number of years in each interval of time period may or may not be the same; for instance, one period may be 21-30 and other $\leq 19$. So, in the second stage, in order to obtain an unbiased uniform effect on each and every data so collected, we transform this initial matrix into an Average Time-dependent Data (ATD) matrix. To make the calculations simpler, in the third stage, we use the simple average techniques and convert the above Average Time-dependent Data (ATD) matrix $A = (a_{ij})$ into a matrix with entries $e_{ij}$ where $e_{ij} \in \{-1, 0, 1\}$. We name this matrix as the Refined Time-dependent Data (RTD) matrix. The value of the $e_{ij}$ corresponding to each entry is determined in a special way. At the next stage using the RTD matrix we get the Combined Effect Time Dependent Data matrix (CETD matrix), which gives the cumulative effect of all these entries. In the final stage we obtain the row sums of the CETD matrix. A program in C is written which easily estimates all these five stages. The C program given in appendix 3 shall give details of the mathematical working and method of calculation. Now we give the description of the problem and the proposed solution to the problem.

We had interviewed 101 people all of them public from city. The persons whom we interviewed formed a very heterogeneous crowed varying from school students, teachers, auto diver, vegetable vendor, professor, housewife, people working in clerical posts, cobblers, tailors, and so on. So the age group also is heterogeneous. However we are first analyzing the problem using the Refined Time Dependent Matrix for they are very simple but a powerful tool to say how each age group feels and the results are also represented as simple graphs so that even a lay man can follow and know what the analysis is about. For this first we divide the data age wise into six classes this is taken along the rows of the RTD matrix and the attributes or feelings related with the disease is taken along the column of the RTD matrix.

We have restrained ourselves by taking only 12 attributes / feeling/ concepts which are described very briefly in a line or two.

The twelve attributes given by the experts are given below with a brief description.

$A_1$ – Averse – People felt very much aversed to discuss about the disease. Some made their aversion only towards the HIV/AIDS patients in particular. As we cannot divide and analyse we have under the head of Averse include both.

$A_2$ – Positive / comforting – In their discussions they show positive feelings or comforting words or expressions towards the patients or discuss positively about the disease.

$A_3$ – Hostile – they may show Hostileness towards the disease or towards the infected persons or towards the persons who discuss about the disease with them.

$A_4$ – Indifferent – They do not react or they are so indifferent towards the disease or towards the infected persons or towards the persons who discuss about the disease. This is gathered from our discussions with them.

$A_5$ – Concerned – When we say concerned, it is they talk with interest and involvement about the help to be rendered to women patients or children or orphaned children and so on.

$A_6$ – Careless – While interviewing / discussing about the disease or the awareness program just they are very causal, talk with no involvement, talk in a way to hurt the patients / infected persons.

Some made such statements one such is why should we bother about HIV/AIDS infected children? It is none of our concern. They reap the acts of their sin and so on.

$A_7$ – Dislike – Some showed positive dislike to discuss or answer our questionnaire. Never even wished to have a look at the questionnaire. Showed positive dislike over HIV/AIDS disease / HIV/AIDS infected persons / the awareness program about HIV/AIDS.

$A_8$ – Unconcerned – Some were totally unconcerned when we talked about persons infected with HIV/AIDS and their rights and how they should be protected. Some said it is because of Karma they suffer the disease so on.

$A_9$ – Fear – Almost all showed fear in varying degrees. They were utterly frightened about HIV/AIDS one of the reasons may be it has no cure.

$A_{10}$ – Social Stigma – The social stigma they attached with the disease and towards the infected knew no bounds. Infact due to our traditional social set-up they were unable to separate the stigma from the disease or vice versa.

$A_{11}$ – Shyness / Conservative: They were shy / conservative to openly discuss about the disease / awareness program / about rehabilitating the children orphaned / infected with the disease because they attached wrong sex with the disease.

$A_{12}$ – Sympathetic – However some of them were sympathetic towards the children infected with the disease some with women infected by their husbands only one out of the 101 showed some sympathy towards his close male relative apart from that only a few school children and rag pickers showed some sympathy. Sympathy towards the infected men was almost absent. Infact more hatred alone was projected towards them.

Now we have picked up these 12 concepts for most of the experts with whom we discussed gave these. More powerful negative approach was given by very few which we felt is not essential of even mention.

Now as said columns will be occupied by the attributes $A_1$, $A_2$, …, $A_{12}$ and the age groups will be taken along the rows. We have taken a non-over lapping age group > 60 years, 50-59 year, 40-49 years, 30-39 years, 20-29 years and ≤ 19 years. We have divided into 6 classes however we can still refine it on experts recommendations.

Now we find number of persons who are averse ($A_1$) with HIV/AIDS disease or patients or the awareness program, thus we place them under head $A_1$ and we say a single person can find his membership under several heads. For he can be averse, fear the disease, suffer social stigma and so on. This is put in the form of a table with 6 rows and 12 columns. This table can be just represented as a matrix, which is known as the raw time dependent data matrix. Now we give the table of the number of persons in each age group.

| **Age Group** | **No of persons** |
|---|---|
| ≥ 60 | 5 |
| 50-59 | 13 |
| 40-49 | 16 |
| 30-39 | 17 |
| 20-29 | 36 |
| ≤19 | 14 |
| Total | 101 |

We see the youngsters in the age group 20 to 29 were more willing to give interview. However the persons in the age group 40-49 and 30-39 showed equal interest. Older people gave interviews but were in most of the cases irrelevant and not upto the mark to be taken up for they boasted of their times and their good living utterly blaming the youth or blaming the government or the education system and so on. Thus though we interviewed over 130 persons only 101 were atleast upto some merit to be considered.

| | $A_1$ | $A_2$ | $A_3$ | $A_4$ | $A_5$ | $A_6$ | $A_7$ | $A_8$ | $A_9$ | $A_{10}$ | $A_{11}$ | $A_{12}$ |
|---|---|---|---|---|---|---|---|---|---|---|---|---|
| ≥ 60 | 5 | 0 | 4 | 0 | 2 | 0 | 1 | 2 | 5 | 5 | 4 | 2 |
| 50-59 | 6 | 4 | 0 | 2 | 8 | 1 | 2 | 3 | 8 | 10 | 6 | 8 |
| 40-49 | 9 | 2 | 7 | 5 | 6 | 1 | 5 | 7 | 15 | 15 | 9 | 7 |
| 30-39 | 7 | 0 | 5 | 2 | 9 | 3 | 4 | 5 | 17 | 17 | 16 | 7 |
| 20-29 | 15 | 5 | 13 | 4 | 21 | 0 | 8 | 8 | 36 | 35 | 25 | 17 |
| ≤ 19 | 7 | 4 | 1 | 1 | 9 | 1 | 1 | 5 | 11 | 14 | 11 | 10 |

A single person can opt for more than one attribute infact in our analysis they have given minimum four to five attributes.

The Raw time dependent data matrix is given in the following

$$\begin{pmatrix} 5 & 0 & 4 & 0 & 2 & 0 & 1 & 2 & 5 & 5 & 4 & 2 \\ 6 & 4 & 0 & 2 & 8 & 1 & 2 & 3 & 8 & 10 & 6 & 8 \\ 9 & 2 & 7 & 5 & 6 & 1 & 5 & 7 & 15 & 15 & 9 & 7 \\ 7 & 0 & 5 & 2 & 9 & 3 & 4 & 5 & 17 & 17 & 16 & 7 \\ 15 & 5 & 13 & 4 & 21 & 0 & 8 & 8 & 36 & 35 & 25 & 17 \\ 7 & 4 & 1 & 1 & 9 & 1 & 1 & 5 & 11 & 14 & 11 & 10 \end{pmatrix}.$$

We now transform the raw data matrix into Average Time Dependent data (ATD) matrix; A = ($a_{ij}$).

$$\begin{array}{cccccccccccc} A_1 & A_2 & A_3 & A_4 & A_5 & A_6 & A_7 & A_8 & A_9 & A_{10} & A_{11} & A_{12} \end{array}$$
$$\begin{pmatrix} 1 & 0 & 0.8 & 0 & 0.4 & 0 & 0.2 & 0.4 & 1 & 1 & 0.8 & 0.4 \\ .46 & .3 & 0 & .15 & .62 & .08 & .15 & .23 & .62 & .77 & .46 & .62 \\ .56 & .13 & .44 & .31 & .38 & .06 & .31 & .44 & .94 & .94 & .56 & .44 \\ .41 & 0 & .29 & .12 & .53 & .18 & .24 & .29 & 1 & 1 & .94 & .41 \\ .42 & .14 & .36 & .11 & .58 & 0 & .22 & .22 & 1 & .97 & .69 & .42 \\ .5 & .29 & .07 & .07 & .64 & .07 & .07 & .36 & .79 & 1 & .79 & .71 \end{pmatrix}$$

Average and Standard Deviation of the columns of the ATD matrix is given below:

| | | | | | | | | | | | | |
|---|---|---|---|---|---|---|---|---|---|---|---|---|
| A | .56 | .14 | .33 | .13 | .53 | .07 | .2 | .32 | .89 | .95 | .71 | .5 |
| SD | .2 | .11 | .26 | .09 | .01 | .06 | .09 | .08 | .14 | .08 | .14 | .14 |

The Refined Time Dependent Data Matrix (RTD matrix) is got using the formula

$a_{ij} \leq (\mu_j - \alpha * \sigma_j)$, $e_{ij} = -1$
$a_{ij} \in (\mu_j - \alpha . \sigma_j, \mu_j + \alpha . \sigma_j)$, $e_{ij} = 0$
$a_{ij} \geq (\mu_j + \alpha . \sigma_j)$, $e_{ij} = 1$, for $\alpha \in [0\ 1]$.

The RTD matrix with respect to $\alpha = 1$ is as follows:

$$\begin{matrix} A_1 & A_2 & A_3 & A_4 & A_5 & A_6 & A_7 & A_8 & A_9 & A_{10} & A_{11} & A_{12} \end{matrix}$$
$$\begin{pmatrix} 1 & -1 & 1 & -1 & -1 & -1 & 0 & 1 & 0 & 0 & 0 & 0 \\ 0 & 1 & -1 & 0 & 1 & 0 & 0 & -1 & -1 & -1 & -1 & 0 \\ 0 & 0 & 0 & 1 & -1 & 0 & 1 & 1 & 0 & 0 & -1 & 0 \\ 0 & -1 & 0 & 0 & 0 & 1 & 0 & 0 & 0 & 0 & 1 & 0 \\ 0 & 0 & 0 & 0 & 1 & -1 & 0 & 0 & 0 & 0 & 0 & 0 \\ 0 & 1 & -1 & 0 & 1 & 0 & -1 & 0 & 0 & 0 & 0 & 1 \end{pmatrix}$$

Row sum of the RTD matrix

$$\begin{matrix} \geq 60 \\ 50-59 \\ 40-49 \\ 30-39 \\ 20-29 \\ \leq 19 \end{matrix} \begin{bmatrix} -1 \\ -3 \\ 1 \\ 1 \\ 0 \\ 1 \end{bmatrix}.$$

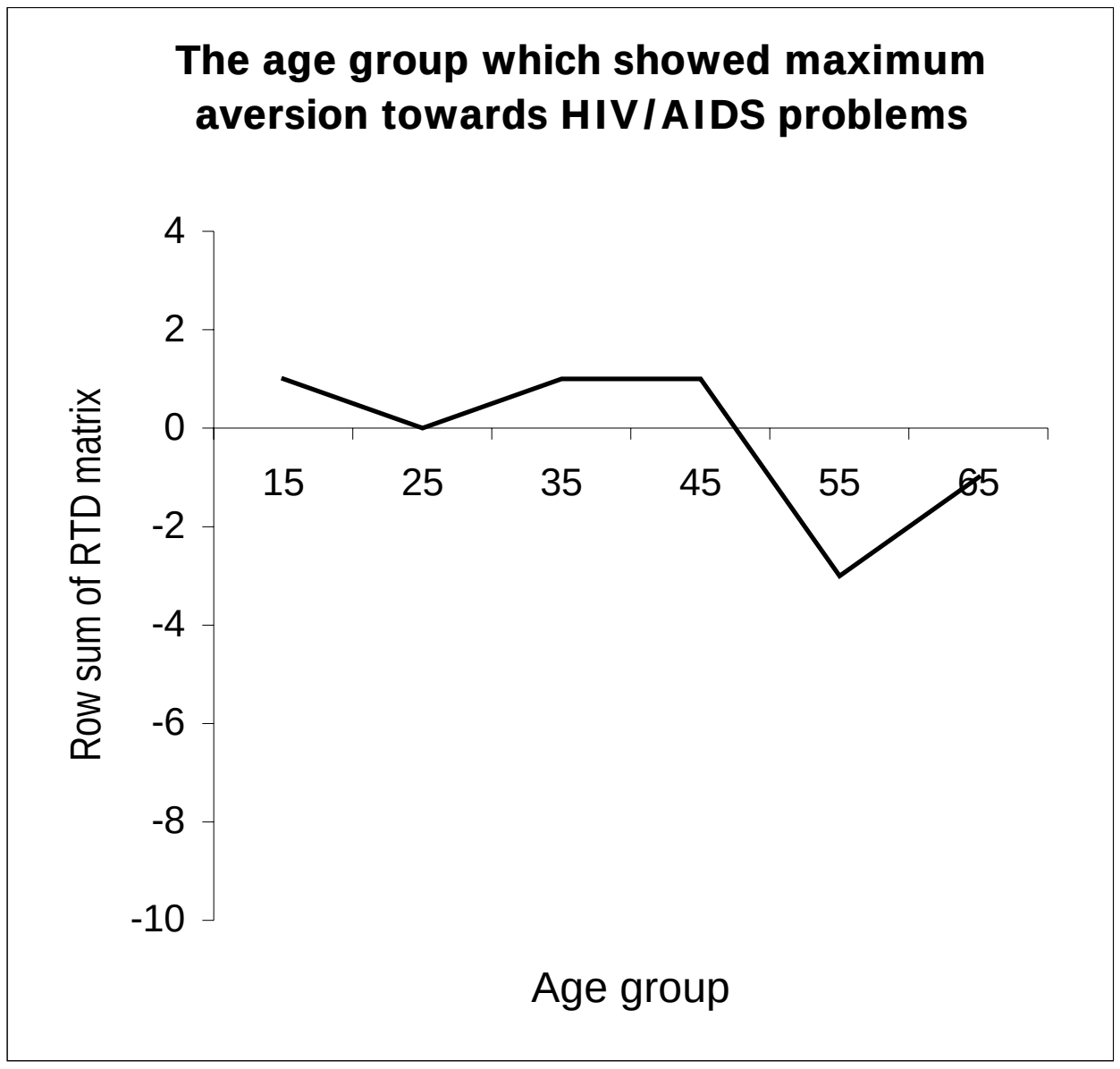

It is very clear from the study when $\alpha = 1$ the age groups 40-49, 30 – 39 and ≤ 19 behave or approach in the same way.

The most repulsive age group being 50 to 59 and 20 to 29 is passive.

Now take $\alpha = 0.7$ we get the following RTD matrix.

$$\begin{array}{cccccccccccc} A_1 & A_2 & A_3 & A_4 & A_5 & A_6 & A_7 & A_8 & A_9 & A_{10} & A_{11} & A_{12} \end{array}$$
$$\begin{pmatrix} 1 & -1 & 1 & -1 & -1 & -1 & 0 & 1 & 1 & 0 & 0 & -1 \\ 0 & 1 & -1 & 0 & 1 & 0 & 0 & -1 & -1 & -1 & -1 & 1 \\ 0 & 0 & 0 & 1 & -1 & 0 & 1 & 1 & 0 & 0 & -1 & 0 \\ 0 & -1 & 0 & 0 & 0 & 1 & 0 & 0 & 1 & 0 & 1 & 0 \\ 0 & 0 & 0 & 0 & 1 & -1 & 0 & -1 & 1 & 0 & 0 & 0 \\ 0 & 1 & -1 & -1 & 1 & 0 & -1 & 0 & -1 & 0 & 0 & 1 \end{pmatrix}$$

Row sum of the RTD matrix is as follows

$$\begin{array}{r} \geq 60 \\ 50-59 \\ 40-49 \\ 30-39 \\ 20-29 \\ \leq 19 \end{array} \begin{bmatrix} -1 \\ -2 \\ 1 \\ 2 \\ 0 \\ -1 \end{bmatrix}.$$

The graph for $\alpha = 0.7$, is given in the next page.

We see when we take $\alpha = 0.7$ we have a variation but still the age group 50 to 59 happens to be repulsive and 20-29 continue to be passive.

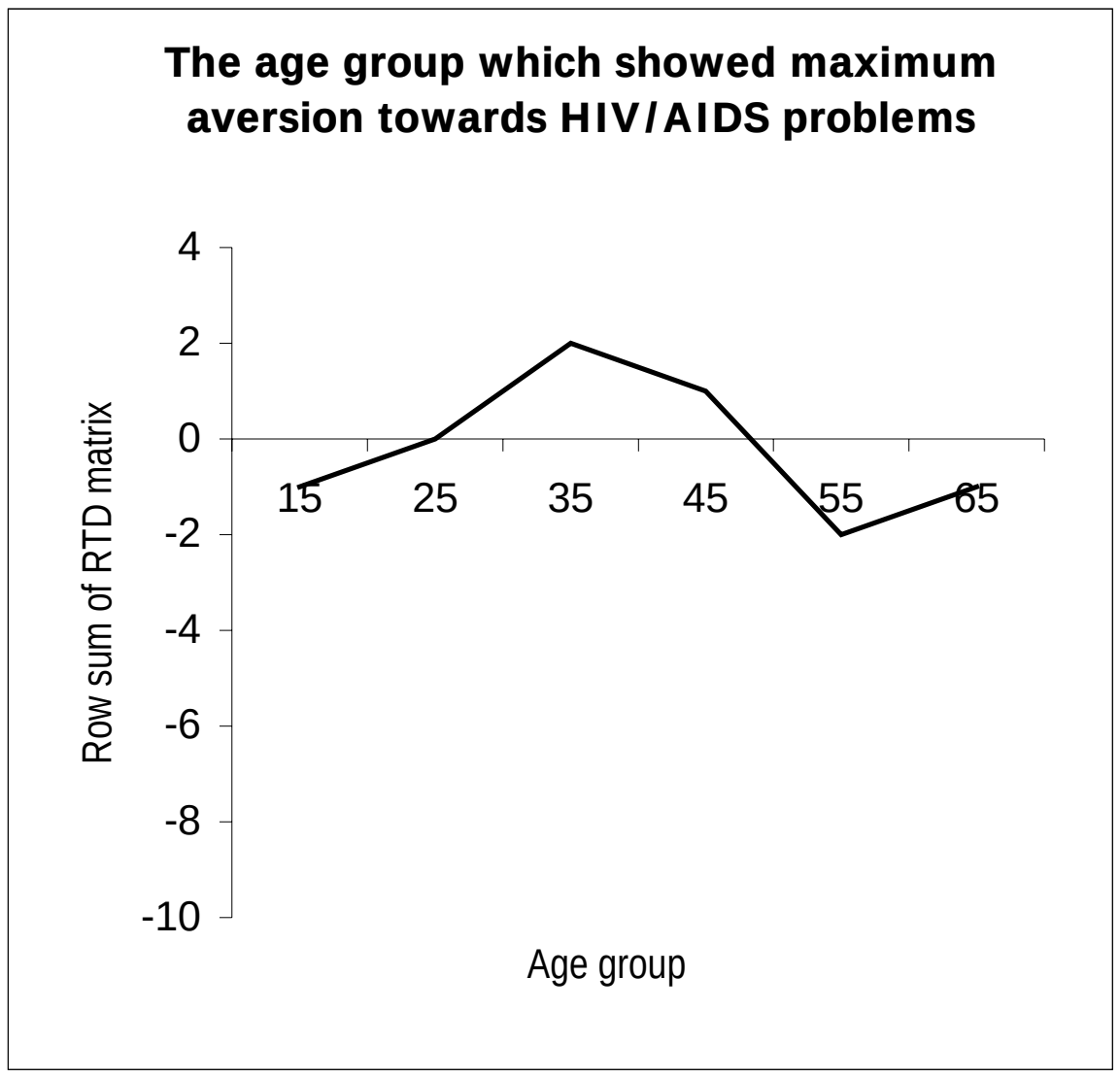

Now we work for α = 0.2 and obtain the following RTD matrix

$$\begin{array}{cccccccccccc} A_1 & A_2 & A_3 & A_4 & A_5 & A_6 & A_7 & A_8 & A_9 & A_{10} & A_{11} & A_{12} \end{array}$$
$$\begin{pmatrix} 1 & -1 & 1 & -1 & -1 & -1 & 0 & 1 & 1 & 1 & 1 & -1 \\ -1 & 1 & -1 & -1 & 1 & 1 & -1 & -1 & -1 & -1 & -1 & 1 \\ 0 & 0 & 1 & 1 & -1 & -1 & 1 & 1 & 1 & 0 & -1 & -1 \\ -1 & -1 & 0 & 0 & -1 & 0 & 1 & -1 & 1 & 1 & 1 & -1 \\ -1 & 0 & 0 & -1 & 1 & -1 & 1 & -1 & 1 & 1 & 0 & -1 \\ -1 & 1 & -1 & -1 & 1 & -1 & -1 & 1 & -1 & 1 & 1 & 1 \end{pmatrix}.$$

The row sum of RTD matrix is

$$\begin{array}{r} \geq 60 \\ 50-59 \\ 40-49 \\ 30-39 \\ 20-29 \\ \leq 19 \end{array} \begin{bmatrix} 1 \\ -4 \\ 1 \\ -1 \\ -1 \\ 0 \end{bmatrix}.$$

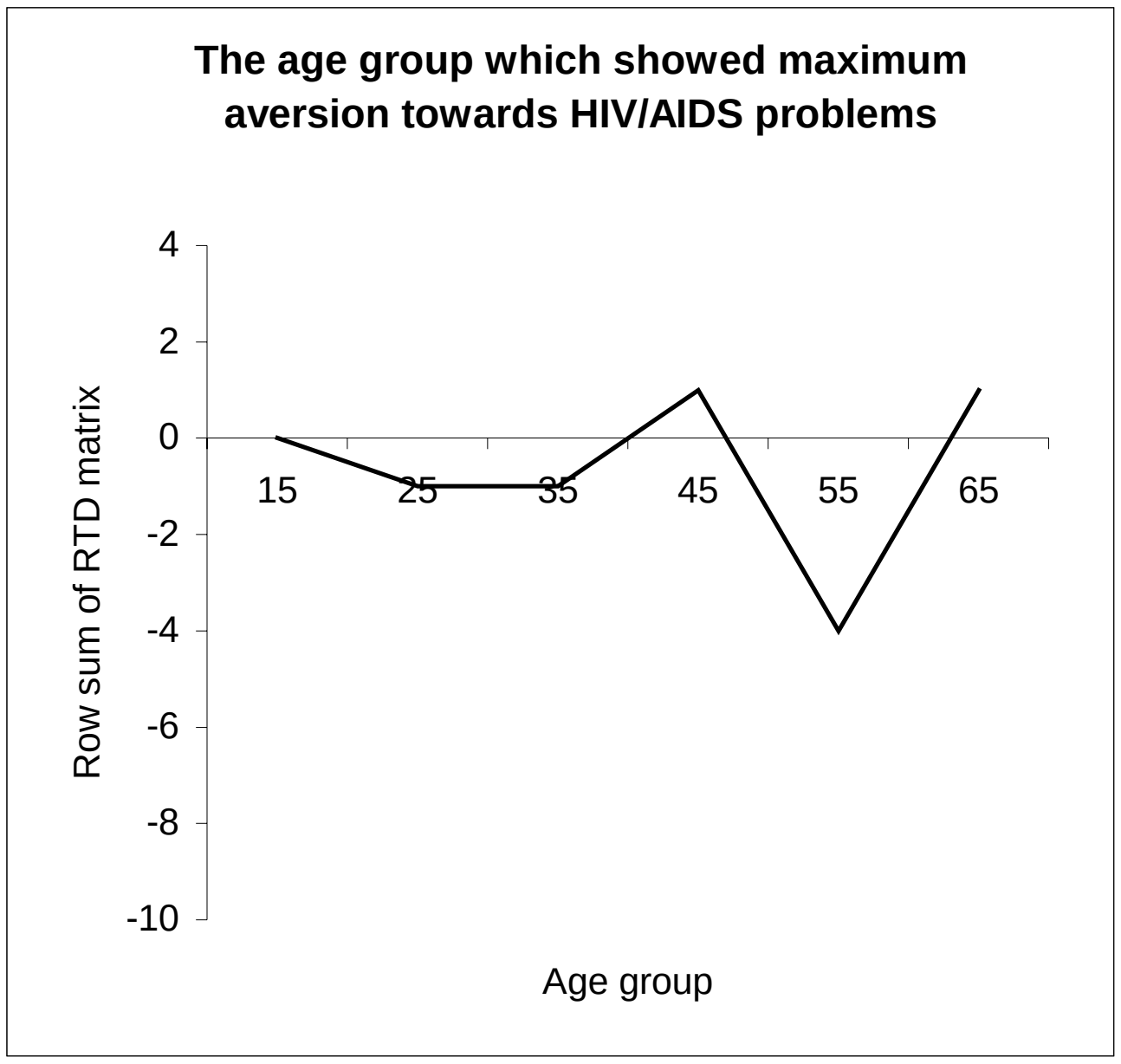

Still the same age group 50-59 happens to be averse and we see several fluctuation; thus one can understand that the choice of the value $\alpha \in [0\ 1]$ has a significant role to play in the resultant. We have worked for $\alpha = 1$, $\alpha = 0.7$ and $\alpha = 0.2$. We now give the combined RTD matrix.

The Combined Effect Time dependent data matrix (CETD matrix) for the 3 values of $\alpha$ is given below:

$$
\begin{array}{cccccccccccc}
A_1 & A_2 & A_3 & A_4 & A_5 & A_6 & A_7 & A_8 & A_9 & A_{10} & A_{11} & A_{12}
\end{array}
$$

$$
\begin{pmatrix}
3 & -3 & 3 & -3 & -3 & -3 & 0 & 3 & 2 & 1 & 1 & -2 \\
-1 & 3 & -3 & -1 & 3 & 1 & -1 & -3 & -3 & -3 & -3 & 2 \\
0 & 0 & 1 & 3 & -3 & -1 & 3 & 3 & 1 & 0 & -3 & -1 \\
-1 & -3 & 0 & 0 & -1 & 2 & 1 & -1 & 2 & 1 & 3 & -1 \\
-1 & 0 & 0 & -1 & 3 & -3 & 1 & -2 & 2 & 1 & 0 & -1 \\
-1 & 3 & -3 & -2 & 3 & -1 & -3 & 1 & -2 & 1 & 1 & 3
\end{pmatrix}
$$

Row sum of the CETD matrix

$$\begin{array}{r} \geq 60 \\ 50-59 \\ 40-49 \\ 30-39 \\ 20-29 \\ \leq 19 \end{array} \begin{bmatrix} 1 \\ -9 \\ 3 \\ 2 \\ -1 \\ 0 \end{bmatrix}.$$

The graph of the CETD matrix is given below:

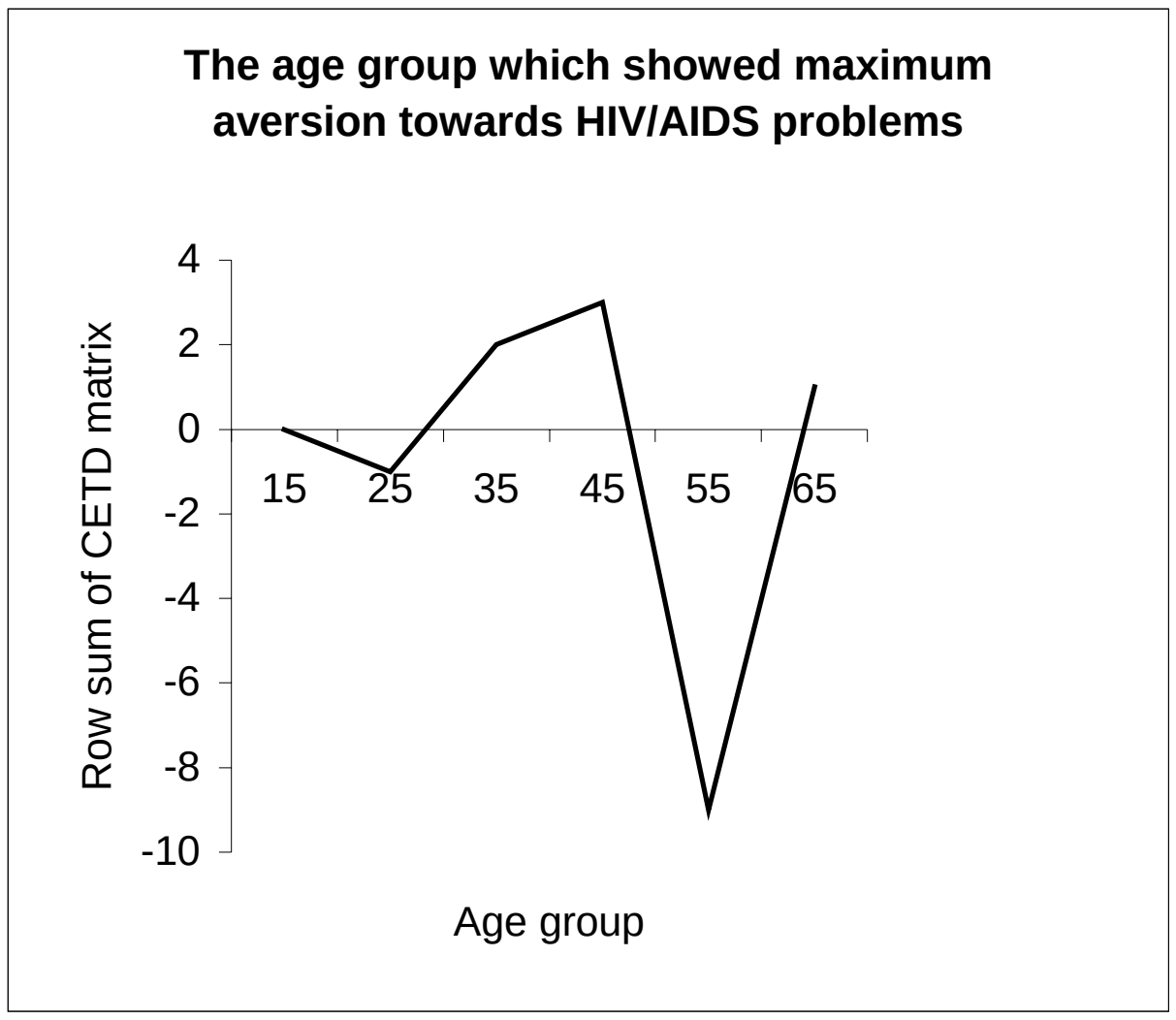

The age group which showed maximum aversion towards HIV/AIDS problems always set in 50-59.

However the age group 40-49 and 30-39 showed positive sign to work with.

The age group 20-29 however showed aversion, ≤ 19 did not react i.e., passive.

Now we proceed on to refine the age group and study so that the breath of the interval may not be 10 but just 5. Now we analyse the data with interval 5.

| Age Group | No of persons |
|---|---|
| ≥ 65 | 5 |
| 60-64 | 0 |
| 55-59 | 7 |
| 50-54 | 6 |
| 45-49 | 10 |
| 40-44 | 6 |
| 35-39 | 10 |
| 30-34 | 7 |
| 25-29 | 15 |
| 20-24 | 21 |
| 14-19 | 14 |

Now we obtain the table with 12 attributes $A_1$ to $A_{12}$ and 11 age groups. This is mainly done to make the solution more sensitive.

| | $A_1$ | $A_2$ | $A_3$ | $A_4$ | $A_5$ | $A_6$ | $A_7$ | $A_8$ | $A_9$ | $A_{10}$ | $A_{11}$ | $A_{12}$ |
|---|---|---|---|---|---|---|---|---|---|---|---|---|
| ≥ 65 | 5 | 0 | 2 | 0 | 1 | 0 | 1 | 2 | 5 | 5 | 4 | 2 |
| 60-64 | 0 | 0 | 0 | 0 | 0 | 0 | 0 | 0 | 0 | 0 | 0 | 0 |
| 55-59 | 6 | 1 | 2 | 1 | 4 | 1 | 1 | 1 | 6 | 6 | 3 | 4 |
| 50-54 | 0 | 3 | 0 | 1 | 4 | 0 | 1 | 2 | 4 | 4 | 3 | 4 |
| 45-49 | 4 | 2 | 5 | 3 | 5 | 1 | 2 | 4 | 8 | 9 | 4 | 6 |
| 40-44 | 5 | 0 | 2 | 2 | 1 | 0 | 3 | 2 | 5 | 5 | 5 | 1 |
| 35-39 | 4 | 0 | 2 | 2 | 4 | 2 | 4 | 3 | 10 | 10 | 10 | 3 |
| 30-34 | 3 | 0 | 3 | 1 | 5 | 1 | 0 | 3 | 7 | 7 | 7 | 4 |
| 25-29 | 6 | 2 | 5 | 2 | 10 | 0 | 5 | 5 | 15 | 15 | 13 | 10 |
| 20-24 | 8 | 2 | 8 | 2 | 10 | 0 | 2 | 3 | 18 | 17 | 11 | 7 |
| 14-19 | 9 | 4 | 1 | 1 | 9 | 1 | 1 | 8 | 14 | 14 | 13 | 11 |

Now we give the Raw Time Dependent Data matrix.

$$
\begin{pmatrix}
5 & 0 & 2 & 0 & 1 & 0 & 1 & 2 & 5 & 5 & 4 & 2 \\
0 & 0 & 0 & 0 & 0 & 0 & 0 & 0 & 0 & 0 & 0 & 0 \\
6 & 1 & 2 & 1 & 4 & 1 & 1 & 1 & 6 & 6 & 3 & 4 \\
0 & 3 & 0 & 1 & 4 & 0 & 1 & 2 & 4 & 4 & 3 & 4 \\
4 & 2 & 5 & 3 & 5 & 1 & 2 & 4 & 8 & 9 & 4 & 6 \\
5 & 0 & 2 & 2 & 1 & 0 & 3 & 2 & 5 & 5 & 5 & 1 \\
4 & 0 & 2 & 2 & 4 & 2 & 4 & 3 & 10 & 10 & 10 & 3 \\
3 & 0 & 3 & 1 & 5 & 1 & 0 & 3 & 7 & 7 & 7 & 4 \\
6 & 2 & 5 & 2 & 10 & 0 & 5 & 5 & 15 & 15 & 13 & 10 \\
8 & 2 & 8 & 2 & 10 & 0 & 2 & 3 & 18 & 17 & 11 & 7 \\
9 & 4 & 1 & 1 & 9 & 1 & 1 & 8 & 14 & 14 & 13 & 11
\end{pmatrix}.
$$

Now we give the Average time dependent data matrix $A = (a_{ij})$.

$$
\begin{pmatrix}
1 & 0 & .4 & 0 & .2 & 0 & .2 & 4 & 1 & 1 & .8 & .4 \\
0 & 0 & 0 & 0 & 0 & 0 & 0 & 0 & 0 & 0 & 0 & 0 \\
.86 & .14 & .29 & .14 & .57 & .14 & .14 & .14 & .86 & .86 & .43 & .57 \\
0 & .5 & 0 & .17 & .67 & 0 & .17 & .33 & .67 & .67 & .5 & .67 \\
.4 & .2 & .5 & .3 & .5 & .1 & .2 & .4 & .8 & .9 & .4 & .6 \\
.83 & 0 & .33 & .33 & .17 & 0 & .5 & .33 & .83 & .83 & .83 & .17 \\
.4 & 0 & .2 & .2 & .4 & .2 & .4 & .3 & 1 & 1 & 1 & .3 \\
.43 & 0 & .43 & .14 & .71 & .14 & 0 & .43 & 1 & 1 & 1 & .57 \\
.4 & .13 & .33 & .13 & .67 & 0 & .33 & .33 & 1 & 1 & .87 & .67 \\
.38 & .1 & .38 & .1 & .47 & 0 & .1 & .14 & .86 & .81 & .52 & .33 \\
.64 & .29 & .07 & .07 & .64 & .07 & .07 & .57 & 1 & 1 & .93 & .79
\end{pmatrix}.
$$

The average and standard deviation of each column of attribute:

| | $A_1$ | $A_2$ | $A_3$ | $A_4$ | $A_5$ | $A_6$ | $A_7$ | $A_8$ | $A_9$ | $A_{10}$ | $A_{11}$ | $A_{12}$ |
|---|---|---|---|---|---|---|---|---|---|---|---|---|
| x | .49 | .12 | .27 | .14 | .45 | .06 | .19 | .25 | .82 | .82 | .66 | .46 |
| $u_j$ | .33 | .15 | .17 | .10 | .23 | .07 | .15 | .16 | .28 | .28 | .30 | .23 |

Take α = 1 and using the mean and SD of the columns we obtain the RTD matrix.

$$\begin{array}{c} \begin{matrix} A_1 & A_2 & A_3 & A_4 & A_5 & A_6 & A_7 & A_8 & A_9 & A_{10} & A_{11} & A_{12} \end{matrix} \\ \begin{pmatrix} 1 & 0 & 0 & -1 & -1 & 0 & 0 & 0 & 0 & 0 & 0 & 0 \\ -1 & 0 & -1 & -1 & -1 & 0 & -1 & -1 & -1 & -1 & -1 & -1 \\ 1 & 0 & 0 & 0 & 0 & 1 & 0 & 0 & 0 & 0 & 0 & 0 \\ -1 & 1 & -1 & 0 & 0 & 0 & 0 & 0 & 0 & 0 & 0 & 0 \\ 0 & 0 & 1 & 1 & 0 & 0 & 0 & 0 & 0 & 0 & 0 & 0 \\ 1 & 0 & 0 & 1 & -1 & 0 & 1 & 0 & 0 & 0 & 0 & -1 \\ 0 & 0 & 0 & 0 & 0 & 1 & 1 & 0 & 0 & 0 & 1 & 0 \\ 0 & 0 & 0 & 0 & 1 & 1 & -1 & 1 & 0 & 0 & 1 & 0 \\ 0 & 0 & 0 & 0 & 0 & 0 & 0 & 0 & 0 & 0 & 0 & -1 \\ 0 & 0 & 0 & 1 & 0 & 0 & 0 & 0 & 0 & 0 & 0 & 0 \\ 0 & 1 & -1 & 0 & 0 & 0 & 1 & 0 & 0 & 1 & 1 & 1 \end{pmatrix} \end{array}.$$

The row sum of the RTD matrix is as follows

$$\begin{matrix} \geq 65 \\ 60-64 \\ 55-59 \\ 50-54 \\ 45-49 \\ 40-44 \\ 35-39 \\ 30-34 \\ 25-29 \\ 20-24 \\ 14-19 \end{matrix} \begin{bmatrix} -1 \\ -10 \\ 2 \\ -1 \\ 2 \\ 1 \\ 3 \\ 3 \\ -1 \\ 1 \\ 3 \end{bmatrix}$$

The graph for α = 1 gives the age which showed the maximum aversion towards HIV/AIDS problems is given below. The age group is taken along the X-axis and the Row sum of the CETD matrix is taken along the Y-axis.

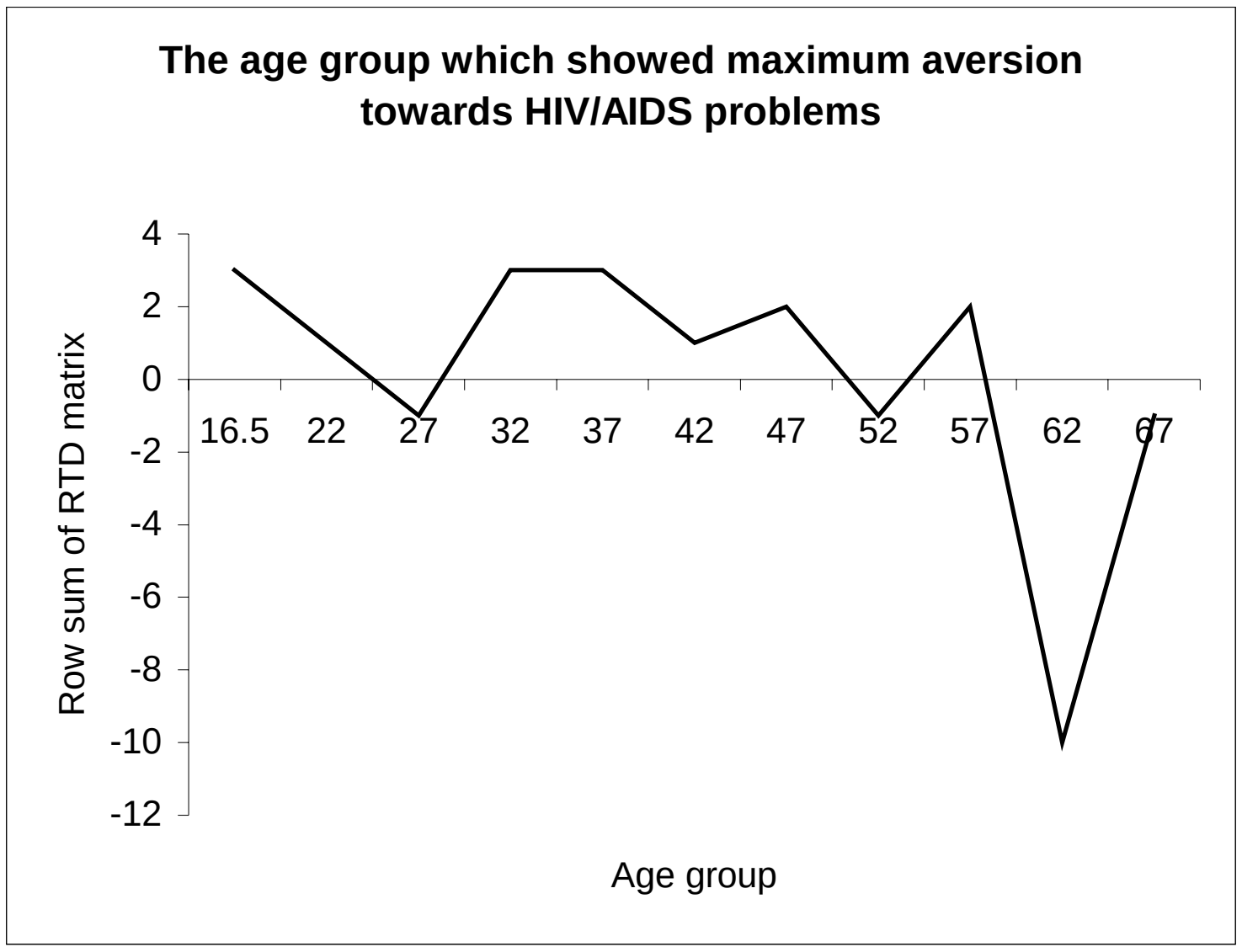

We have obtained for α = 1 a refined data which shows the age group, 35-39 and 14 –19 show the maximum concern followed by the age group 55-59, 45-49 and 30-34.

However the age groups 40-44 and 20-24 show some concern. The worst is 60-64 were none even opted to answer our questionnaire.

Now we work the same data with α = 0.7 ∈ [0 1]. Further this shows 20-24 show some positive inclination it is only the age group 25 to 29 showed negative feeling towards the patients, disease and the awareness programs.

For α = 0.7 we get the following Refined Time Dependent matrix

$$
\begin{array}{cccccccccccc}
A_1 & A_2 & A_3 & A_4 & A_5 & A_6 & A_7 & A_8 & A_9 & A_{10} & A_{11} & A_{12}
\end{array}
$$

$$
\begin{pmatrix}
1 & -1 & 1 & -1 & -1 & -1 & 0 & 1 & 0 & 0 & 0 & 0 \\
-1 & -1 & -1 & -1 & -1 & -1 & -1 & -1 & -1 & -1 & -1 & -1 \\
1 & 0 & 0 & 0 & 0 & 1 & 0 & -1 & 0 & 0 & -1 & 0 \\
-1 & 1 & -1 & 0 & 1 & -1 & 0 & 0 & 0 & 0 & 0 & 1 \\
0 & 0 & 1 & 1 & 0 & 0 & 0 & 1 & 0 & 0 & -1 & 0 \\
1 & -1 & 0 & 1 & -1 & -1 & 1 & 0 & 0 & 0 & 0 & -1 \\
0 & -1 & 0 & 0 & 0 & 1 & 1 & 0 & 0 & 0 & 1 & -1 \\
0 & -1 & 1 & 0 & 1 & 1 & -1 & 1 & 0 & 0 & 1 & 0 \\
0 & 0 & 0 & 0 & 1 & -1 & 1 & 0 & 0 & 0 & 1 & 1 \\
0 & 0 & 0 & 0 & 0 & -1 & 0 & -1 & 0 & 0 & 0 & 0 \\
0 & 1 & -1 & -1 & 1 & 0 & -1 & 1 & 0 & 0 & 1 & 1
\end{pmatrix}
$$

Now we find the related row matrix sum for $\alpha = 0.7$

$$
\begin{array}{r}
\geq 65 \\ 60-64 \\ 55-59 \\ 50-54 \\ 45-49 \\ 40-44 \\ 35-39 \\ 30-34 \\ 25-29 \\ 20-24 \\ 14-19
\end{array}
\begin{bmatrix}
-1 \\ -12 \\ 0 \\ 0 \\ 2 \\ -1 \\ 1 \\ 3 \\ 3 \\ -2 \\ 2
\end{bmatrix}
$$

The graph of the row sum for $\alpha = 0.7$ is given below:

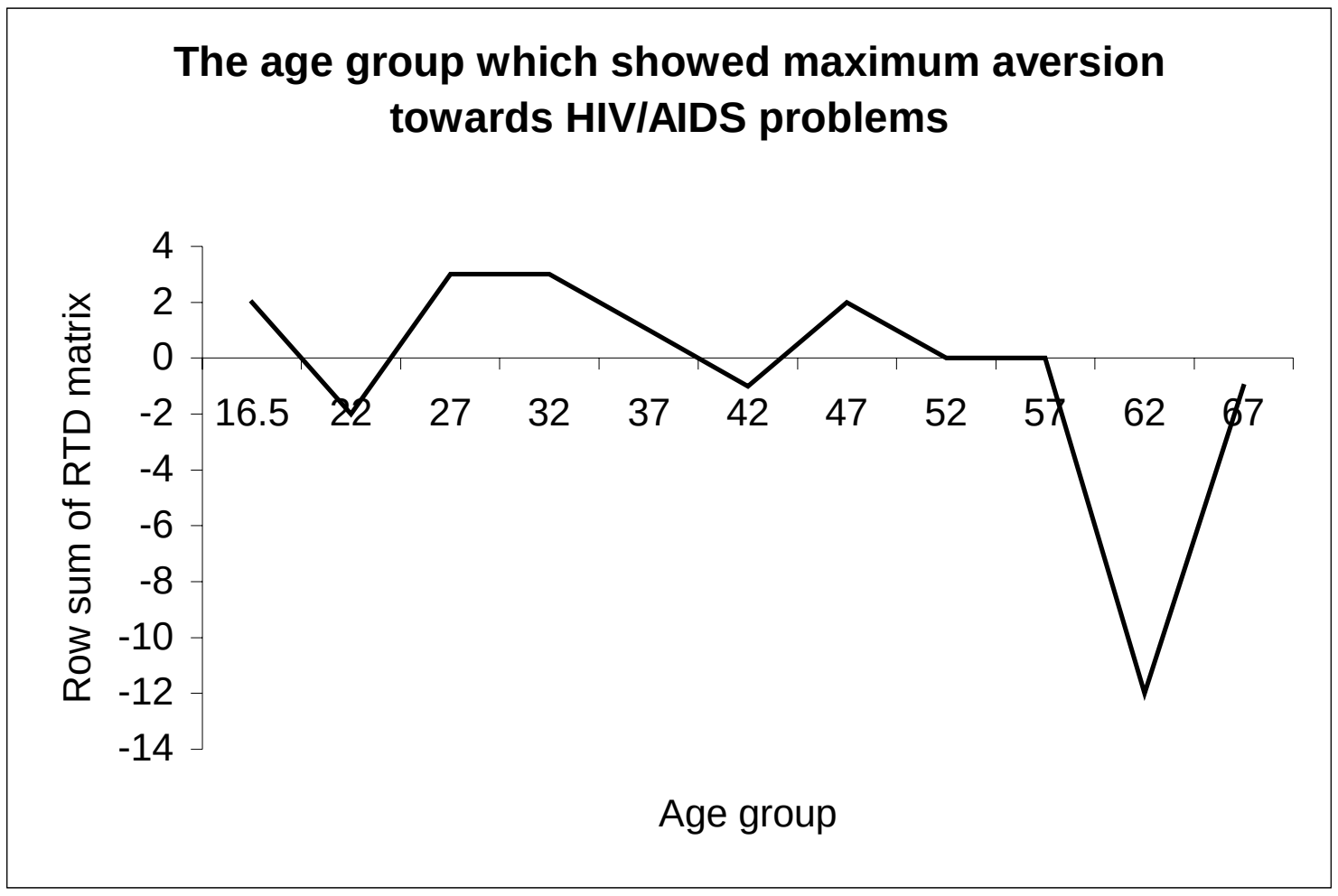

The age group 60-64 showed the maximum aversion followed by the 20-24. 55-59 and 50-54 remained passive, 30-34 and 25-29 showed maximum interest. There is some inconsistency yet not much of deviation so we proceed on to work with $\alpha = 0.2 \in$ [0 1].

The RTD matrix related to $\alpha = 0.2$ is as follows:

$$
\begin{array}{c} \begin{array}{cccccccccccc} A_1 & A_2 & A_3 & A_4 & A_5 & A_6 & A_7 & A_8 & A_9 & A_{10} & A_{11} & A_{12} \end{array} \\
\begin{pmatrix}
1 & -1 & 1 & -1 & -1 & -1 & 0 & 1 & 1 & 1 & 1 & -1 \\
-1 & -1 & -1 & -1 & -1 & -1 & -1 & -1 & -1 & -1 & -1 & -1 \\
1 & 0 & 0 & 0 & 1 & 1 & -1 & -1 & 0 & 0 & -1 & 1 \\
-1 & 1 & -1 & 0 & 1 & 1 & 0 & 1 & -1 & 0 & -1 & 1 \\
0 & 1 & 1 & 1 & 1 & 1 & 0 & 1 & 0 & 1 & -1 & 1 \\
1 & -1 & 1 & 1 & -1 & 1 & 1 & 1 & 0 & 0 & 1 & -1 \\
0 & -1 & 0 & 0 & -1 & 1 & 1 & 1 & 1 & 1 & 1 & -1 \\
0 & -1 & 1 & 0 & 1 & 1 & -1 & 1 & 1 & 1 & 1 & 1 \\
0 & 0 & 1 & 0 & 1 & 1 & 1 & 1 & 1 & 1 & 1 & 1 \\
0 & 0 & 1 & 0 & 0 & 1 & -1 & -1 & 0 & 0 & -1 & -1 \\
0 & 1 & -1 & -1 & 1 & 1 & -1 & 1 & 1 & 1 & 1 & 1
\end{pmatrix} \end{array} .
$$

Now we obtain the row sum of the RTD matrix for α = 0.2.

$$\begin{array}{c} \geq 65 \\ 60-64 \\ 55-59 \\ 50-54 \\ 45-49 \\ 40-44 \\ 35-39 \\ 30-34 \\ 25-29 \\ 20-24 \\ 14-19 \end{array}\begin{bmatrix} -1 \\ -12 \\ 1 \\ 1 \\ 7 \\ 4 \\ 3 \\ 6 \\ 9 \\ -2 \\ 5 \end{bmatrix}.$$

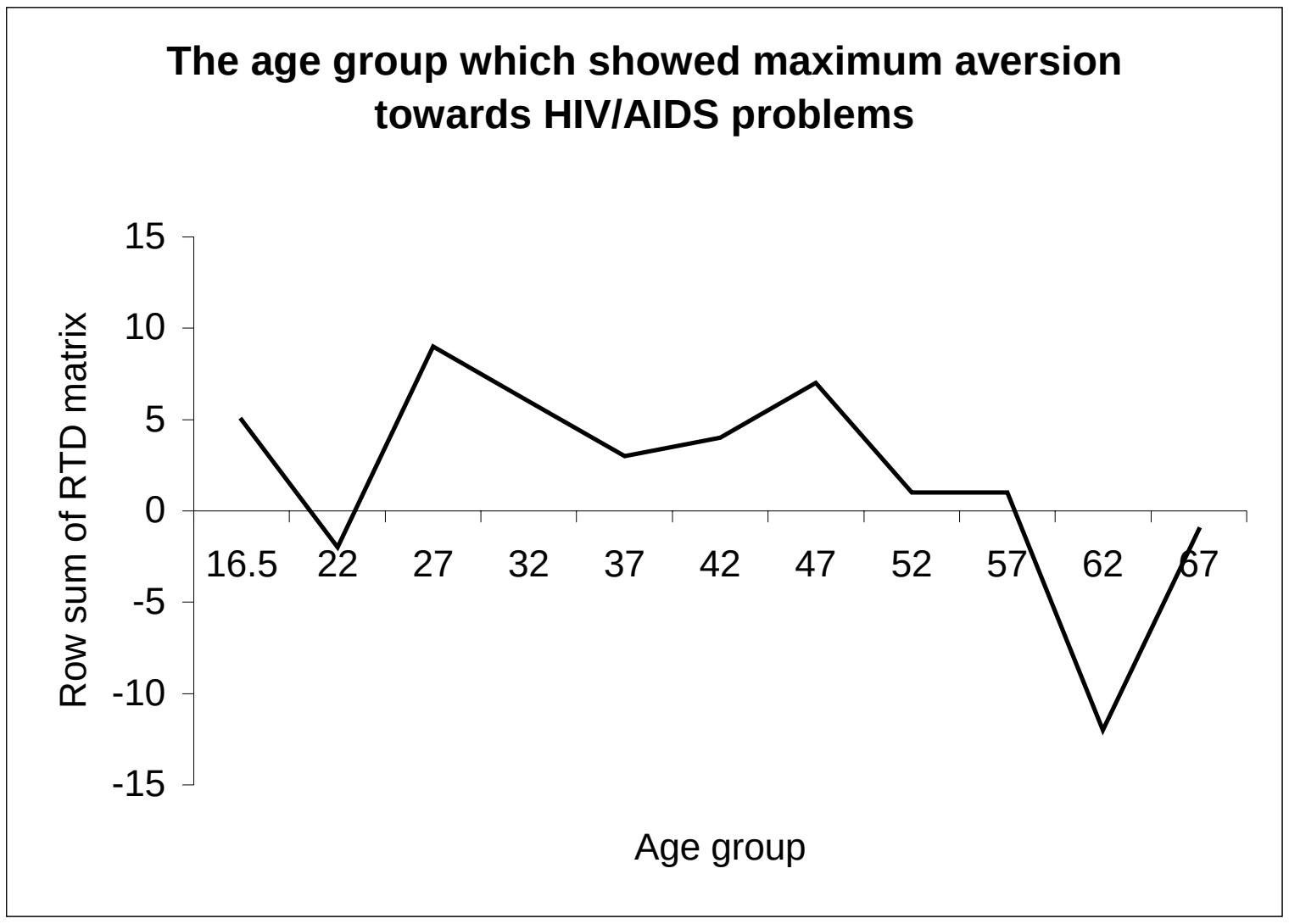

The row sum shows that 60-64 is the people who are maximum aversed with the HIV/AIDS patients or its awareness program. The age group 25-29 showed maximum interest followed maximum interest followed by the age groups 46-49, 30-34 and 40-44. As we have very varying row sum we obtain the CETD matrix for these 3 parametric values $\alpha = 1$, $\alpha = 0.7$ and $\alpha = 0.2$.

The CETD matrix is as follows:

$$
\begin{array}{cccccccccccc}
A_1 & A_2 & A_3 & A_4 & A_5 & A_6 & A_7 & A_8 & A_9 & A_{10} & A_{11} & A_{12}
\end{array}
$$

$$
\begin{pmatrix}
3 & -2 & 2 & -3 & -3 & -2 & 0 & 2 & 1 & 1 & 1 & -1 \\
-3 & -2 & -3 & -3 & -3 & -2 & -3 & -3 & -3 & -3 & -3 & -3 \\
3 & 0 & 0 & 0 & 1 & 3 & -1 & -2 & 0 & 0 & -2 & 1 \\
-3 & 3 & -3 & 0 & 2 & 0 & 0 & 1 & -1 & 0 & -1 & 2 \\
0 & 1 & 3 & 3 & 1 & 1 & 0 & 2 & 0 & 1 & -2 & 1 \\
3 & -2 & 1 & 3 & -3 & 0 & 3 & 1 & 0 & 0 & 2 & -3 \\
0 & -2 & 0 & 0 & -1 & 3 & 3 & 1 & 1 & 1 & 3 & -2 \\
0 & -2 & 2 & 0 & 3 & 3 & -3 & 3 & 1 & 1 & 3 & 1 \\
0 & 0 & 1 & 0 & 2 & 0 & 2 & 1 & 1 & 1 & 2 & 1 \\
0 & 0 & 1 & 1 & 0 & 0 & -1 & -2 & 0 & 0 & -1 & -1 \\
0 & 3 & -3 & -2 & 2 & 1 & -2 & 3 & 1 & 1 & 3 & 3
\end{pmatrix}.
$$

Now the row sum of the CETD matrix

$$
\begin{array}{r}
\geq 65 \\ 60-64 \\ 55-59 \\ 50-54 \\ 45-49 \\ 40-44 \\ 35-39 \\ 30-34 \\ 25-29 \\ 20-24 \\ 14-19
\end{array}
\begin{bmatrix}
-1 \\ -34 \\ 3 \\ 0 \\ 11 \\ 5 \\ 7 \\ 12 \\ 11 \\ -3 \\ 10
\end{bmatrix}.
$$

The graph of the row sums of the CETD matrix is as follows:

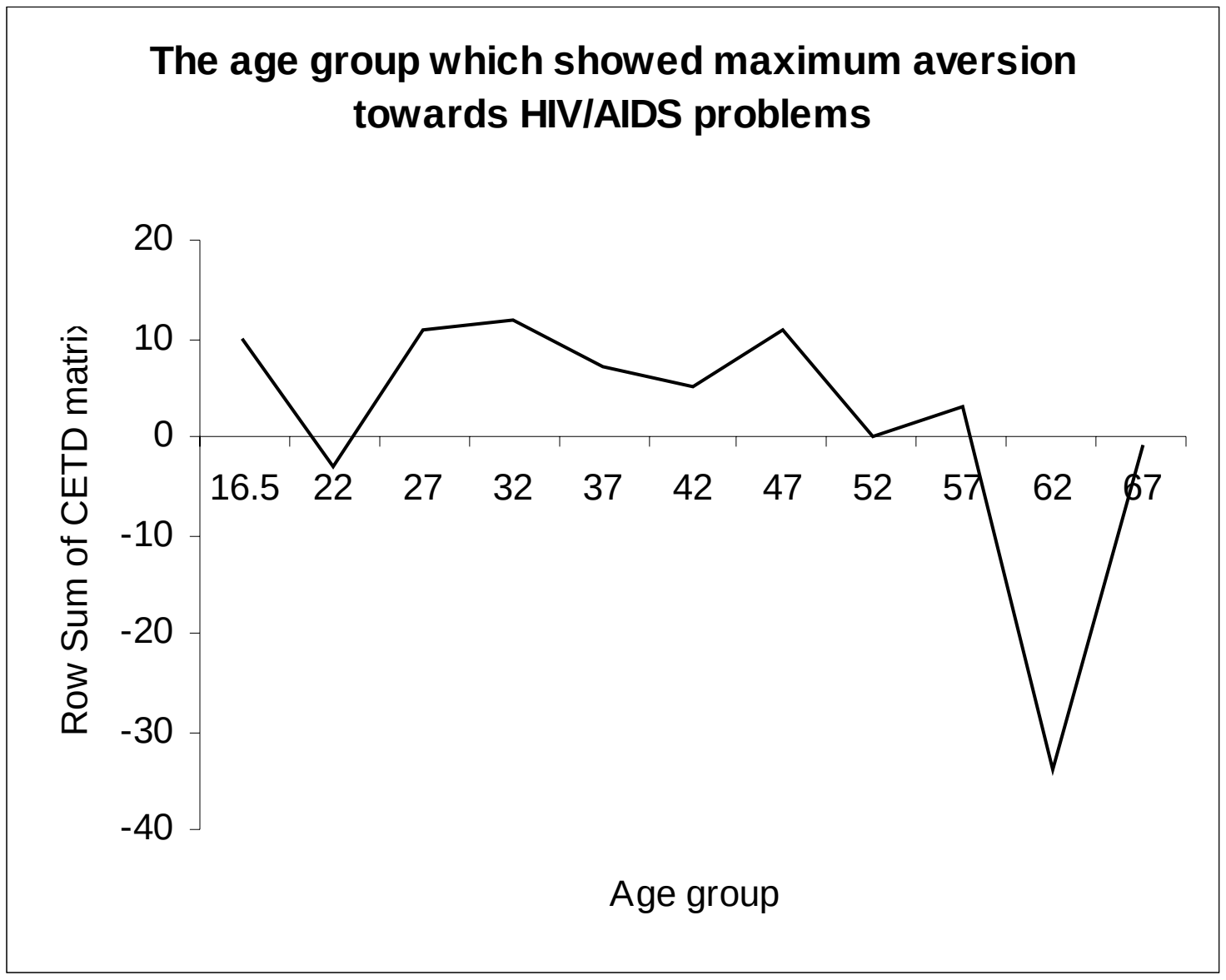

Thus we see the age group 20 to 24 showed aversion. But the age group 60-64 and ≥ 65 maximum aversion.

The study revealed that the age group 60-64 never wanted to involve with HIV/AIDS awareness or so. Also the age group ≥ 65 were averse but less aversed than the age group 20-24. The people in the age group 30-34 showed maximum concern equally followed by the age group 45-49 and 25-29. The age group 14-19 showed equal interest. However the age group 50-54 did not like to involve themselves. So people after 50 did not show any concern over HIV/AIDS awareness program or about the patients.

Similarly the age group 20-24 were not interested might be more involved in their getting employment or further studies or future plans or so and are aversed to talk or think about HIV/AIDS as they were engrossed about themselves. Similarly the age group 50-54 were more serious about their family affairs after retirement plans like building house, marriage of their children and so on.

From the graph it is very clear that the age-group from 24-57 showed positive response. Using the C-program given in appendix 3 Conculsions derived are given in chapter 4.

## 3.2 Use of 2 - adaptive fuzzy model to analyze the problem

In this section we for the first time use the 2-adaptive fuzzy model $C_2$ [$M_1$, $M_2$] (where $M_1$ is a fuzzy matrix or a CETD matrix and $M_2$ is a BAM model) to analyze the problem of sociological and psychological feelings of public towards the HIV/AIDS disease, related concepts and the awareness programes.

We have already discussed about the functioning of this model in chapter 2 page 41. We have already analyzed using the single model viz. fuzzy (CETD) matrices in section 3.1.

Now we use the 2-adaptive fuzzy model $C_2$ [$M_1$, $M_2$] to analyze the views of public about HIV/AIDS and the awareness program.

We have taken $M_1$ to be the fuzzy matrix model and $M_2$ to be the BAM model.

Now we have already worked with 3 parameters. $\alpha = 1, 0.7$ and 0.3.

Now the expert has worked with two more parameters and gives the following CETD matrix for 5 values of $\alpha$ given by $M_1'$.

$$M'_1 = \begin{array}{c} \begin{matrix} A_1 & A_2 & A_3 & A_4 & A_5 & A_6 & A_7 & A_8 & A_9 & A_{10} & A_{11} & A_{12} \end{matrix} \\ \begin{pmatrix} 5 & -5 & 4 & -5 & 0 & -5 & -3 & 0 & 5 & 5 & 3 & 0 \\ 0 & -1 & -5 & -4 & 2 & -5 & -4 & -2 & 2 & 3 & 0 & 2 \\ 1 & -4 & 0 & 0 & 0 & -5 & 0 & 0 & 5 & 5 & 1 & 0 \\ 0 & -5 & -1 & -4 & 1 & -3 & -1 & -1 & 5 & 5 & 5 & 0 \\ 0 & -4 & 0 & -4 & 1 & -5 & -3 & -3 & 5 & 5 & 2 & 0 \\ 0 & -3 & -5 & -5 & 2 & -5 & -5 & 0 & 3 & 5 & 3 & 3 \end{pmatrix} \end{array}$$

The row sum is as follows:

$$\begin{array}{r} \geq 60 \\ 50-59 \\ 40-49 \\ 30-39 \\ 20-29 \\ \leq 19 \end{array}\begin{bmatrix} 4 \\ -2 \\ 3 \\ 1 \\ -6 \\ -7 \end{bmatrix}$$

The graph for the CETD matrix is given as follows

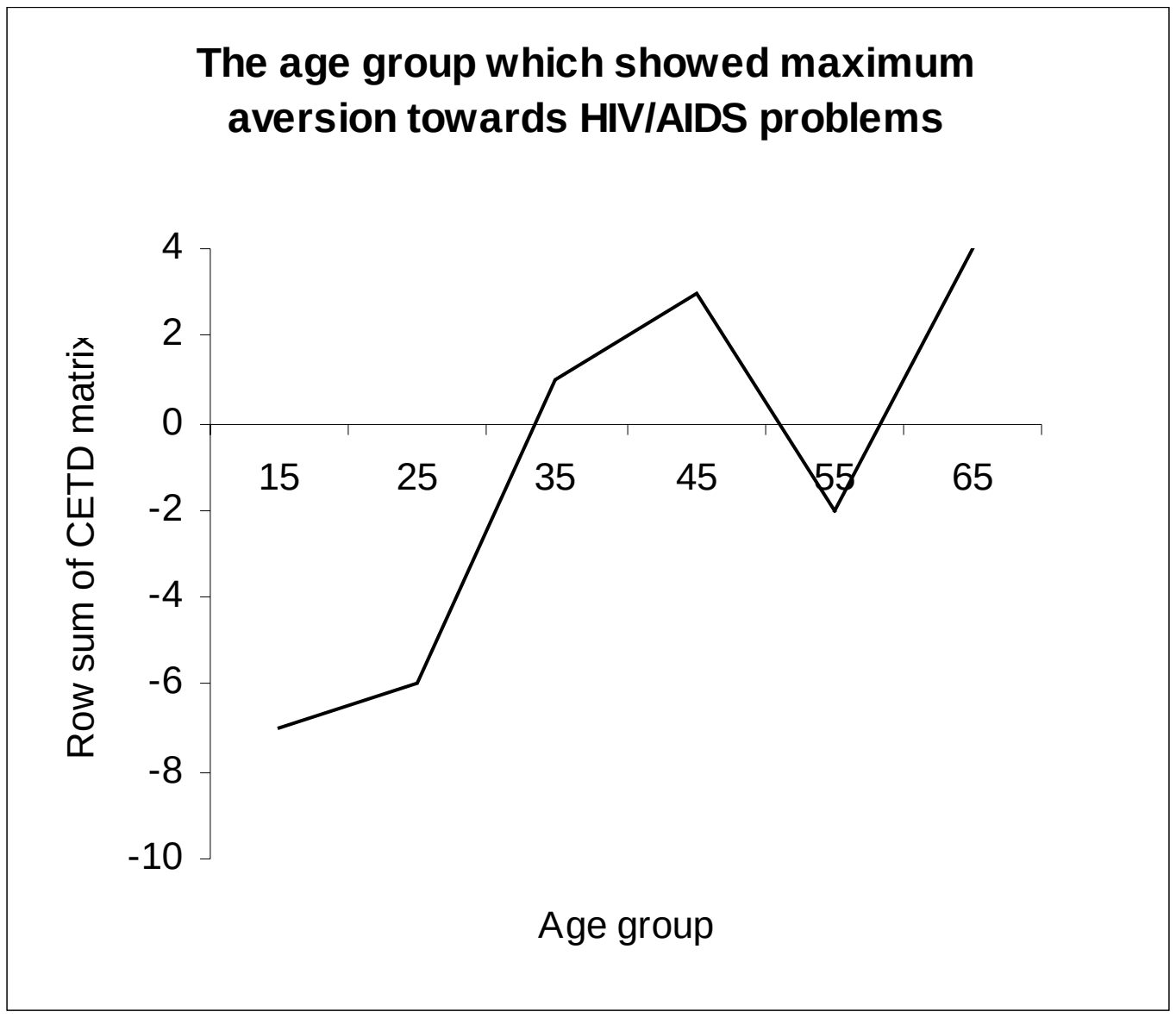

Now the expert can use this matrix $M'_1$ itself as a BAM model and find the effect of any initial vector on the BAM model.

Take $M'_1$ as the BAM model. Put $M'_1 = M$

$$M = \begin{pmatrix} 5 & -5 & 4 & -5 & 0 & -5 & -3 & 0 & 5 & 5 & 3 & 0 \\ 0 & -1 & -5 & -4 & 2 & -5 & -4 & -2 & 2 & 3 & 0 & 2 \\ 1 & -4 & 0 & 0 & 0 & -5 & 0 & 0 & 5 & 5 & 1 & 0 \\ 0 & -5 & -1 & -4 & 1 & -3 & -1 & -1 & 5 & 5 & 5 & 0 \\ 0 & -4 & 0 & -4 & 1 & -5 & -3 & -3 & 5 & 5 & 2 & 0 \\ 0 & -3 & -5 & -5 & 2 & -5 & -5 & 0 & 3 & 5 & 3 & 3 \end{pmatrix}$$

(columns labelled $A_1$ $A_2$ $A_3$ $A_4$ $A_5$ $A_6$ $A_7$ $A_8$ $A_9$ $A_{10}$ $A_{11}$ $A_{12}$)

Let $X_k$ = (3 –1 0 –2 0 1 –2 1 0 5 0 4) at the time k, the initial vector is given such that social stigma is given the highest value 5 followed by sympathetic towards the HIV/AIDS patients, public is averse about HIV/AIDS patients, followed careless and unconcerned with $A_3$, $A_5$, $A_9$ and $A_{11}$ taking 0 and $A_4$ and $A_7$, taking the value –2.

The binary signal vector got by using the activation function gives

$S(X_k)$ = (1 0 0 0 0 1 0 1 0 1 0 1)
$S(X_k)M^t$ = (5 –2, 1, 1, –3, 3)
= $Y_{k+1}$ .

$S(Y_{k+1})$ = (1 0 1 1 0 1)
$S(Y_{k+1})\,M^t$ = (6, –17, –2, –14, 3, –18, –9, –1, 18, 20, 12, 3)
= $X_{k+2}$.

$S(X_{k+2})$ = (1 0 0 0 1 0 0 0 1 1 1 1)
$S(X_{k+2})\,M^t$ = (18, 9, 12, 16, 13, 16)
= $Y_{k+3}$.

$S(Y_{k+3})$ = (1 1 1 1 1 1)
$S(Y_{k+3})M$ = (6, –22, –7, –22, 6, –28, –16, –6, 25, 28, 14, 5)

$$= X_{k+4}.$$

$$S(X_{k+4}) = (1\ 0\ 0\ 0\ 1\ 0\ 0\ 0\ 1\ 1\ 1\ 1)$$
$$= S(X_{k+2}) .$$

Thus we see the stability of the dynamical system is given by {(1 0 0 0 1 0 0 0 1 1 1 1), (1 1 1 1 1 1)}. Thus all the age group people feel the same. They are averse towards HIV/AIDS at the same time they are concerned about it. They fear it, they associative social stigma, they are shy and conservative about it, but however they are sympathetic towards the patients. This is what the model built using only the views of them say. Interested readers can study with other $X_k$'s. The most important thing is that this change will take place with in four years for stability is reached at $X_{k+4}$.

Thus we have seen using the 2 adaptive fuzzy model we have derived the conclusions.

Suppose one is given the BAM model in the scale [–t t] we can make use of the same matrix using it as the number it as the number of parameters t used in the CETD matrix, as this is a 2 adaptive fuzzy model.

It is still important to mention that using any the ATD matrix also one can use a desired transformation function and get the BAM.

It is in the hands of the expert and on the nature of the problem to get the BAM model having the ATD, RTD and CETD matrices. Even the transformations functions can be different also the scales can be different.

Now using the RTD matrix given in page 65 we use a function f to convert it into the BAM model, which is given by $M'_2$. It is important to mention here that even the values given in the matrix $M'_2$ are at the discretion of the expert and the function he defines.

The BAM model $M_1$ in the scale [–5, 5] is as follows: This is obtained from the RTD matrix with refined age intervals.

$$
\begin{array}{cccccccccccc}
A_1 & A_2 & A_3 & A_4 & A_5 & A_6 & A_7 & A_8 & A_9 & A_{10} & A_{11} & A_{12}
\end{array}
$$

$$
\begin{pmatrix}
5 & -5 & 0 & -5 & -4 & -5 & -4 & 0 & 5 & 5 & 4 & 0 \\
-5 & -5 & -5 & -5 & -5 & -5 & -5 & -5 & -5 & -5 & -5 & -5 \\
4 & -4 & -3 & -4 & 1 & -4 & -4 & -4 & 4 & 4 & 0 & 1 \\
-5 & -5 & -5 & -4 & 2 & -5 & -4 & 0 & 2 & 2 & 0 & 2 \\
0 & -4 & 0 & -3 & 0 & -5 & -4 & 0 & 4 & 4 & 0 & 3 \\
4 & -5 & 0 & 0 & -4 & -5 & 0 & 0 & 4 & 4 & 4 & -4 \\
0 & -5 & -4 & -4 & 0 & -4 & 0 & -1 & 5 & 5 & 5 & -1 \\
0 & -5 & 0 & -4 & 3 & -4 & -5 & 0 & 5 & 5 & 5 & 1 \\
0 & -4 & 0 & -4 & 2 & -5 & 0 & 0 & 5 & 5 & 4 & 2 \\
0 & -5 & 0 & -5 & 0 & -5 & -5 & -4 & 4 & 4 & 1 & 0 \\
2 & -3 & -5 & -5 & 2 & -5 & -5 & 1 & 5 & 5 & 4 & 3
\end{pmatrix}
$$

Now we take the same initial vector $X_k = (3, -1, 0, -2, 0, 1, -2, 1, 0, 5, 0, 4)$ at the time k mainly for comparative purposes as the initial vector. The effect of $X_k$ on the dynamical system

$S(X_k) = (1\ 0\ 0\ 0\ 0\ 1\ 0\ 1\ 0\ 1\ 0\ 1)$

$S(X_k)\, M_1^t = (5\ {-25}\ 1\ {-6}\ 2\ {-1}\ {-1}\ 2\ 2\ {-5}\ 6) = Y_{k+1}$.

$S(Y_{k+1}) = (1\ 0\ 1\ 0\ 1\ 0\ 0\ 1\ 1\ 0\ 1)$

$S(Y_{k+1})\, M_1^t = (11\ {-25}\ {-8}\ {-25}\ 4\ {-28}\ {-22}\ {-3}\ 28\ 28\ 17\ 10) = X_{k+2}$

$S(X_{k+2}) = (1\ 0\ 0\ 0\ 1\ 0\ 0\ 0\ 1\ 1\ 1\ 1)$

$S(X_{k+2})\, M_1^t = (15\ {-30}\ 14\ 3\ 11\ \ 8\ 14\ \ 19\ 18\ \ 9\ \ 21) = Y_{k+3}$.

$S(Y_{k+3}) = (1\ 0\ 1\ 1\ 1\ 1\ 1\ 1\ 1\ 1\ 1)$

$S(Y_{k+3})\, M_1^t = (10\ {-45}\ {-17}\ {-38}\ \ 2\ {-47}\ {-31}\ {-8}\ 43\ 43\ 27\ 7)$

$= X_{k+4}$.

$S(X_{k+4}) = (1\ 0\ 0\ 0\ 1\ 0\ 0\ 0\ 1\ 1\ 1\ 1) = S(X_{k+2})$

Thus the resultant binary pair is a fixed point given by {(1 0 1 1 1 1 1 1 1 1 1) (1 0 0 0 1 0 0 0 1 1 1 1)}. But when we refine the age group we see that the age group 60-64 alone are the group not interested about any awareness of HIV/AIDS.

**Chapter Four**

# PUBLIC ATTITUDE AND AWARENESS ABOUT HIV/AIDS

In this chapter we give a brief introduction to the psychological and social problems concerning AIDS, which were analyzed by us while doing fieldwork. This forms the first section of this chapter. The second section gives an outline of the individual interviews.

## 4.1. Introduction

We have earlier interviewed HIV/AIDS patients (migrant labourers and rural women) and based on our interviews we analyzed the data using some type of fuzzy model. As our interviews of these patients clearly revealed that the patients suffered a lot of social stigma especially in their villages we felt it deem fit to interview the public because we wanted to know how people in general viewed the disease. The area of study was Chennai city. We had taken an interview of about 150 people from various walks of life, and here we have given 101 interviews. Those interviews that did not convey any message or contribute anything to our data were dropped. Many of them were reluctant to give interviews and in fact made some derogatory remarks about HIV/AIDS patients.

We have chosen all categories of people: barbers and businessmen, post-graduates and pavement dwellers. Our interviews include the opinions of educated and uneducated, the poor and the rich, professionals and unemployed. We are happy that the sample is as representative of the entire society and as heterogeneous as possible. The linguistic questionnaire used by us is enclosed in the book. We observed that women spoke more freely than men and they did not shun us. Public were unanimous that the government should help the HIV/AIDS infected children. Unlike men no women ever showed displeasure to be interviewed or feared somebody would see them. Conclusions based on our study and mathematical analysis is given in the final chapter.

At the outset we want to record that getting the questionnaire filled or getting an interview from the public was a very painful process. Some people treated us worse than HIV/AIDS patients. Some used very abusive terms; some tried to even manhandle my students. Thus we spent the maximum period of time only interviewing the public. People dodged us. They wanted a private place to talk about HIV/AIDS. All our experience only showed the fear that the public had over the disease. They had the fear of being branded as 'characterless' men. The social stigma can be only attributed to the way one becomes infected by HIV/AIDS. They discriminated HIV/AIDS infected persons in their minds, even if they did not do so openly.

The fear about AIDS makes people discriminate HIV/AIDS patients, ill-treat and harass them. Can anyone with a little rationalism ever trouble a person who is already seriously suffering with both mental and physical problems? The discrimination of the HIV/AIDS patients is often based on prejudice resulting from an irrational fear. In India, the public fears the disease as a social stigma. Thus ignorance has fueled irrational fears of infection, leading to prejudice and discrimination against and stigmatization of people with HIV. Moreover, religious leaders do not talk about HIV/AIDS. Some feel that it is a product of karma. People in general hold the opinion that men and CSWs affected by HIV/AIDS are "guilty perpetrators."

Even in developed nations like the United States, this disease was stigmatized in the early 1980s. In fact local school officials have expelled children from schools and parents groups have organized boycotts by keeping their children out of schools when court ordered that HIV-positive children have a right to attend classes. Children and infants have been rejected from home placement and seropositive parents have been denied child custody and visitation rights despite the antidiscrimination policy resolutions in late 1980s by American Medical Association. (Encyclopedia of AIDS, Fitzroy Dearborn Publishers, Chicago (1998)).

Even as early as 1980, healthcare professionals in the US refused to treat HIV-positive patients at emergency rooms, clinics, and medical and dental offices, denied them access to ambulances, nursing homes, dialysis treatment, surgery and general nursing care. Medical labs have refused to perform work on specimens from HIV-sero positive persons. There are numerous instances where healthcare professionals have breached patient's confidentiality about their HIV-status leading to harassment and stigmatization in the workplace, termination or reassignment of job duties and/or loss of insurance coverage. Funeral-home workers have been charged with refusing to prepare bodies for burial, leaving the deceased wrapped in body bag, or refusing burial outright. Many insurance companies either deny coverage to potential HIV-positive policyholders or provide only limited coverage and cancel insurance policies to existing holders once a person's seropositive status is known.

Workplace discrimination also runs high, stemming in part from employers' irrational fear of transmission or from their fear about loss in productivity, or the possibility of higher insurance costs or that the entire business will become stigmatized by association. Employers often fire, reassign or do not hire employees based solely on their HIV status. On several occasions people have been tested for HIV against their will or without their knowledge – as a basis for employment, entry into a hospital or insurance coverage. Landlords have evicted or denied rentals to tenants, a number of airlines have been sued for refusing to allow persons with AIDS to fly and people have been denied services in restaurants, beauty salons, barbershops and gyms.

In fact the Mayor of a West Virginia town ordered the municipal public swimming pool emptied, scrubbed, disinfected and refilled after a person with AIDS swam in it and forbade the person from returning.

Seropositive prisoners are often segregated and denied access to appropriate medical and legal resources. Sometimes they are forced to wear differently coloured uniforms, denied visitation rights, excluded from participation in programs granted to other inmates, including work release and parole programs, vocational and recreational programs and even religious services, denied probation bail or fair trial based on their HIV status.

In government programs they face delays and face unnecessary obstacles. In addition general harassment and violence have been directed against people with HIV/AIDS. Often victims of this violence are deterred from filing complaints and seeing a case through legal channels for fear of retribution and loss of anonymity.

We see that in Tamil Nadu, the behavior of city dwellers and rural people towards HIV/AIDS patients is markedly different in several ways. It is very difficult for an HIV/AIDS patient to live in his village once the villagers come to know that he/ she is seropositive. They are so barbaric and ruthless that they deny them all basic needs like ration, water, milk, kerosene etc. and ultimately chase them out of the village. This arises because of fear of infection and the social stigma associated with this disease.

On the contrary, people in the cities may not make a big fuss. But once they come to know that someone is seropositive, they will immediately cut all contact with that person.

People living with HIV/AIDS in cities in Tamil Nadu do not face many problems because they have a safe anonymity. They can buy anything they want from any shop. Only their own relatives or friends, who know the seropositive status, might have the capacity to discriminate them. Only private hospitals or clinics might turn them out. The economic status of the HIV/AIDS patient, plays a major role. If the patient is rich, all doctors are willing to serve him.

It is pathetic state of affairs that people living with HIV/AIDS are subject to discrimination in their workplace. If the house-owner comes to know that the tenant has HIV/AIDS, he/she is mercilessly treated and asked to vacate the premises as soon as possible. There have been cases where lawyers have refused to take up legal matters of HIV+ people alleging that they will be labeled AIDS lawyers.

But it must be said to the credit of the government and NGOs that in the past five years, there has been a remarkable shift in public opinion and attitude. The increased sensitization has ensured that people are more aware and more tolerant.

It is pertinent to mention that several colleges have taken up the mission to spread awareness about HIV/AIDS. The stigma associated with AIDS has visibly reduced to a great level. The efforts of the Tamil Nadu State AIDS Control Society (TANSACS) in this regard have been monumental. It has successfully pioneered a major shift in public perspective, changing an attitude of contempt into an attitude of care.

## 4.2. Interviews

This chapter gives the verbatim interviews of 101 people from various walks of life, like professors, teachers, students, women, government staff, workers, uneducated daily wagers, drivers, housewives and so on. The interviews were taken in Chennai. We present their views about HIV/AIDS patients, problems of the orphaned children whose parents died because of HIV/AIDS, children living with HIV/AIDS and their suggestions for a better living. We also concentrate on what should be done by the government to prevent the further spread of HIV/AIDS and types of awareness programmes that must be given. Nearly half of interviewed persons talked openly and positively about helping people living with HIV/AIDS while the remaining showed utter contempt and dislike for these patients. A vast majority was of the opinion that HIV/AIDS is a curse and punishment for immoral activities. Many members of the public, roughly 40%, were even reluctant to speak with us about HIV/AIDS. They felt bad even to utter the word 'AIDS'.

At times, the opinions of some of the individuals were contradictory, inconsistent and incoherent, yet we have given them as it is. Our main motivation was not to gloss up their opinion but to present it faithfully. During the interviews, if we pointed out contradictions, the interviewees were most likely to discontinue expressing their opinion. For instance, while we approached people to answer the questionnaire, several of them ran away from us. Sixteen of them got irritated with us and abruptly broke away the interview. We therefore let them speak without any interruption. People were reluctant to speak about HIV/AIDS. When we supplied our questionnaire and requested for an interview, many of them stared at us and even refused to look at the questionnaire. We felt that a person who speaks about, or serves, people living with HIV/AIDS is also treated with contempt. Some people insulted and ill-treated us while some chased us away as if we were unseeables, when we approached them for interviews. A few of them asked us to 'get lost', while few of them ran away from us saying "you have no other job." Some ill-treated us by saying that people get HIV/AIDS due to bad character. They could not understand why we were highly bothered about the patients. To be accurate, only 42 people gave a second thought about the plight of these people affected with HIV/AIDS. Only 8 of the interviewed felt little pity.

In many cases, we could not convince them that everybody has erred in life and that we should ensure that people with HIV/AIDS should not suffer psychologically. Three out of the 101 people interviewed said that they preferred a murderer to a person infected with HIV/AIDS. Even literates failed to fill the questionnaire. They said it would be easier if they answered our questions. Interviewing migrant labourers and rural women who were living with HIV/AIDS was easier compared to getting public opinion. We can clearly say that the public shuns the topic of HIV/AIDS and they are very reluctant to discuss it. We now give a brief summary of interviews with the general Public. We also give the additional observations made by us in some cases while interviewing them.

## 1

On 10-3-2003, we interviewed a professor employed in a private college in the city. He was 53 years old. Both his children are engineering graduates. His wife was also a professor employed in a government college. He is from the upper middle class section of society. He was also a social worker, helping the poor and needy. He said several facts about the HIV/AIDS affected patients. He said that the social stigma faced by the people living with HIV/AIDS patients is much more than the stigma faced by Dalits. It is impossible for an HIV/AIDS patient even to buy a cup of tea from any shop in his native village. Even the shop owners boycott them. They fear to sell provisions to them. Even a tumbler of water will be refused to HIV/AIDS patients in remote villages. This professor is from a remote village but is settled in Chennai for the past three decades. He says the stigma associated with HIV/AIDS is so strong that when the village comes to know that a person has HIV/AIDS, their families have no option but to leave the village.

He says this is one of the reasons why men affected with HIV/AIDS do not inform their wives or other family members about it. In fact many persons who have contracted the HIV virus die untreated fearing stigma. He says it is impossible for a traditional society like ours to openly accept men who visit CSWs and end up becoming HIV/AIDS patients. In his opinion, the public will not easily accept people living with HIV/AIDS. The only option is to strive to gain people's sympathy and remove the stigma so that those affected with HIV/AIDS can lead a peaceful mental life.

He suggested that the government should provide some percentage of reservation for people living with HIV/AIDS, in the public and private sectors. The government should reserve some medical and engineering seats for children whose parent(s) are/were affected by HIV/AIDS. Education must be provided free of cost to children whose parents suffer(ed) from HIV/AIDS.

All these measures will not only rehabilitate the people living with HIV/AIDS but also make the public think about AIDS related problem. Further when the government helps these patients, men/women affected with HIV/AIDS may come out openly and their psychological problems due to social stigma would be reduced. Only the rural uneducated people with HIV/AIDS suffer this social stigma to the maximum extent. In the urban areas, people living with HIV/AIDS take treatment secretly. People in the age group 18 to 24 years, especially men, when they came to know that they are affected by HIV/AIDS committed suicide without any second thought.

He suggests that children affected by HIV/AIDS should be provided free education, boarding, medical aid and nutritious food by the school. Similarly, women who were infected by their husbands should be given free medical aid and some training to sustain them.

*Additional Observations:*

*He is one among a handful who openly discussed all the problems faced by people living with HIV/AIDS. He gave several suggestions.*

**2**

We interviewed a 49-year-old homemaker. She was a graduate in fine arts. She has two sons—one of them is married. The younger son is in the final year of school. Her husband is a bank employee. She is from an upper middle class background. Her first son is an engineering graduate working in the PWD. She wanted her identity to be kept confidential. She became aware of HIV/AIDS only from TV and radio advertisements. She says HIV/AIDS is a very vulgar disease and a lot of social stigma is associated with it. She refused to be interviewed about HIV/AIDS in the presence of her family members. We interviewed her when she was alone.

She says the root-cause for men becoming victims of HIV/AIDS is CSWs and alcohol. She says that school and college students get the HIV/AIDS infection because of their deviant behaviour caused by movies and porn magazines and the Internet culture. She blames media for its pathetic portrayal of man-woman relationships.

She says that movies play a major role in making people commit errors. She had seen one of her relatives seeing highly pornographic pictures on the Internet. He is in his late fifties. She felt difficult to talk with him after knowing that he has such bad tastes. We cannot say anything about anyone and even age is not a bar for any bad action. She feels that most of the cinemas only kindle the base instincts of people. After seeing such films, students in their adolescent age can make any mistake.

She does not spare educated men who are affected with HIV/AIDS. She believes that rich men who suffer from HIV/AIDS, hide the news from their family and only the concerned doctor knows about it. She is of the opinion that these men live with the disease for over 15 years, because they are very careful and have learnt a lesson for life. When we asked her why most of the poor rural men when affected by HIV/AIDS died in a short span of less than 5 years, she said it was because they led careless lives, didn't eat nutritiously and also visited more CSWs. They not only increase their chance of getting AIDS, but also spread the disease among the various CSWs. She said that only god could save these rural men. She wanted the state to protect the children affected by HIV/AIDS. She wanted them to be given good education and care. Government must help these children or else they will have no future.

She felt sad for the plight of women infected by their husbands. She feels that these women who suffer greatly in their last stages should be given good food, free medicine and volunteers must be employed in hospitals to take care of them. Women who worked all their lives for others are abandoned with nobody to even care for them. These women must be protected and they should be allowed to spend the last days on earth peacefully.

She feels that the advertisement and propaganda given about HIV/AIDS is not sufficient. It could be improved greatly and students must be involved in spreading the awareness about HIV/AIDS. She pities CSWs because in her opinion only the erring men infect CSWs. In most cases, CSWs are helpless, they cannot ask or advice these men to have safe sex since these men are in a drunken state most of the time. It is a sad state of affairs that most CSWs suffer for no fault of their own.

She says the social stigma associated with HIV/AIDS cannot be wiped out in a day. It will take over decades of care, cure and prevention of HIV/AIDS. Even though she sympathizes the CSWs, women and children living with HIV/AIDS, she is very critical and harsh towards men who are living with the disease. She feels that they are antisocial elements worse than a thief or a murderer because their act is a cold sin. We felt that it was impossible for us to make her realize that these people should be given some protection.

She feels that it is the duty of the parents to educate their children and give them awareness about HIV/AIDS through open discussions. She feels that it is better if this subject is broached in every family, than through classrooms.

**3**

On 21-3-2003, we met an assistant professor of an engineering college in Chennai. He is a doctorate. He is 35 years old, is married and has a son. His wife is a homemaker. He believes in yoga and conducts classes on spirituality in the after hours of his college. When we asked about HIV/AIDS he said it is karma of ones acts. Nobody needs to help a person living with HIV/AIDS because only through undergoing physical suffering and mental torture one may be purified. He further added that it was the karma of his own choice. So he feels that the public need not pity people affected by HIV/AIDS. When we put forth the question as to why they were being ill-treated or avoided, he said it was mainly because of the fear and stigma associated with the disease. When we asked why most of the Hindu religious leaders do not talk about HIV/AIDS or help these patients, he said that the only reason is that the disease is associated with sex; so Hindu spiritual leaders do not wish to discourse about it.

He did not give any suggestions to help children who were infected by HIV/AIDS. He attributes that also to their deeds in past lives. Likewise, he says that women are infected by their husbands because of their karma. Thus we felt that it was not possible to receive any answer or suggestions from him.

*Additional Observation:*

*This interview is very special and unique case because the interviewee attributes the cause of HIV/AIDS to present/ past Karma. We felt it was not worthwhile to discuss any relevant information with him.*

# 4

On 1-4-2003, we interviewed a 68-year-old grocery store owner. He is a Muslim from Bangalore settled in Chennai. He lamented that even the boys who work in his shop are intoxicated by cinema. More than prayer or food they crave only for movies. He says that movies are the root-cause for men or women leading promiscuous lives. We said that according to a rough estimate and our own statistics "4% of the rural uneducated women and 5% of migrant labourers who were infected with HIV/AIDS were Muslims" He said that nowadays Muslims have become money-minded. They did not value faith in god or the code of Koran. They easily visit CSWs when they get enough money. He deeply sympathized the women (and children) who have become HIV/AIDS infected because of their husbands (and parents). He said that everyday, the mullah must talk about HIV/AIDS in the mosque after prayers. This would spread awareness. This will make those on the 'wrong' paths to think twice before committing a crime. He said that a strong code of law must be formulated in order to severely punish erring husbands. He suggested that it could even be imprisonment for up to 5 years (with medical treatment given in the prison). When they are put in special prisons with medical aid, certainly these men cannot visit CSWs or infect others. If this type of strong law is made and implemented people will be afraid to visit CSWs. Secondly this will drastically bring down the number of HIV/AIDS patients. Although his views are very severe and cruel, we didn't argue with him. We simply let him express all his feelings. He wants the government to give pension to the CSWs so that they will no longer practice the trade. He also wants the children and women who were affected by HIV/AIDS to be given monetary help by philanthropists. Regarding public opinion, which is often not favorable to the HIV/AIDS patients, he says that the public will not sympathize or accept the patients because they suffer due to their own heinous acts. Consequently, he feels that it is not a problem if they lead a life of desperation and dejection. If the public did not scorn them, they would go on erring. More people will be tempted to follow in their footsteps. He says there is a lack of awareness of the seriousness of HIV/AIDS. He says that in the major railway stations every half an hour an advertisement about HIV/AIDS must be displayed on the giant screens. He feels that this will reach a large number of people who will attentively see it.

He says the same sort of advertisement can be done in large bus terminuses where buses from several states go after a halt since it will also reach a larger mobile population.

He is against sex education or HIV/AIDS awareness education in schools. He feels it will only make the children carry out experiments with their own bodies. Parents must take up that responsibility. He says that parents urge their children to excel in life, but very few of them impart value education. He says it is the most important and the primary duty of the parents to educate their children in this matter. Most of the parents do not have moral stand because they themselves lead a bad life. He felt that if they led a honest and upright life, how would their children ever be bad? We did not counter him. As a last question, we asked his education qualification. He said he was a $9^{th}$ standard dropout. Since he had to assist his father in business, he had come to Chennai with his father to run this provision shop when he was just 14. He says that he has two sons who are textile engineers in a Coimbatore mill. He says they are well settled. He also added that he had only one wife, and likewise his sons were also monogamous.

*Additional Observation:*

*Several of his views need second thought. At the outset we were surprised to get such views from a grocer who had studied only up to $9^{th}$ standard. Throughout the interview, he was sincere and serious. He feels character building must start at home.*

## 5

On 4-4-2003 we interviewed a clerk who worked in the Secretariat in Chennai. He is an arts graduate from a city college. When asked his age he asked us to put it down as 40 years. He is married. He declined to give his wife or children's ages. He is settled in Chennai for two generations. He was not very willing to be interviewed but when several other secretariat staff chased us away saying that we had no better job etc., he took pity on us and asked why we were doing this and what was our qualification etc. in a sympathetic way. Our interview was conducted in a lonely place since he preferred it. We record this information here, just to point out that because of stigma people shun speaking about HIV/AIDS publicly.

When we asked why men became infected with HIV/AIDS, he said there are several reasons. He gave the reasons as follows:

1. Due to bad company and bad habits, men visit CSWs and become infected. He is well aware of HIV/AIDS. (One of his close friends was an HIV/AIDS patient).

2. In the early years of the marriage, when their wife leaves to her maternal home for delivery of children, most men visit CSWs and become infected.

3. Men, who stay away from their family for many days, seek CSWs.

4. Due to frustration and quarrels at home, a few men visit CSWs.

He says whatever be the cause one cannot talk about it openly. He says that even in his office there are some people with HIV/AIDS who say that they suffer from TB for the fear of losing their jobs and ill-treated by their colleagues.

He says in the Indian social setup, it is impossible to accept people living with HIV/AIDS and treat them respectfully. Yet, the public mindset can be changed regarding children infected by HIV/AIDS and orphaned children. He is certain that people will be willing to help them. Women infected by their husbands can be rehabilitated. But he feels men infected by CSWs cannot be accepted or pardoned by the public. They will be stigmatized even after their death. This plight cannot be changed. The government must give free education, medical aid and free food and shelter for HIV/AIDS affected children. It should extend support to women infected by their husbands.

When a person comes to know that he/she has HIV/AIDS, the immediate reaction is to commit suicide. Government should counsel these patients before they are informed about this disease. More importantly, the family members of these patients must be counseled. Otherwise, no one can save a life. It is important to see that these men infected with HIV/AIDS are isolated so that they do not continue to frequent CSWs. Unless this is done systematically, it is impossible to control the spread of HIV/AIDS. He is of the opinion that once these men know that they are HIV positive, they refuse to have safe sex (use condoms) with CSWs. They are so deceitful and infect their innocent wife ruthlessly. He repeatedly said that isolation (quarantine) of these HIV/AIDS patients is mandatory. When we asked him how we can isolate them (because it is a very inhuman thing to do), he said that just like Yerwadi for the mentally-ill patients or like care-homes for spastic children, all HIV/AIDS affected men can be put in the outskirts of the city where they will be given free medical aid, food and also allowed to take part in activities like carpentry, basket-making etc. Since homosexuality is not a major problem in India, such quarantine will stop the spread of AIDS to the general populace. These infected patients must not be allowed to leave the place. They can be given all sorts of entertainment like films; T.V, dailies, magazines etc; sports and yoga can also be included in their daily regimen. According to him this is the only solution to stop the spread of HIV/AIDS. We did not argue with him because this was just his opinion. We tried our best to convince him against his idea, but he got upset, and said that it was getting late for him and he left us.

*Additional Observation:*

*He was very firm on the isolation of these infected men. He was full of support for HIV/AIDS infected children and women. Some of his views need a second thought.*

## 6

On 5-4-2003, we interviewed a 50-year-old Sri Lankan refugee. She has been living in Red Hills for the past 30 years. She knows about HIV/AIDS. Her parents married her to a Tamilian in north Chennai. She has 3 children; all of them are married and well settled. She is uneducated. Her children are also school dropouts. They do some petty business in the same area. Economically, they belong to the middle class.

She spoke of women becoming passive victims of HIV/AIDS. Her husband had all bad habits: drinking alcohol, smoking and visiting CSWs. When she came to know from some neighbours and friends that her husband had visited CSWs, she restrained from having any sort of sexual relations with him.

Several times, he used to beat her badly. Once or twice, he had even attacked her with a knife. She feels her husband may have died of HIV/AIDS, but she is not sure. She gave him the best of medical treatment but no doctor could cure him. Everyone said his days were counted.

She felt that men consider women as objects. As long as women passively obey men, they cannot be saved from HIV/AIDS. She says it is 'macho' behaviour that makes men beat their wives and be violent on them. If women hit back or resort to violence to protect themselves, certainly men will learn a lesson. This is her contention. She is free of diseases. She earns her living by tailoring. She has earned enough money to support her children and marry them off.

She feels that the movies/ serials on TV are of bad taste and has ruined women. Her daughter-in-laws are addicted to TV and she detests it. Women are portrayed either in tears or in complete happiness.

They weep and laugh so spontaneously, in contrast with men who never show out any emotion. Most of the TV serials are exaggerations. This has ruined young minds. It is very difficult to liberate their mind from this artificial and unreal stuff. She has failed in her attempts to change her own grand children. She feels vulgarity in TV and movies bred wrong behaviour which leads to HIV/AIDS.

Because people want to experience sex at a young age, even adolescent boys seek sex. Boys bring CSWs to their hostel rooms or visit lodges with CSWs. Some girls leave hostels in the weekends on the pretext of visiting their parents but they do not go home. Instead they are with men in some lodge. She has come to know about all this because of gossiping with her neighbours. If such practices continue, youth would be the largest population group with HIV/AIDS.

## 7

On 15th April 2003, we interviewed a 19-year-old boy who is an auto-driver by profession. He is a 9th standard dropout. He felt that education would not help him earn a livelihood or support his family. He lives in an urban slum near Adyar. He was happy to share his views on HIV/AIDS patients. He is fully aware of how HIV/AIDS spreads. He was not shy or ashamed to speak about it openly. He acknowledged that one of his close relatives had become infected with HIV/AIDS and was now an inpatient of the Tambaram Sanatorium. He visits that patient regularly because he likes him. That man had been a cab driver, and he had taught driving to this boy.

All his friends know about HIV/AIDS. He feels that almost all city dwellers are well aware of HIV/AIDS. Only people from rural areas may not be aware of it. He adds that unlike in Bombay where CSWs are concentrated in some red-light areas, in Chennai they are spread uniformly all over the city. They are available in residential areas, in slums, in five-star hotels and every customer is aware of it. He says that the public might not know how to identify them, but those who seek CSWs are well aware of it. He heard stories of the notorious Auto Shankar from his friends, who was an expert pimp running the flesh-trade in Chennai. Though Auto Shankar was hanged, there are several such people in the city. He feels that the nexus between the police and pimps is the main reason for the flourishing of CSWs. Ultimately these CSWs do not get even proper medical aid or a good meal every day. The pimps take a major portion of their earnings. It was alarming that a young boy knew so much about CSWs. He further says that some upper class women also work as CSWs for very rich clients. Most of their work is based in luxury hotels. He feels they may also have HIV/AIDS. It will appear to the world that these women live in luxury. He feels that though he is poor, he is better than these women. He is very proud that he earns by very honest means and he is happy with his conscience. He says that he has at times driven these ladies from their homes to the hotel and back again, late in the nights.

He said that the disease is publicized and stigmatized whenever poor people fall prey to it. He is sure that several rich people suffer from AIDS, but they have the money and means to take treatment secretly. When we questioned him about children affected with HIV/AIDS, he said they need special care. They should be given all material and moral support. When they become adults they can be given jobs by the government according to their capacities so that they can be self-sufficient. He feels that children affected by cancer, physical and mental disabilities are helped by support groups, likewise some groups should help HIV/AIDS affected children. For instance, children affected with cancer have special care centers, they are given good food. But children affected with HIV/AIDS lack access to nutritious food. Although good food has the capacity to prolong their lives and help their immune levels, they rarely get to eat it.

This interviewee showed sympathy and pity for HIV/AIDS infected children, women and CSWs. However he did not show any form of support for the rich, since he is of the opinion that they can manage everything with money.

He feels that awareness programs, government subsidy and public support must first reach rural, uneducated poor people and next the uneducated poor migrants in the city. He says awareness programs in rural areas are not sufficient.

The first step to AIDS eradication must be a ban on liquor. Only after getting drunk, men visit CSWs. He says even before any step is taken ban must be made on liquor. Unless liquor is banned, or the cost of liquor is increased exorbitantly (in order to keep it out of the purchasing power of poor people), the awareness program will not be successful. Most rural men take pride in stating that they are drunkards.

A psychological study of HIV/AIDS infected men is a necessity for any awareness program to be successful. According to this boy, infected men hold a grudge against all women. In a HIV/AIDS affected man’s uppermost mind, he always blames women (mostly CSWs) as the reason for getting infected. To avenge them, he tries to infect as many women or CSWs as possible. The infected man forgets that by this rash behaviour, he only losses what little remains of his health.

His relative in the Tambaram Sanatorium used to visit CSWs even when he was an inpatient. It was his only past time. This boy came to know this through the ayahs and ward boys in the hospital.

This auto-driver was very matured for his age. His experience and awareness for HIV/AIDS was more than what is usually expected.

When we asked for a solution to stop the spread of HIV/AIDS, he says that the only solution is that each individual should have self-control. No other solution is possible.

## 8

On 16$^{th}$ April 2003, we interviewed a vegetable vendor. She is from Tanjore. She came to Chennai five years back along with her family since she could not support them by working as an agricultural coolie. She is happy with the way her business has developed in these years. She is uneducated. Her husband died recently due to liver failure: he was an alcoholic and the migration to Chennai helped him to be stone drunk on all days. She has a son and two daughters. Her son is a fruit vendor. Her two daughters help her with her work.

She is fully aware of HIV/AIDS. She says it changes the basic appearance of a person. When we asked her why men get HIV/AIDS, she says it is their arrogance that makes them seek CSWs and ultimately makes them patients. HIV/AIDS spreads because men are not faithful to their wives. Often married men seek sex from CSWs and other women. Only self-realization and self control can prevent this.

She was also aware of the children affected by HIV/AIDS. She feels first free and good medical care is essential for them. They should be counseled and given moral, physical and monetary support. She says they are innocent. For no fault of their own they are made into a stigmatized element. Schools fail to admit them. Society views them as sinners. She was very sympathetic towards them.

Then we questioned her about women who were infected by their male partners. She felt very sad for their plight in the hospitals. Only women organizations should take up their issue and get them some help from the government.

She thinks HIV/AIDS is like leprosy because it disfigures the patient externally. They often have skin ailments. Sometimes the big boils they suffer cannot even be looked at. The scabies stinks so badly that no one comes forward to touch the patients for the fear that they may also get it. We asked her how she was aware of all these symptoms. She says one of her husband's distant relatives was diagnosed to be suffering from AIDS. When she was admitted to a Chennai hospital, this woman used to go and visit her occasionally taking food and clothes for her. We asked her about what sort of awareness program must be given. She says people in the cities are well aware of HIV/AIDS. If they do some mistake it is not out of ignorance. Even she came to know about HIV/AIDS only after she came to Chennai. So she feels outreach programs in rural areas that are specifically targeted at women will help people like her. Even if one woman knows about the disease, she will spread the information to all her friends. During elections, people go around in vans and autos fixed with mike sets and ask for votes. In the same way, they can campaign and spread awareness for HIV/AIDS. People should also be given adequate information about the symptoms of HIV/AIDS. Otherwise, they will not even go and take a test. This should be a weekly or fortnightly program. Songs about the disease can be played to spread awareness in tea shops, bus-stands and so on.

She is of the opinion that even doctors are afraid to come near certain patients who have skin diseases and they stink horribly. How will anyone ever be prepared to serve them? To avoid further spread of the disease these persons can be killed with injections. We were shocked and told her it is wrong to take a life. She said that they were going to die soon, so why they should suffer and infect others? She was not listening to any of our arguments in this matter.

When we asked her about HIV/AIDS awareness education to school children, she replied that there was no need. Children these days learn everything from cinema and films. She pointed to a cinema poster in front of us and asked us, ‘how can a youth be self-controlled when he sees such posters? Even married men are affected, so what about the youngsters?’ She says 30 years ago cinemas were not so vulgar. When an HIV/AIDS patient is considered a social stigma, why are lewd posters not considered with contempt? She says that unless cinema producers and magazine owners think of public welfare and not of money; no one can control the spread of HIV/AIDS. She cursed that film actors and actresses, producers, magazines and TV channel owners should become infected with HIV/AIDS for polluting the minds of so many people and making them patients of a dreaded disease.

## 9

On 30 April 2002 we interviewed a post-graduate student. He is just 21 years old. He is from an educated, middle class family. His father is a bank officer and his mother is a teacher in a private city school.

When we entered the college campus with our questionnaire, students did not run away from us (like government officers or policemen did). On the contrary, they welcomed us with tea and snacks. They were very friendly, with an urge to contribute for the cause. They flocked together in groups and seriously participated in the discussion. Their collective opinion was taken under the name of one single student because the other students shared and voiced the same thoughts. We have suppressed the name of the college to protect identities. It was a blessing in disguise that the students whom we met were studying sociology. We put forth the question about awareness of HIV/AIDS among students. He says almost all boys are well aware of HIV/AIDS.

We asked them what they did to spread awareness of HIV/AIDS being sociology students. They said that the only means for them to make the whole college participate was to grandly celebrate AIDS day on December 1. They would organize programmes at least for 5 hours. It would have songs, skits, debates etc. They had even invited an HIV/AIDS patient and a doctor to give their firsthand experience about the disease.

When we counter-questioned him, saying that a day was very short and that it may or may not have a good impact on students and faculty; they replied that perhaps they could organize the programmes for a whole week. They were even enthusiastic about the idea of going on camps to rural areas and spreading awareness there. As sociology students, they were perturbed by the obscenity and vulgarity in the representations of women in the mass media. Such a devaluation of women will spoil and mislead youngsters. Only because women are not given dignity and respect, prostitution is able to sustain as a profession. He says that all this contributes to the spread of HIV/AIDS. A student activist said that he knew even students from a premier institution in the city were caught in brothels when a police raid took place. The authorities of the institution had a very difficult time saving the name of the institution from being ruined.

When some students come to know that they have HIV/AIDS, they commit suicide. But in their suicide notes they blame parents/ teachers for scolding or bad marks as the reason. He says even school students are not free from fear of HIV/AIDS because of their risk behaviour. He cites instances of school children turning up in government hospitals wearing their own school uniform to take the HIV/AIDS test. This indicates that risk behaviour is possibly prevalent among some sections of students.

He says each college should have a hundi for collecting funds for HIV/AIDS patients. This will certainly have an indirect impact, because not only will the students think of AIDS, but they will also try to help. He also feels that students should organize lot of fund-raising events for HIV/AIDS. They could also periodically visit these patients. As a sociology student, he feels that HIV/AIDS is more a social problem than a medical problem.

## 10

On 2nd May 2002, we interviewed a 45-year-old barber from the Anna Nagar area. He owned a saloon. It was his hereditary profession for the past four generations. His son is doing his engineering degree in a city college and plans to go to the United States for higher education. He views himself as a middle class man. He studied only up to 8th standard. In fact he was a dropout from school, he had repeatedly failed in English. Now he regrets that he had wasted his school days.

We asked his opinion on HIV/AIDS. He says it is a ‘wrong’ disease that one gets by following a ‘wrong’ path. He was not ready to talk about HIV/AIDS. He refused to fill the questionnaire. His eyesight was poor so he was not able to read. He wanted to go for an eye checkup because sight is very important for his profession. We asked him what was the reason for the spread of HIV/AIDS? He said the first reason is cinema and TV. It spoils both uneducated and educated public, whereas magazines and porn novels ruin only the educated.

He feels that invariably after watching such shows, most men visit brothels. Some of them bring a CSW to see the movie with them. We asked him what should be done to prevent this? He says the public and government should join together to ruthlessly ban all vulgar cinemas, magazines, novels, and serials; and also impose a fine on producers, editors and writers who indulge in producing such works. He says that whether a film is a box-office hit depends on how scantily dressed the heroine is, and how vulgar the dance sequences are.

He says it is impossible to walk in the streets of Chennai looking at the lewd posters. Majority of the posters carry skimpily clad actresses. He feels even the Old Stone Age culture was better, that is why they did not have HIV/AIDS.

He says more rural men suffer from HIV/AIDS than the urban men because they do not have safe sex. He adds that using a condom is not the only solution; it can only decrease the spread of HIV/AIDS. He lamented that Chennai has become a very cosmopolitan city and consequently, it has become easier for CSWs to practice their trade. Like the others, he was supportive of children affected by HIV/AIDS.

## 11

On 10th May 2002 we went to a police station to take collective interviews of policemen. They were practically frightened to see us even to receive the questionnaire. With no other option we left the station. But fortunately, a constable followed us. He wanted to get more information about HIV/AIDS. He said that he was drastically losing weight. He also felt very tired to even do his daily routine. He indirectly acknowledged that he had once visited a CSW and had unprotected sex. We advised him to go for a blood test as soon as possible. He was practically fear ridden; he feels he is infected by HIV/AIDS. He was very tensed. He did not talk about any awareness program for HIV/AIDS. He is vaguely aware of how HIV/AIDS spreads. He blames his present state on “one minute of carelessness”. He was very uncomfortable and in no proper frame of mind to talk about HIV/AIDS.

*Additional Observation:*

*Normally, we had heard that police used to take bribes from CSWs. It was a surprise to us that this police constable is not aware of how HIV/AIDS spreads and the problems related with it. It is clear that he fears the dreaded disease and the associated social stigma.*

## 12

We met a 22-year-old housewife Kanchana. She has studied up to 6th standard. She is a native of Vettaikara Kuppam, Kanchipuram. Her father is a farmer and mother is an agricultural coolie.

She says an HIV/AIDS patient must be taken care of by his relatives and friends, only then we can think of public support. These days even if a man has HIV/AIDS he never says anything to his family members. On the other hand, if a woman has HIV/AIDS, in majority of the cases her in-laws chase her out of her house if her husband is dead. Unless people tolerate their own kith and kin who are affected with HIV/AIDS and provide them with proper medical aid, how can we expect the public to be kind or sympathetic? She welcomes the idea of free HIV/AIDS for everybody and she feels that not everybody will come forward because of the stigma associated with the disease.

She says people are affected only due to lack of awareness. She is willing to come out and speak about HIV/AIDS so that people will become aware of it. She is not shy about spreading awareness and campaigning in her neighbourhood. She first came to know about HIV/AIDS through TV when she was just 15 years old. She says HIV/AIDS patients are very pitiable: first, they know they will die soon; second, they feel invariably guilty.

She feels HIV/AIDS patients should come forward and speak about their experience, only then it will have more effect on the public.

She says CSWs are not any new clan. They are very poor people who have taken up this profession due to their circumstances. She sympathizes with them.

She thinks it is not essential for a government employee to inform his/her HIV/AIDS status to anyone. It is ones own right, and they need not reveal their HIV/AIDS status to their employer.

## 13

We met a 22-year-old science graduate from Chennai. She is aware of the many ways in which HIV/AIDS spreads. She blames the rapid spread of HIV/AIDS due to lack of basic awareness. Even after knowing that they have some STD they continue to hide it from one and all due to guilt and shame. If they take proper treatment before the disease becomes chronic, certainly they could extend their span of life. Because of not taking treatment, they don't live more than 3 to 4 years after the infection. She feels that advertisements must portray real HIV/AIDS affected patients themselves and not models. Only then people will feel that anybody can get HIV/AIDS.

Everybody knows that it is easy to become an addict of bad habits or take up to bad ways than to remain good. She feels elders don't inculcate proper values in children. That is why children are irresponsible.

She says that a person is influenced not by reading or hearing but mainly by seeing (which is a subconscious level), so when he/she sees all eve-teasing scenes, rape scenes, sex scenes with vulgar song sequences, they will certainly ingrain that and try to imitate the same. She feels that a sensitized media will portray women not as objects but as living beings. This would reform youngsters.

**14**

We interviewed a 21-year-old student who is doing her graduate course in Siddha medicine. She is staying in a hostel in Chennai, her native place is Palayamkottai. Her father is a bank officer and her mother is a teacher. She says HIV/AIDS spreads because of multi-partner sex. She has talked about HIV/AIDS with friends and family members. HIV/AIDS patients suffer lots of mental problems that they don't discuss with anybody. If we shun and ill-treat them it would double their problems. If a medicine is invented to cure HIV/AIDS it will not spoil society because the good people will always remain good and only the persons who commit mistakes will continue to do so. Even if the government brings free blood test for HIV/AIDS it is doubtful whether all people will accept it. She feels that sex education can be given from the age of 13. Parents and not teachers should impart HIV/AIDS awareness. She feels government should ban porn movies and films. She has seen HIV/AIDS patients and she felt very sad when she saw them. If given an opportunity to work for awareness of HIV/AIDS she would willingly utilize it. She will certainly come forward to spread awareness about HIV/AIDS, now she has no time, as she is busy studying. She feels the number of infected persons is certainly not decreasing. The main reason people don't talk about HIV/AIDS is because they associate a social stigma since they feel that immorality causes the disease. She feels the steps taken by the government to prevent HIV/AIDS are not sufficient. She is supportive of a scheme of reservation for HIV/AIDS patients in the public sector, but she doubts if the patients will have the necessary stamina to do the work.

**15**

We interviewed 26-year-old student preparing for the civil services exam. She is unmarried. She resides in Anna Nagar in Chennai so that she can have easy access to books. She is a native of Mettur. She has not talked about HIV/AIDS with her family members because they would begin to doubt her; but she has discussed about AIDS with her friends. She feels HIV/AIDS infected children should be given good food and free medical aid. She supports free HIV/AIDS test but says that not many people will come forward to take it because taking the test is a sign that they have committed mistakes. She feels that the number of HIV/AIDS patients is only increasing day by day. The only reason is the ignorance about the disease in rural areas.

She learnt about HIV/AIDS at the age of 16 through her friends. She feels this disease cannot be cured by Siddha or Ayurvedic medicines. She does not know the difference between HIV/AIDS/VD/STD.

**16**

We interviewed a vegetable shop-owner in Chennai. He is 46-years-old and has studied up to the 10th standard. He earns about sixty rupees a day by selling vegetables. His marriage was an arranged one. He gives fifty rupees a day to his family and spends Rs.10/- on cigarette and tea. Due to his friends, he took to the habit of drinking and smoking at the age of 15. He has nothing to think about HIV/AIDS. He feels people see advertisements but forget about it immediately. He has learnt from advertisements that one should not go to CSWs. He thinks that by touching HIV/AIDS patients one may get the disease. He thinks this disease spreads through the air. He advices that one should be very hygienic, for instance one should never walk barefoot in the streets. For the past 10 years he has not seen any movie. Since every film is available on TV, there is no necessity for him to go to theatre. He feels prayer cannot cure any disease. Only medical means and not god can cure diseases. He says CSWs are the root cause of HIV/AIDS. He has heard about condoms on TV but has never used them.

**17**

We interviewed a 43-year-old graduate who runs a petty shop in Chennai. He is married. He is a native of Tiruchi who has lived in Chennai for the past 10 to 15 years. He says he has no bad habits. He has not seen any porn movie. If one is very careful, one will not get the disease. He feels visiting CSWs cannot be prevented because everyone cannot get a family life. So those who don't get family life visit CSWs, nothing is wrong in going to them but they should use condoms when they visit CSWs to save themselves from HIV/AIDS. He feels nothing can be said for sure about the spread of HIV/AIDS. He thinks that it can spread through air, water or even by touching the patient. He said that if a HIV/AIDS patient had TB we would get both TB and HIV/AIDS through breathing the same air. To control HIV/AIDS everyone should lead a proper life. Government cannot do anything to control it. Unless everyone is careful, no one can control the spread of the disease. Even if a man visits a CSW, only if she has HIV/AIDS he would get it, otherwise he would not get it. He says HIV is an initial symptom before one gets AIDS. He feels that it cannot be called a killer disease because it has a treatment.

*Additional Observation:*

*He has a major misconception that HIV/AIDS spreads by air. We could not convince him to change his opinion.*

## 18

We interviewed a 22-year-old salesgirl working in a shoe-company in Chennai. She has studied up to 7th standard in her native place, Tanjore. She says she first learnt about AIDS through advertisements on TV and radio. She says that HIV/AIDS spreads when a man/ woman has several sex partners.

She is unaware of other way by which HIV/AIDS spreads. She is very sure that it is not a contagious disease that spreads through casual contact. She says traveling in Chennai after nine in the night is very risk for women like her, who are originally from villages. She says one of the main reasons for the spread of AIDS is the fall in moral standards. She cites the example of male customers who have tried to misbehave with her. This type of impropriety will cause them to slip off into the path of sin and they will commit more and more mistakes. She says that people must always think before they act. They should not be impulsive and do something which they will regret all their life.

## 19

We met a middle-aged man who is self-employed as a mobile tailor. He is in his late forties. He has lived in Chennai for the past 25 years. Fifteen years ago, he had a tailoring shop in Adyar. It was difficult to pay the rent for that shop so he became a mobile tailor. He has roamed in Adyar, Indra Nagar, Thiruvanmiyur and Anna University and Saidapet. He says he has seen rich and poor customers. He usually spends his leisure time gossiping local politics with the people of that area. He earns over Rs.200 and on lean days just about Rs.20 to Rs.25, which is sufficient just for the day's expenditure. He earns enough to maintain his family.

We asked him his opinion about HIV/AIDS. He says the nation has become spoilt. All the time, movie songs are heard everywhere right from teashops to electrical repair shops. These songs have become part and parcel of everyone. It is unfortunate that the lyrics are vulgar. The main way to eradicate HIV/AIDS is to eradicate vulgar movies and vulgar item numbers. Next he wants a ban on prostitution. Now CSWs are hungry for customers. He wants the CSWs and the pimps who negotiate the customers to be put behind bars. The increase in AIDS patients has been relative to the increase in number of CSWs. He observed that CSWs catered to exclusive social strata: rich, middle class and the poor. The trade flourishes because of the weakness of men. Therefore it is difficult to totally control the spread of HIV/AIDS.

## 20

We interviewed a 27-year-old graduate named Wilson. He is a native of Madurai who is working in Chennai in the commando squad. He says that HIV/AIDS is a disease of the scientific age; it was invented to control man from loose conduct. In the 21st century, the world has shrunk due to advancement and hence it has become victim of several evils. He says HIV was first discovered in 1981. He says HIV is not AIDS, and only the advanced state of HIV is AIDS. It is a killer disease that spreads in India mainly through sex.

A complete cure has not been found for AIDS till date. He says HIV/AIDS is more prevalent in developing and under-developed nations. Developed nations know to control it and stop the spread. He occasionally smokes and drinks. He got these habits ‘automatically’ once he was 17 years old. He has seen porn magazines and films. He feels 75% of the people live in villages and they do not know much about HIV/AIDS.

Since STD and HIV spreads via sex, people cannot talk about it openly due to shyness and social stigma.

## 21

We interviewed a 38-year-old man from Korukkupettai who works in a upholstery company. He has studied up to 10th standard. When he was asked to give his idea of how people got HIV/AIDS, he said it was through women: CSWs and multiple sex partners. He came to know about the disease only after he met two people who were affected with HIV/AIDS.

We asked him why people shy away from talking about HIV/AIDS but are willing to talk about diseases like diabetes or cancer. He feels this is a sexually transmitted disease that cannot be talked about publicly.

He believes that this disease can be cured by faith in god. He cannot say anything about it since everything is in the hands of god. When he sees advertisements about HIV/AIDS he feels the need to be very careful. Whether we watch AIDS awareness advertisements with the family or alone, we have to be cautious! He says that when we talk about AIDS with a person living with HIV/AIDS, he/she does not reply, he/she only laughs. Government repeat the same advertisements—only if it is changed from time to time, people will listen to it. Otherwise they will consider it to be a routine matter and nothing will get into their mind about HIV/AIDS.

He has never heard about CSWS. Also he has not heard about condoms. It was an unforgettable experience when he came to know through his friend that some teenage girls in Andhra Pradesh were having HIV/AIDS. Those girls were not even married, so he was shocked to know this story. He has not seen these girls, and he was sad when he heard their stories from a friend. Another direct experience when he was involved with HIV/AIDS was when he had to admit a close relative to the Tambaram Sanatorium.

## 22

We interviewed a 22-year-old man named Mahesh. He has studied up to 10th standard and is unmarried. He is a native of Madurai, he went to Paramakudi in search of a job, stayed there for four years and has now come to Chennai.

He thinks he will never get HIV/AIDS. According to him sex is the only cause for the spread of HIV/AIDS. To save one self from HIV/AIDS one should have control of feelings and abstain from sex with strangers. Once in a while he consumes beer. He does not know the full form of HIV/AIDS. He thinks that HIV/AIDS might spread by using things like plates, tumblers, etc. that were used by someone who already had the disease. We asked him how he would react if he saw a person affected with HIV/AIDS. He answered that he will be kind but he would also maintain a safe distance. He is not afraid to discuss HIV/AIDS with people who are not infected with it.

## 23

We interviewed a 32-year-old police constable named Saravana Kumar. He has studied up to 12th standard. He gave a very long interview. He said that getting the disease was like falling into a deep sewage. The government should take more steps to help the HIV/AIDS affected persons. Only self-control can save somebody from HIV/AIDS because there is no medicine to cure it.

Government should also soon discover medicine for HIV/AIDS. Even discovery of medicine can help only when people can control their feelings. Without self-control

HIV/AIDS cannot be eradicated. It is impossible to fully eradicate HIV/AIDS. It can be controlled to a certain extent.

The policeman said that people are actually scared of the disease and the social stigma it holds. He wondered why people like us work in a HIV/AIDS project. To remove such fear and adverse feeling towards HIV/AIDS, people must conduct seminars and debates about the spread and control of HIV/AIDS.

When someone says that Tamil Nadu is the state that is most affected by HIV/AIDS, does it not make us feel bad because we are from this state? Likewise, if someone says India is the worst affected country with HIV/AIDS will not Pakistanis think very bad about us? So we should have love for our people, our state and our country. Only then we will think of HIV/AIDS not as a problem of others, but as our own problem.

We asked his opinion on introducing sex education in schools? He said "only in India, people make a big fuss to talk about sex. It is good to give sex education to matured children. In our villages, if the sari slips a little and the blouse is seen people look with a craze. For in our nation, there is no openness about sex. So it has become difficult for men to control their urges. That is why they think of sex as a very major thing. That is why sex cannot be talked in a coeducational classroom. People don't have mental maturity and self-control."

We asked him about his friends. "I don't have bad company. I don't have good friends. Each is a different type. Some don't like us giving advice. Some think that anyone who takes good advice is a fool. We can't say anything about them." He was a social drinker.

How did you start to drink? "Just when I finished twelfth standard I went to work. At that time, I attended my friend's marriage and drank beer. Later on, whenever we went out sightseeing we drank beer. But I don't drown myself in drink."

Have you seen porn magazines? "I have only heard about it. One should not go to such a cheap extent as to see it."

Why don't people talk about HIV/AIDS as they do about cancer or diabetes? "Yes, because it holds a social stigma as it spreads through sex. It is not like cancer/TB that anyone (even with a good character) can get. In Tamil Nadu, HIV/AIDS spreads mainly by sex. We don't get HIV/AIDS because we talk or eat or sit or kiss an infected person. We get HIV/AIDS only when we have sex with a man/woman who have HIV/AIDS. Otherwise there is no chance of getting it. It does not spread through mosquito bites. There is no history of anyone having got HIV/AIDS through mosquitoes.

Do you know the full form of AIDS? "It is not my job to know the full form. It is the job of doctors and persons like you who come to take surveys."

What do you think about CSWs? "They take up this profession only due to their circumstances. No one ever enters this job of her own will. It is only a means of income. Family situation is also a cause for it. Suppose one wants money and he is in absolute poverty, he may steal, but he does not take up stealing as a profession. He who takes it up as a profession is a scoundrel. Some people live below the poverty line and have taken up this trade. Some rich persons do this trade only for the money and not due to poverty."

What can you say about the influence of police department in the trade of CSWs? "When they raided a flesh trade center these women run, I was part of the team. We usually advice them not to practice the trade. If they are listening we leave them. But if they continue the same lifestyle we register cases against them."

Suppose you meet a HIV/AIDS patient, how will you conduct yourself? "In my opinion we should not neglect or hate him. HIV is in him and it is in his blood. By being with him or eating with him or sharing clothes with him or by talking with him we are not going to get HIV/AIDS. Suppose we segregate him from society it would only make him more depressed and make him feel as if he were a criminal. By doing so a man who is to die after 60 years will die at 30 years, so we must help these patients to the maximum possible extent."

## 24

We interviewed a security guard named Duraikannu who works for a private company. He is 68-years-old and has studied up to $12^{th}$ standard. He is a native of Maduranthakam. His monthly salary is Rs.1500/-. He stays in a nearby place with his wife. He has two children, one son and a daughter. His wife is employed in the midday meals scheme.

What do you know about AIDS? “It is a cruel disease and one should be very careful not to get it. I do not know how it spreads. I do not know the symptoms of the disease. I very badly want to know the symptoms. Now that you are here, please tell me about it.”

Do you have any form of relation with any other women other than your wife? “Thirty years ago I had such connections.” Was she a CSW? “No, she was a family woman.”

Was it before marriage? “It was only after marriage.”

Did you feel guilty when you have sex with your wife after you had relationship with other women? “Yes I felt bad. It was many years ago.” Did you have such relationship with many women or only one woman? “Only with one woman.”

Did your wife know about your relationship? “Yes, she scolded me. After that I did not go after the other woman.”

Why are you frightened of AIDS? “Suppose I have become infected by HIV/AIDS from someone who visits the building? Will it spread by air? I am very frightened of it. Do you know how to control the disease? Please tell me. I don’t know anything about the disease.

Have you spoken about HIV/AIDS with your family members? “They know very well about it. My daughter says this is a very cruel disease. So my children will not do any wrong. They will not love. The marriages of my children are only arranged marriages. When our children are so good, I regret that I made mistakes in the past.”

What do you feel about a cure for AIDS? “On one side, we are growing in technology. On the other side, certain evils are also growing more than technology. That is why we don’t have medicine to cure HIV/AIDS!”

## 25

We interviewed a 52-year-old porter. He has studied up to 5th standard. He did not give the name of his village. He says HIV/AIDS is a disease that affects anybody (rich or poor) when they visit CSWs. He argued, “Doctors say that smoking causes cancer. Drinking causes damage. But more advertisements must be given stating that visiting CSWs causes AIDS.”

He further said, “if you use condoms you will not get infected. But the quality of the condoms is so poor that one gets infected.” If he gets a chance to meet HIV/AIDS patients he would be happy to be with them. He recollects how in a movie Dr. Siva used to touch the leprosy patients. One should do any action with devotion. These days nobody makes movies like that. In any film of MGR, there will not be a scene where he smokes but nowadays every hero smokes onscreen.

If movies of popular actors include scenes in which HIV/AIDS patients are treated with kindness certainly the stigma associated with it would go.

## 26

We interviewed 32-year-old Sunil who worked in a fancy stores in Chennai. He has studied up to 10th standard. He is working in this job for the past 2½ years. Earlier he was a driver. That job was very tiring and he did not like the company of drivers. So he has taken up this new profession.

When he was young he had all bad habits. In fact all his driver friends had the same kind of habits. Now he has reformed that is why he changed his profession. He feels the awareness program given by the government for HIV/AIDS is more than sufficient. He was reluctant to speak more.

## 27

We interviewed a 45-year-old man who is self-employed as an electronic goods repairer for the past 25 years. He says to control the spread of HIV/AIDS we should try to change our way of life, culture and tradition in order to suit the present context so that people are careful in their sexual behaviour. If one has to change society, one has to change the media. But media and technology is only in the hands of the ‘higher’ castes. That is why a woman is only identified as an object of pleasure. Likewise, TV and newspaper advertisements project half clad women so that their goods have better sales. This sort of advertisement will only make people go for sexual activities. Because of sexual urge and need, men take up to a wrong path in life. Especially drivers, cleaners, and businessmen who go to other states or cities become easy victims of HIV/AIDS.

A week ago he had been to Delhi. "I had to restrain myself from visiting sex workers. Everybody who goes to other cities visits red-light areas. But I tried to follow the correct path. If people are to take up a proper course of life it is only in the hands of media. In olden days, courtesans were available for Kings and ministers. They were a special type of call girls who danced to please the courtiers and satisfied the king and ministers. But nowadays, due to economic problems, women take up prostitution to support themselves. With the aid of computers, prostitution is carried out through Internet. This has affected the younger generation very badly. "Proper economic planning by the government will reduce poverty. This in turn will stop women from entering prostitution. This will directly reduce the number of HIV/AIDS affected persons. Indirectly the government supports vulgarity in cinema because it is not interested in moral or social reformation. The media gives only very little information which is good for reformation of society or revolution. That is why the media and government go united in these things. If government decides to change the media it can certainly do so. But they don't do. For instance in Kerala, they put a rule that MTV should not enter Kerala. On the other hand Tamil Nadu has encouraged MTV. Globalization and liberalization has made society very poisonous. The attitude that AIDS is somebody's problem must be changed. Unless that change takes place it is not an easy task to control the spread of HIV/AIDS. One cannot easily dispose the problem of HIV/AIDS. Everyone should show deep concern over it. Unconcern will not help anyone, since AIDS can affect anyone at some time or the other. It is a people's problem."

## 28

We met a 19-year-old boy Karthik who has studied up to 12$^{th}$ standard. His only past time is going out with friends. He knows about porn films, but he has never seen one. He knows that HIV/AIDS spreads through sex with strangers (CSWs). He says one will first get only VD, only then AIDS. It is in the hands of each and every individual to control HIV/AIDS. If one has self-control certainly one will not get the disease.

He says it is the right of every individual to think anything about HIV/AIDS. But he will be very friendly with a HIV/AIDS patient, if he meets them.

He feels sex education can be given from 6th standard onwards. He has no faith in god. "If God can cure the disease, he could have first cured immorality."

## 29

We interviewed Senthil, a 26-year-old computer graphics artist who has studied up to 10th standard. He has no past time or hobbies because he has no leisure. He has seen porn movies with his friends. He is a bachelor. He used to be a chain smoker and alcoholic but now he has given up those habits. He feels HIV/AIDS is a 'bad' disease. "We pay money simply to get the disease. To control HIV/AIDS the only way is to avoid sex with CSWs. If no one goes to CSWs, automatically their number will decrease." What do you feel about CSWs? "I do not know about CSWs. Only I should visit CSWs and then find out information about their way of life."

He refused to talk more about HIV/AIDS.

## 30

We interviewed a 24-year-old man who has studied up to 10th standard. He owned a watch repair shop. He is married and has a daughter. He thinks of HIV/AIDS as a killer disease. "To prevent HIV/AIDS one should not covet another's wife. If one goes to CSWs, condoms should be used. Above all being faithful to ones partner prevents HIV/AIDS."

He claims to have led a proper life, has never been to a CSW or been with any women other than his wife. He goes for movies once a month and views only ordinary movies. He has never read porn magazines. To know about sex, one need not read such books. He says one cannot ban such books because ban of porn books cannot cure or reduce the disease. Actually such books help people control their emotions and sexual urges.

He thinks CSWs exist only due to horrible circumstances of life. Many young girls are kidnapped or cheated at a young age and later made to work as CSWs. Because of this, the population of CSWs is steadily increasing. No women come to this profession of her own will.

## 31

We met an 18-year-old boy working in a fancy shop that sells handbags and suitcases. He had studied up to 8th standard and is originally from Kerala. He has been settled in Chennai for the past ten years. To be free from HIV/AIDS, one should be very conscious of ones actions. By seeking CSWs one will certainly get HIV/AIDS.

He has seen porn movies and has read porn magazines. He philosophically says, "I don't want to lie or hide the truth, this is the weakness of youth." He has so far not had sex and he is sure he will remain that way till he gets married. He has seen condoms but does not know how to use it. "I have not told my parents about seeing such porn movies. In two years I am going to get married. So why make such small issues big? It is natural. If one does not feel so, he is not a man!"

## 32

We interviewed a 59-year-old man named Krishnamoorthy. He has studied up to 10th standard and works in the Chennai Corporation. He is married and has two children. He confesses that he smokes, drinks and visits CSWs. He started doing this kind of stuff at the age of 18. He got all these habits only through his friends. We asked him what pleasure he found in all these? He said that he healed his wounded soul by fire. He has visited CSWs before and after marriage. "I fear about HIV/AIDS. I had seen VD in one or two ladies. I too was affected with VD and I took treatment for it in GH. I have fear when I visit CSWs. I have also had sex with other married women. I always pay them. Sometimes it is as cheap as Rs.50, sometimes it is up to Rs.300. Suppose we ask them why they do this job they will immediately ask back, are you going to support me or provide me with food?"

He blames the current generation. "Today school children are falling in love. In my days, even college students did not dare to do such things."

**33**

We interviewed a 24-year-old lottery shop and STD booth owner. He is a bachelor who has done various jobs. He thinks HIV/AIDS is a 'killer' disease because no one has survived it. He blames ignorance as the reason for the increasing spread of HIV/AIDS. He feels posters and wall writing will only reach out to the educated. But in reality, only the uneducated people have to be targeted. To spread awareness among them there should be advertisements on television and radio. He thinks that not many people in the villages know that condoms can be used to prevent HIV/AIDS, they only think that it is a method of birth control. He feels that nobody discusses HIV/AIDS in public because the topic of sex is never openly discussed in Indian society. "Even though men visit CSWs in private, they will only condemn prostitution in public. This kind of dual attitude and hypocrisy is responsible for the spread of HIV/AIDS in our state."

**34**

We interviewed a 35-year-old cassette shop owner named Krishnakumar. He has been doing this business for the past 15 years. He says HIV/AIDS spreads due to multi-partner sex, blood transmission and sharing injection needles. He thinks HIV and AIDS are two different diseases caused by similar germs. He said, "It does not matter if you have sex with an infected person once or ten times. A crime is a crime. So you will always land up getting the disease." He thinks that people avoid talking about HIV/AIDS because of two reasons: first, they fear the disease; second, they want to hide their own faults.

He is supportive of the idea of free treatment to HIV/AIDS patients. We asked him what he thinks of sex education in schools. "He says it is good to give sex education from $8^{th}$ or $9^{th}$ standard. But we never know how students will take it. If they take it positively, it is well and good. If they take it negatively and try to experiment, it will only lead to ruin." He feels the awareness advertisements are not natural. They appear artificial in certain aspects.

If a particular area is more affected by HIV/AIDS one should first study the environment, like presence of brothels etc. He criticized the advertisements that promote condom usage to protect oneself from HIV/AIDS. "This advertisement supports their faults and provides them a means to protect themselves from HIV/AIDS. This will make people think in the same way. Such advertisement by government is very morally flawed because it amounts to saying that one can steal but do not get caught by the police. You can copy in the exams but make sure the invigilators or supervisors do not notice you. Such advertisements must be banned." He does not know how long it takes for a HIV-infected person to become an AIDS patient. If he meets an HIV/AIDS patient he would behave as usual since this disease does not spread by touch. He would not fear or shun them.

He says women never purchase porn movie cassettes. Only men in the age group 20-30 buy these porn CDs or videocassettes. He is of a varied opinion about CSWs. In any social setup CSWs are also essential. In fact they do a service to the nation. Otherwise infatuated men would be misusing family/ married women by force and this would certainly ruin the society.

## 35

A 24-year-old Keralite businessman, settled in Chennai was interviewed. He is unmarried. He said that in the world, India is leading in HIV/AIDS and his state Kerala has the maximum number of HIV/AIDS patients in comparison with its population. The major reason for the spread of HIV/AIDS in Kerala is that men leave their wife and go with another man's wife. Each and everyone must reform himself if AIDS is to be eradicated. If a single person thinks about it nothing can happen. In recent days, because of cell-phone and e-mail the flesh trade in Chennai is carried out in a very modernized way. HIV/AIDS affected persons are segregated even by their own relatives. Even though this disease is not contagious, people shun away from it because it is associated with social stigma.

He says by talking with a HIV/AIDS affected person we are not losing anything.

Any one who sees a HIV/AIDS patient in the right perspective would never commit a mistake. He is of the opinion that everyone should know about HIV/AIDS. There is nothing wrong in seeing the AIDS awareness advertisements with ones family. When people can see a movie together, they can surely see advertisements together without feeling embarrassed.

## 36

We interviewed a 26-year-old man, owner of a garment shop. He was very unwilling to give the interview. He said that HIV/AIDS is something CSWs should worry about, it is not a topic for shopkeepers like him. We told him that even CSWs are human beings and we should not condemn them. He said, "A CSW can never be rehabilitated. Even if we set her up as a florist or vegetable vendor, what will she say about her illicit children? It is better they continue in the flesh trade itself. They earn through the heinous method of selling their bodies. They have to work that way to survive."

He says that there is a disease called HIV/AIDS but he does not know much about it. He is married and has a son. He has not visited CSWs even before marriage. "I was leading my life like Swami Vivekananda" he says very proudly of his celibacy. He has seen porn movies only after marriage. He says everyone has certain concepts about sex and their opinions are diverse and different. Even if someone is a close friend one cannot discuss sex with him or her. He says, "Having sex with a stranger will not make us HIV/AIDS infected. Only if we repeatedly have sex, we might get the infection. That too, if the woman is free of any disease, we will not get any disease. If she has HIV/AIDS, you are sure to get it. But why should we seek CSWs, when we have our wives?"

**37**

We met Senthil, a 23-year-old man working in the fish market. He has studied up to 7th standard. He says he has seen a porn movie every week for the past five years. Only his friends used to take him to these movies in the beginning but now he goes on his own. When he has sex with women known to him, he does not use condom. But with women who are strangers he has always used condoms. He has usually spent Rs.150 on these CSWs. They are his good friends. They confide all their problems to him.

**38**

We interviewed Prabhakar, a 56-year-old motor mechanic. He is staying in a rented house paying Rs.1200/- p.m. He is married and has two children. He smokes and at times drinks beer. On social occasions he drinks rum/ whisky. He took to all these habits at 18 years. He says, "AIDS can be eradicated only through self-control. Only if thieves control themselves, theft can be stopped. This is a similar situation." "As long as CSWs exist, HIV/AIDS will continue. Nowadays flesh trade takes place in Tata Sumo and Cielo cars. Prostitution is at its peak. Earlier, women became CSWs due to poverty. Nowadays they become CSWs to earn fast money and become rich easily. Prostitution will not differentiate between rich or poor. Both will certainly get HIV/AIDS. The police department must eradicate the flesh trade. We see in the newspapers the names of actresses and extras being flashed for doing prostitution. Are they so poor to do this trade? Only greed makes them do all this sin." People don't talk about HIV/AIDS even to their doctors. They feel shy when they have to talk about their behaviour. They keep visiting CSWs and prolong the disease till it comes to a very chronic stage. To decrease the number of CSWs, "the government must first solve the unemployment problem. Unless unemployment is properly tackled it is impossible to decrease the number of CSWs. The government brought the CSWs from Bombay and tried to rehabilitate them but could not succeed due to the problem of finding them a suitable job. Women do prostitution in order to educate their children; some women do it because their husbands have deserted them."

**39**

We interviewed Murugavel from Ariyalur who works in a teashop. He is 29 years old and has been settled in Chennai for the past 16 years. He is married, and has a son and a daughter. He says one gets affected by HIV/AIDS only by visiting CSWs. He has not visited CSWs, but he has seen porn movies. He says he has seen them before marriage for enjoyment. He has no habits other than drinking tea often.

When we asked him if AIDS prevention is possible, he laughed very loudly and said it was impossible. He feels that a person's mindset is the factor, which determines whether they will mingle with a HIV/AIDS patient, or not. He is very sure that HIV/AIDS is incurable. Even gods cannot cure it. He says "I have not seen god and I know nothing about him. I have only seen people keep variety of statues made of stones and worship it."

**40**

We interviewed a pavement vendor Suresh. He has studied up to 9th standard. His native place is Tirunelveli but he has been living in Chennai for the past 30 years. Earlier he was selling plantains and he got into lots of debt, so he has changed his business.

He thinks HIV/AIDS cannot be eradicated. No one should get this disease because nobody does the mistake consciously. He has one major doubt: "Prostitution has existed for so many centuries. But why do people get HIV/AIDS only now? So, it does not come from CSWs who have always existed in our land. It is only a mysterious curse on the world."

He wants more publicity to be given to the symptoms of HIV/AIDS. He said, “If we know the symptoms we can go and get ourselves tested. If diagnosed in the earlier stages, the disease can be controlled. But who in the world knows what the symptoms are?”

“Today prostitution takes place in star hotels and on the platforms. Those who need a CSW know where to find one. But the majority of people are ignorant of it! Even today if one goes to Parry’s corner (Flower bazaar) in the night, one can see police arresting several women. Even to come out on bail, they don’t have people to arrange the money. Only the pimps have to take them out. Thus they are constantly in servitude.” Only by finding a proper solution to the CSWs problem, one can find methods to control HIV/AIDS.

“Mixing one drop of lemon in milk will curdle it. Likewise having sex with an AIDS affected CSW is enough to get HIV virus. But nothing will happen if we move with patients.”

He says our society is very conservative; this is one of the causes for HIV/AIDS. “Even a wife and husband do not discuss about sex openly among themselves. That is the state of affairs. How will they ever discuss sex with others? Another factor is, if a man holds hands with a woman (even his wife) and walks in the roads, people look at them scandalously; this is another reason for viewing HIV/AIDS as a stigmatized disease. Even in buses, men don’t think of women as sister or mother or as an equal, that is why there are ladies seats!”

He says sex education can be given in the plus two classes. He says from his 18 years he has the habit of smoking and drinking. He says he has seen a porn movie only once. The story looked very imaginary and not life-like so he never again went to such movies. He says people think of HIV/AIDS only when they are infected by it. Otherwise, they just ignore the problems. This lackadaisical attitude must change.

## 41

We met a 38-year-old government employee named Babu. He is a graduate who resides in Medavakkam. He is married and has a son and a daughter. He says to be free from HIV/AIDS self-control is important. Even foreign countries that have good healthcare and wealthy nations are not able to control HIV/AIDS.

He suggested that to spread awareness about HIV/AIDS, government officials, students and NGOs should camp in different parts of the nation and educate the rural people about the disease. According to him the steps taken by the government to stop the spread of HIV/AIDS and the awareness programs are not sufficient. He feels awareness about HIV/AIDS in the society especially in the villages is insufficient.

He feels that all adult education centers especially in rural areas can be converted as HIV/AIDS awareness centers. These adults in turn can teach the youngsters who are unaware of HIV/AIDS. People know about the seriousness of the disease and its evil effects. Certainly it is not like cancer or diabetes. He says there must be a law that everybody must undergo a free HIV/AIDS test. He feels that persons who say ‘HIV/AIDS is not my problem’ suffer from ‘complex.’ It only implies they do not have awareness that HIV/AIDS is not a contagious disease.

He thinks the symptoms of HIV/AIDS is that one feels tired and does not have appetite etc. He says he is not aware of the exact symptoms. He says, “it is impossible to cure of HIV/AIDS by mere prayer. Faith in god is only a means to make people respect each other. It cannot bring about cures.”

### 42

We interviewed, 25-year-old Das who works in an auto company. He earns Rs.1500/- per month. He says the money he earns is just enough for his boarding and lodging. He says he sees movies, listens to music and reads books and helps others in need. He does not see porn movies or read porn magazines because he feels his health would be affected and he would soon become a HIV/AIDS victim. He feels that only men with no self-control seek CSWs. It is more important to control ones senses than control ones budget. He says HIV/AIDS is not a contagious disease like TB or jaundice. One does not become infected even if they move with an affected patient. In India, particularly in TN, HIV/AIDS infection spreads only through sex. If the husband has HIV/AIDS, the wife and young children have chances of getting it. Likewise if a wife has HIV/AIDS, there is chance of her husband and her child getting it. As soon as anybody comes to know that he/she is affected with HIV/AIDS, they should take treatment immediately and there is nothing to feel ashamed or shy about it. Ignorance is the root cause for the spread of HIV/AIDS in this modernized world. If one is infected with HIV/AIDS we should treat him/ her kindly.

He has seen porn movies from the age of 16. He used to go along with his friends to such movies. He says he did not know that they were going for such a movie but by seeing it he understood what is the institution of marriage. He says he got many opportunities to visit CSWs but he controlled his feelings and did not go for sex. People have ill-feeling about HIV/AIDS because this disease has a lot of social stigma associated with it. He thinks it takes one year for the HIV virus to become AIDS, he has seen this information in a notice given by NGOs at the Egmore railway station. In his school sex education was given in the 6th standard. He says only the fear that AIDS is a killer disease prevents people from talking about it. He says he would move with HIV/AIDS patients in such a way that they don't become affected by his actions.

He says most of the CSWs have taken up this trade only due to poverty. Government can give them alternative employment.

Since most of the HIV/AIDS patients are weak and some are terminally ill there is nothing wrong in giving them the first preference in treatment. He feels AIDS may be cured if there is god's grace. He has not talked about HIV/AIDS with his family members. When he watches the HIV/AIDS advertisements alone, he would think that no one should get this disease. But when he is with his family, he would wonder why such advertisements are put.

## 43

We met a 30-year-old man working as a mechanic. He has studied up to 7th standard. He has been doing this job for the past 13 years. He has gone abroad and visited Abu Dhabi three times. His marriage was an arranged one. He has two children and he belongs to a middle class family. He says he earns from Rs.3500/- to Rs.4500/- per month. He gives his full income to the family. He has seen advertisements about HIV/AIDS on TV and radio. In his opinion, there is a possibility of spread of HIV/AIDS through mosquito bite. He concluded by saying that everyone should lead a clean and controlled life to avoid the spread of HIV/AIDS.

## 44

We interviewed a rice merchant. He has studied up to 5th standard. He is married. He is in his late forties. He has only one son. His wife has studied up to the 10th standard but is only a housewife. He says that 90% of the people in Chennai are aware of HIV/AIDS and are very careful to protect themselves from HIV/AIDS. But people in Delhi or Bombay are not so careful. Also, people from Tamil Nadu who visit Delhi and Bombay may visit CSWs there and catch HIV/AIDS. He says that uneducated villagers who visit Chennai for pleasure trips might acquire HIV due to ignorance if they are not careful. Many migrant labourers who work in construction sector in Chennai, do not take home money alone with them, they take back HIV/AIDS also. He says that men visit CSWs only because it is impossible for them to bear the separation from their families. “HIV/AIDS affected persons hold their life in their hands. I do not know the other means by which this disease spreads. I am not educated enough to know about these things.”

“There is nobody in the whole world who has not seen a porn movie. I have seen them when I was 18 years old. I have seen lots of these films. However, I have not read porn books.” He also acknowledged that he has all habits like drinking, smoking and taking drugs like ganja. He has visited CSWs and only his friends introduced him to CSWs. He says he does not wish to talk about his visit of CSWs.

## 45

We interviewed a man in his late twenties who works as a construction labourer. He had a love marriage. He stays in his own house. He does not see porn movies; he only sees movies that feature Tamil actors. He was very offended by our question. When he lived in a different surrounding before marriage, he used to see such movies. He says HIV/AIDS advertisements given on TV are not very clear. They only say in a roundabout manner "Why go for sex when you have your wife?" or "Have only one partner" or "Use condom when you visit CSWs."

He says advertisements are not complete! Only if it is pointblank, people like him will follow things properly. "Do the police not know about CSWs? Police are hand in glove with them. If you have courage, go to a red light area and take an interview. You would be torn into pieces!" he added.

He said, "We have not made any mistakes to get HIV/AIDS. So how can we say anything about HIV/AIDS? Does one get HIV/AIDS only due to women? Is it not due to change in food habits? When people don't have resistance to disease, does it not imply that they have no good food to eat? Earlier, the poor people ate beef and were free of HIV/AIDS; they were strong and sturdy. The BJP government banned cow slaughter, and people started getting HIV/AIDS! Let the government provide a balanced diet to everybody in the nation. There will be not a single case of HIV/AIDS."

## 46

We interviewed a spare-part salesman from Chennai. He is a bachelor. He earns Rs.2600 per month of which Rs.1500 goes towards the rent. He gives the rest of the money to his mother who takes care of him. He says he does not have any money left to himself, so he cannot afford CSWs. He admits that he has seen porn movies. Since these movies are not screened in the theatres, he used to get the CD and view it when nobody else was at home. He thinks it is perfectly natural to see porn movies. He also knows about condoms. "I have seen a big model of a condom near Padmanabha theatre in North Chennai." He has not discussed HIV/AIDS with his friends. He says many HIV/AIDS advertisements are given only in Pothigai TV and not on other channels.

## 47

We interviewed a clerk from the transport department. He is a 29-year-old graduate. He is married and has one son. He said at the outset that he didn't want to give us an interview. When we persisted, he relented and started talking to us.

He said that HIV/AIDS couldn't be prevented even by condom usage. He has not visited CSWs after he got married. "Because I am a married man, I have family responsibilities. It is unnecessary for me to go and visit CSWs. Some people who have money think that they can buy anything in the world. That is why they try to buy sex."

He has heard of a family in Kerala where a boy and a girl (They are brother and sister) got HIV/AIDS from their parents. The children's grandfather took pity on these children and educated them. The school authorities sent back the children once they came to know that these children are HIV/AIDS infected. So he admitted them in a government school. Even in the government school they were given separate seat. "I think India needs Marxism. At that time, Indian society needed Ambedkar and now it needs Marx. We were divided on caste lines, we are divided on class lines. This must go."

He says even advertisements are not going to help the people. "They will not see the HIV/AIDS advertisement but they will excitedly watch cricket matches played between India and Pakistan. Why do they not take the same interest and see these advertisements. To be more sensible they should talk about HIV/AIDS with the family. They don't talk about HIV/AIDS with family."

He says that he enjoys the Vijay TV program "Kathai Alla Nijam" which has often featured CSWs.

He says since his family is an educated one if they see advertisement about HIV/AIDS they change the channel.

## 48

We interviewed a 22-year-old man named Senthil. He is educated up to the $10^{th}$ standard. He works as a daily-wager in a teashop. He is a bachelor and not interested in marriage. He has heard about HIV/AIDS, but does not know much about it. He is not interested to know about it either. He says only when we go to CSWs we get HIV/AIDS.

He says he has seen porn movies. He says there is no chance that belief in god will cure the patients, only medicine can cure HIV/AIDS. But scientists will never reveal medicine to AIDS because it will indirectly approve free sex. So the patients will always continue to suffer. He says he has not discussed about HIV/AIDS with his friends.

He has also not seen HIV/AIDS advertisements in the presence of his family members. If such TV advertisements come in the presence of his father he would walk out of the room. He says he is very disinterested in life.

## 49

We interviewed a 23-year-old man who repairs watches in Chennai. He has studied up to 9th standard and is doing this job for the past 10 years. He is unmarried. He occasionally sees porn movies. He says that the Americans have done research about HIV/AIDS and found that it has no cure. To protect oneself from HIV/AIDS one should use condoms. He admits that he had sex with CSWs and used condoms. CSWs need money and only their poverty forces them to take up this trade. He wants a ban on the practice. He has often discussed about HIV/AIDS with his friends.

God will cure all patients, even AIDS patients, because it is his duty to protect his devotees. He says that nobody disrespects the HIV/AIDS patients. Only they shun away from others because of guilt. They might also leave us because of reduced immunity they will get our diseases very easily. One of his neighbours died of HIV/AIDS so he can never forget the experience. It was frightening to see him suffer with the disease, especially during his last stages.

## 50

We interviewed a 56-year-old man named Devaraj working as a construction labourer. He has studied up to 4th standard. He has four children: three girls and one boy. HIV/AIDS is not only sexually transmitted but it comes because of our change of food habits. The youngsters think that the disease spreads only through sex. Youth of today lack stamina because of changed food habits. Protein-rich food would certainly save people from HIV/AIDS.

People have no responsibility these days—they spend Rs.20,000 to Rs.30,000 for a HIV/AIDS test itself. They keep the matter a secret from others. They don't like to go to government hospitals to take the test, because of laziness and pride. He has never even looked at any woman other than his wife.

He says that one can get to know about many things by reading books, but one has to analyze the right and wrong.

Hospitals don't give proper counseling and follow-up to patients who are diagnosed to be HIV+. He knows two college students who got themselves tested, the result showed they had HIV sero-prevalence. So these boys quit the college, ran away from their homes. Because they felt that they will not live long, they went on visiting CSWs to have all the fun they could. He says now no one knows how many CSWS they will be infecting and these CSWs would be infecting in turn several men. If we know that someone might have HIV/AIDS we should take him to a doctor first give him proper councilling and then give him medicine and treatment. This should be every Indian's duty. He feels the awareness among people is not sufficient. He says the only means to spread awareness about HIV/AIDS is cinema. Nowadays the medical profession has become a business. That is why doctors don't say enough about any disease to the patient since it would reduce their income. In those days doctors gave medicine for cure and were interested in curing their patients. Nowadays they don't bother about the real disease or cure; they develop a patient as a customer, so in the end most of the patients who are clients become chronic patients. If the doctor runs a hospital, he admits him in his own hospital. Thus if they get a rich customer they never leave them. That is they never cure them. He has lots to say about CSWs. If a migrant labourer goes to Calcutta or Delhi or Bombay and visits CSWs it is mainly out of loneliness and confidentiality. When he is totally away from the family no one will ever know what he is doing. That is why most of the poor labourers get the disease after they visit other states. But if they visit CSWs in Chennai it shows their ability to pay for sex. "Today, if a man gets HIV/AIDS he suffers both internal and external conflict. The problem he faces from public is that he cannot face anyone. If there were a social service center like the one run by Mother Theresa, all problems of HIV/AIDS patients would vanish the same day. But today's NGO work for money and there is very little service motive in them."

He went on to talk about the food habits. "Earlier the food had nutrient values. This is because the manure used was natural and not chemical. So, fertility of the land was preserved. Now the land has lost its fertility. So the modern crops do not have true food value. These crops grow within three months. But in those days the crops grew for 6 months and had nutrient values. So people in those days had strength and stamina. Everything has been ruined by the chemical pollution. Our medicine and manure had good effect on humans and land; we were completely cured of our diseases by our natural medicine. But, in allopathic medicine if the doctor cures the liver ailment, the same patient develops kidney ailment, then by trying to cure the kidney problem the patient becomes a heart patient. In majority of the cases no disease is curable. We can live with the disease, so they are livable diseases. If we can think and avoid taking all these modern medicines no disease will come near us. We will improve and above all our nation will improve."

## 51

We interviewed a diploma holder from Chengai who is currently employed in Chennai. He said that right from his birth he has never done any mistake. "I live with respect and dignity in my village and I do not want to spoil that repute and good name. But I have heard about HIV/AIDS."

He started to berate men who visited CSWs. "If any man visits CSWs, he can never be corrected. He will shift from one women to another women. He will not have any thoughts but think only of CSWs. But if a man is in love and thinks only of his lover and marries her then there is no chance of the spread of HIV/AIDS." He knows the meaning of HIV/AIDS but has forgotten the full form.

He has an unforgettable memory of a person from his own village who was strong and sturdy but who became very weak after he was affected by AIDS, he married for the second time. "Perhaps his first wife also should have died of HIV/AIDS. Now the second wife is also infected. Both of them suffer from HIV/AIDS." This has made a deep impression in his mind. He says, view everyone in the same way whether it is a HIV/AIDS patient or otherwise. He does not know to discriminate between people. 'If someone says HIV/AIDS is not his problem one should immediately learn that he is frightened of it because he has led an improper life. The reason why men visit CSWs is their mental setup and the circumstances under which they live. Some of the prostitutes spoil young men and keep them as regular customers. That is the main reason why some young men become chronic HIV/AIDS patients. Being a man, I am not in a position to earn as much as these CSWs!"

In conclusion he gave us the following advice: "There are hundreds of TV channels. The government must pass a law that every channel must allocate a certain period of time to display HIV/AIDS advertisements. The government should also make sure that these ads are displayed in the primetime slot."

## 52

We met a 33-year-old farmer from Tiruvannamalai, now settled in Chennai. He is a science graduate. He feels that HIV/AIDS is a very indecent and fearful disease. Because it is an indecent disease people don't talk about it. Infections occur because of immorality. He says most people do not know the way HIV spreads but he knows how it spreads. According to him, 'it spreads only because man lacks self-control and seeks sex like a dog'. People go for very cheap sex. The rate of CSWs in some places is as low as Rs.20 to Rs.50. They buy the disease at a very cheap rate.

He feels the present generation is more inclined to sex than studies. When he looks at a HIV/AIDS patient he pities them and thinks that they have lost their life. People become HIV/AIDS either due to ignorance, or knowing fully well that they can contract it if they have unprotected sex. In colleges and schools they should make HIV/AIDS awareness a compulsory, separate subject and teach it. Only then students will know the basic values and would be careful in life. The message about HIV/AIDS can be conveyed to the public only through drama / cinema and the steps taken by the government and other agencies has not been sufficient. They only distribute condom but don't give importance to awareness programmes.

So far he has not taken up any bad habits. It is foolish to think that HIV/AIDS can be saved by faith in god because God is just one's mind and nothing more than that.

He has not seen any porn movies. The disease does not spread through air or casual contact. Only the police can control CSWs. But they don't take proper action since they are easily bribed. He feels little uneasy when he sees HIV/AIDS advertisement with family. However he has not discussed about HIV/AIDS with his friends. Some people shun talking about HIV/AIDS mainly because they fear it. Might be they feel guilty.

## 53

We interviewed a 28-year-old carpenter who has studied up to $10^{th}$ standard. His family is in Chennai for the past 40 years although he is a native of Kancheepuram. He said that one gets HIV/AIDS only when they visit 'call-girls' i.e. CSWs. One can talk and chat with call girls but never opt for sex because we will also get HIV/AIDS.

He started smoking and drinking when he was 15 years old. He learnt it from his friends. He suffered from a bad cough, so he started to drink in order to save himself from that tension.

"Everyone should come forward and take treatment for HIV/AIDS if they doubt that they suffer from it. Because they hide it, the children born to them are also affected by it."

He has seen porn movies. He first saw them out of a sexual urge. He says seeing porn has helped because he learned about true sex. It has really helped him to be careful. In general these films induce sex urge too. One should have a lot of self-control to watch porn movies; otherwise it will be easy to visit CSWs.

He sympathizes all the people living with HIV/AIDS, he won't mind helping them. If he gets a chance to meet them, he would make friends. He would learn from them about the disease and symptoms because only they can supply the first hand information.

He does not know whether the disease is contagious or not. 'If not contagious why do they pack the dead body of HIV/AIDS patient in a polythene cover?'—this was his question to us. He feels the government should increase the number of advertisements that deal with AIDS in TV and radio. This disease has come to India from foreign countries and now it is spreading very fast. The more advertisements he sees he will only resolve not to do mistakes. He has not talked about HIV/AIDS with his friends, but he has talked about it with his family members. He feels CSWs practice this trade because of their hunger and poverty.

## 54

We interviewed a barber aged 49 years. He has studied up to 8th standard. He says he is a native of Chengalpattu district but is currently settled in Chennai. He is married and has two children. He started smoking and drinking when he was 25.

He feels HIV/AIDS is not so alarming in Tamil Nadu. He says that before the advertisements and awareness programmes for HIV/AIDS came out, lots of people used to come to his shop for shaving and tonsuring. After these advertisements, the number has drastically decreased. He has suffered monetarily because of it. He feels that this is because of the cunningness of the companies that manufacture blades. They advertise that using the same blade, there is a chance for people to be infected with HIV/AIDS. They want to increase their sales, so they have done like this. He feels all these are only to alarm the people of the disease while in truth the disease spreads only by CSWs.

Rarely has anyone been hurt while he shaves, thus he wonders how HIV/AIDS spreads because of blades! We could not argue with him, we just recorded his statements.

He says that in his 22 years of service as a barber, he has so far not seen even a single HIV/AIDS patient. Yet, most of the customers who come for shaving demand him to use a new blade. They would openly say that they would get infected if he uses the old blade. He feels that if the government helps barbers like him in some other way it would be nice. He says that in the Tirupathi temple, thousands of devotees (both men and women) tonsure their heads and the barber uses the same knife. But no one says they would get infected with HIV/AIDS. He says that more HIV/AIDS patients visit the Tirupathi temple (praying for a cure) than the number of HIV/AIDS persons who visit an ordinary barber like him. He is not able to comprehend the nature of people—how in one place they feel they will get infected and in another place they don't even talk about it! He says that if the disease really spreads by blades or knives then why does it not spread when they tonsure their heads in the Tirupathi temple? He firmly feels that HIV/AIDS spreads only through unprotected sex.

He says that people visit CSWs even though they know they might get HIV/AIDS. That is why commercial sex work is flourishing. But the barbers' profession is suffering due to these advertisements. He said, 'without demand how can supply increase? The number of CSWs is increasing day by day solely because of the demand among the men.

*Additional Observation:*

*He holds a unique view about his profession. We could not counter him. Only for the sake of unbiased representation we have included this interview.*

## 55

We interviewed an uneducated man who was in his early thirties, working as a street-side hawker. He complimented us for taking interviews on such a neglected topic. He feels that people unreasonably fear HIV/AIDS. They have to use their reasoning ability and try to understand the disease. To wipe out the fear of HIV/AIDS one must use his brains. Without reason, the fear in people cannot be changed. He earns Rs.100 a day. He gives half of his earning to his family and the other half is spent on alcohol. He does not drink with his friends because some of them have cheated him. So he drinks all alone. He started to drink at the age of 15. The practice of drinking has made him a thinker, so he never commits errors. He admits to seeing porn movies regularly with friends. He sees one such movie every week. He does not see porn magazines because he cannot read.

**56**

We met a 29-year-old man working in the Navy. He has studied up to 10$^{th}$ standard. His native place is Virudhunagar, but he has been staying in Chennai for the past 8 years. He says one gets HIV/AIDS due to CSWs or due to homosexuality. He says nowadays injection needles are used and thrown away, so one cannot get HIV/AIDS by using needles.

He drinks and sometimes beats up his wife and children. He took to the habit of drinking beer and smoking at the age of fourteen. Only his friends made him take up these habits; now these habits don't go out from him, even if he wishes to leave them. To avoid HIV/AIDS one should never visit CSWs. He has so far never visited any CSW. Even before marriage he did not have such a practice. He lives only with his wife and children.

If each and every one observes his duty, AIDS would never spread. If a married person restricts his life to his wife and children certainly spread of HIV/AIDS can be restricted.

He cannot say whether the sex education is good or bad only if he attends two or more classes he can comment on it. He says he has taken free HIV/AIDS test when he had been to the village. He feels that if an employee is affected by HIV/AIDS he can inform the boss because others can be alert about it.

**57**

We interviewed a 22-year-old graduate from Chennai. She is unmarried. Her father works in a multinational corporate. She has not discussed about HIV/AIDS with her family although she has discussed it with her friends. She feels that finding medicine to cure HIV/AIDS will not spoil morality. Children born to HIV/AIDS patients are affected with AIDS, only the discovery of medicine can cure them. She feels that free HIV/AIDS test for one and all by the government is good, viable idea but she does not know how far people will accept it. She feels it is the social responsibility of everyone to help HIV/AIDS patients. She would come forward to help them; she would also come forward to spread awareness about the disease if given a chance. She supports the idea of giving sex education to school children. She says they can start sex-education courses from 8th standard. She says there is no relation between porn magazines and HIV/AIDS. She thinks by looking at these porn magazines people only satisfy their sex urges. She has seen HIV/AIDS patients and their plight made her sad. The number of HIV/AIDS patients is increasing day by day; the main reason is that most men fail to live with one woman. She feels that any person who gets HIV/AIDS is only responsible to a scale of 60%, the remaining 40% is because they are mislead by their friends.

She feels HIV/AIDS patients cannot be given preference in treatment because all patients are equal, irrespective of the disease from which they suffer. Free consultation about HIV/AIDS should be advertised with phone number in bus terminus, railway stations and in all AIDS advertisements.

## 58

We interviewed a 26-year-old post-graduate student. She was working in a private firm. Her sole ambition was to become an IAS officer. Both her parents are teachers. She feels sex is a main reason for the spread of HIV/AIDS in Bombay, Andhra Pradesh and Tamil Nadu. She has discussed about HIV/AIDS with her family members and friends. She says the public will never accept free and compulsory HIV/AIDS test. Perhaps if the results are made absolutely confidential people may take it secretly. She says it is not anyone’s duty to take care of HIV/AIDS affected persons. She thinks that people who have complete self-control can see porn magazines and cinemas. Nothing can shake them. She has seen patients infected with HIV/AIDS and was very upset by it. She hates the fact that public opinion is not in favor of HIV/AIDS patients. Even those who speak about AIDS are considered by the general public to be people with some defect. No one realizes there can be service-minded people who want to help these patients. However government / private sector can give employment to these patients. She feels preference must be given to HIV/AIDS patients, since this disease has no cure.

## 59

We interviewed a 29-year-old housewife from Chennai. She has studied up to 6th standard. She has two children and her husband is working in a private company. She says this disease cannot be cured that is why there is a great deal of social stigma associated with it. She is well informed about HIV/AIDS mainly through TV advertisements. The spread of disease can be controlled if men don't go to CSWs or even if they visit CSWs they use condoms. She says her husband has never suspected her and she also does not suspect him. She belongs to a middle class family and lives in a rented house.

She feels problems of CSWs are entirely different from other people's problems. They take up this profession due to poverty, frustration, and love failure. In some cases they drink, they act as extras in movies and are entirely misused by their boss and co-workers. To maintain themselves in the same social status they take up prostitution; for according to them it is a very well paid job for they also use the glamour they have earned through acting in movies. She sadly remarked, "This is an easy means of becoming rich. All that one needs is white skin and glamour." She feels that prostitution is modernized after the advent of cell phones and Internet. She has heard that over the Internet customers can see the photos of these CSWs. They are taken in AC cars and these CSWs without fear travel in night do the job and are returned to their place by the customer. Her brother is in the police department in Bombay and has said about the red light area there and its method of functioning. Small girls in the age group 12 to 16 are used mercilessly. These young girls are given hormones and steroid injections to look fat and sexually attractive. The side effect of all these drugs on these small girls is that their lifespan is affected and they are easy victims of HIV/AIDS.

Before they reach the age of 30 they are affected by all kinds of STDs. All of them cry in their heart and hate doing this work. Several of them end up as cervical cancer patients in the latter days of their life. After 16 years of entertaining men for sex makes them totally invalid in their late twenties or early thirties.

They are usually beaten up and forced to entertain men when they are just in their early teens. She says her neighbour who is just 32 years old was rescued from Bombay red light area and is now her close friend. She has told several horror stories of her life.

This woman told her friend's lifestory, "When she was just 11 years old an elderly lady took her from her parents and promised to find her employment as a servant maid. But to her surprise she was not taken to any house but to Bombay. As she had no knowledge about the place and language she could do nothing. For several days she was closed in a room without food and beaten up. She was given some injections. She does not know why these injections were given. Now that woman cannot move due to cancer. She may not live long. She goes for regular treatment. She is given some protection by her maternal uncle. Only if women like her are rehabilitated and commercial sex work is banned, HIV/AIDS can be prevented."

## 60

We interviewed a 32 year-old man working as an accountant and assistant in a pawnbroker's shop. He feels the basic reason for the spread of HIV/AIDS is vulgar scenes in cinema. "Though in Tamil Nadu they speak of woman as a god, they are only portrayed as an object of pleasure in movies. Every man who sees such movies will only feel like having sex with women, nothing more. How will we not get HIV/AIDS?

Unless such movies are banned it is impossible to make people behave with self-control."

He added, "Most of the Hindi films have a dance scene shot in a bar or a red light area or a song sequence portraying a CSW. Our men stand in long queues in the sun to get ticket and see this vulgarity. What happens after seeing? They visit CSWs."

He feels spiritual involvement will certainly channelize them into good path. They will not become a victim of bad deeds. "Young women should dress properly, their dress should not provoke the emotions of men. We should behave in accordance with our culture. Because of western culture our women have started to dress very scantily which is the reason for several problems like eve teasing, rape. Nowadays girls in the age group of 10 to 30 dress up like film actresses. They forget that these film actresses show their body only to earn money. What they see on TV serials or movies is only a story and not a reality. Movies must inculcate social reforms and ideals to educate and eradicate the base thoughts in youth.

Present day movies look more like porn movies, and how does the censor board pass all that? Song-sequences and dream-sequences in most of the movies are sufficient to make even a strong-minded man, a sexually obsessed man. They directly contribute to the spread of multi-partner sex. That leads to AIDS."

He further remarked, "All directors and produces see to it that their movies have enough obscene poses and sexual activities so that the film can run for at least 100 days. What national responsibility do these directors, producers, actors and actresses have? They should be put behind the bars for such acts."

He is unmarried. He has no bad habits, he has not seen porn movies or magazines. He has not been to CSWs, so he has no idea about sexual relations with CSWs. He feels that if someone refuses to share his / her ideas about HIV/AIDS, it is their right so no one can comment upon it.

He believes in god, but adds that god will not save AIDS patients. He does not know anything about STD or condoms or CSWs. He said in conclusion, "No religious leader or politician talks about HIV/AIDS. This is because they themselves lead a loose, immoral life."

## 61

We met a 26-year-old art graduate from Chennai. He says the only reason for the spread of HIV/AIDS is the absence of self-control in people. He has no bad habits. When he was 15 years old he saw porn movies. "They are not bad at all. They educate us and relieve us of our inhibitions."

When asked about AIDS, he said "first one will get STD, afterwards HIV and in the final stage AIDS. It takes 5 to 10 years for HIV to become AIDS."

He feels that if the government bans CSWs then all people have to suffer. He says fear about HIV/AIDS exists because of secrecy about sex. Openness about these things will remove all fear of HIV/AIDS.

"Only through advertisements the spread of HIV/AIDS can be controlled. The government has succeeded to a great extent. When I was a school boy I didn't know about AIDS. But nowadays, first standard children know about AIDS."

He said with a note of tension in his voice, "even if someone really became HIV/AIDS infected by blood transmission, society will only think that the person got the disease through CSWs. This is the view of the world."

## 62

We interviewed a 22-year-old B.Com graduate working in a liquor shop in Chennai. His father owns this shop for the past 30 years. He is rich. He declined to answer questions about CSWs. He never goes to theatre for movies. By faith of god no disease can ever be cured. It is a matter of destiny. He says the main reason for the spread of HIV/AIDS is migration from villages to cities. He says even if he wishes to visit CSWs in Chennai certainly one day or other everyone will come to know about it. So he would not go but this is not the case with migrant workers. Thus only migrant labourers freely visit CSWs for they know that nobody from home will ever know it!

## 63

We met a 23-year-old man. He was an 8th standard dropout. He smokes and drinks. He watched a porn movie at the age of 12. For a long time he was addicted to it. Now he does not see them. He says, "even if we see porn movies we should think of the good we can take from it. If everyone views movies this way and try to learn a lesson then certainly the world will improve. For five minutes pleasure, a man becomes a HIV/AIDS victim by visiting CSWs. If such things continue one day the world will be destroyed by HIV/AIDS."

He says he has not seen any HIV/AIDS patient.

He emphatically said, "Only religion and religious practices are the root cause for the spread of HIV/AIDS and misconduct among married couples. In some families women abstain from sex saying that they are doing some prayer and sacrifice. At that time the men think their ego is hurt and go for other women. Likewise few men also abstain from sex with their wife saying that they are under some vrtha (vow), so these women visit their neighbour for sex. In due course of time they develop new contacts, thereby HIV/AIDS spreads."

## 64

We met a 47-year-old man who owns a petty shop. He sells condoms. He has seen porn movies and porn magazines. He sells them even today. Youth in the age group 20-30 buy these books; CSWs also buy these porn books. Porn books are sold out within two days of their arrival. He has special and regular customers who take these books from him. He never sells these porn books to unknown youths. They are the most fast-moving commodities of the shop. Further all usual magazines, which have indecent covers, are sold in a day or two.

He has seen both porn books and films. He saw such movies mainly because of his friends. "If a man has self-control he can see any number of porn films and it is not wrong. Only a man without self-control must not see it, because he would immediately go to CSWs." He also sells condoms, so he appreciates condom promotion campaigns. He says condom sales in his shop are mediocre. People in the age group 20 to 60 buy condoms, but certainly people in the age group 40-60 buy more. Even CSWs buy condoms from his shop. He finally said, “We should demand chastity for men in the same way as it is demanded for women.”

## 65

We interviewed a 25-year-old policeman from Sathoor who is currently posted in Chennai. He has studied up to +2. He says he has no bad habits. He is unmarried and celibate. He learnt about HIV/AIDS only through advertisements, he has never heard about STD or VD. He has not discussed about HIV/AIDS with his friends. He says that HIV/AIDS can be completely cured through faith and prayers to God.

He said that CSWs must correct themselves; they must feel disgraced to carry out such a low trade. For the sake of money, one cannot sell their body to so many men. He asked, “When a man cannot go home and tell his wife or parents that he visited a CSW, why should he do that mistake in the first place?”

## 66

We interviewed a final-year B.Sc. student from a government college in Chennai. She has not talked about HIV/AIDS with any one. She feels that one must certainly help HIV/AIDS patients because they suffer mental tension because of the disease. She has so far not seen any HIV/AIDS patient.

She says god cannot cure the disease. “The creation of god itself is only for our mental peace and we have created him. How can our own creation cure us? God is only for man’s consolation.” She refused to talk more about AIDS.

## 67

We met a 22-year-old lady who was handicapped. She has studied up to 8$^{th}$ standard. She is self-employed and sustains herself by making wire bags. She comes from a poor family.

Her father was a drunkard, always coughing because of TB. She used to fear that her father might have HIV/AIDS. Her mother sells flowers for their living. Her father is not earning, he is sick with TB. She is well aware of HIV/AIDS.

Her mother's youngest brother died of HIV/AIDS in Kancheepuram. Since she was handicapped he did not wish to marry her. When he died he was just 27 years old. He was taking treatment for 1½ years in the Tambaram hospital. Some brown colour liquid was oozing out of his body which was wrapped in plastic sheets after his death. They were advised not to open the cover and warned that it would infect them. As the body could not be taken to Kancheepuram they took it straight to the burial ground. Her father refused to bring his body to their Saidapet home for the fear of social stigma. She said that her deceased uncle looked dull blue in colour. They said it was due to medicines. Further he weighed less than 20 kg at death. He was so thin and the body had no flesh.

She says that he died because of his arrogance. He married another girl just 2 years before his death. Now the wife fears that she too may be having HIV/AIDS.

She said "god is great that he created me as an handicapped woman. Otherwise my mother would have surely married me off to him and now I would be a HIV/AIDS patient."

## 68

We interviewed a 24-year-old woman who has studied up to MCA. She stays in T.Nagar and is looking out for a job. She knows about HIV/AIDS. She says that in Chennai, people get HIV/AIDS only through sex. Only in other countries, people get HIV/AIDS through injections or blood transmissions. Even god cannot cure the disease but these men will start to pray only after getting HIV/AIDS. Men in rural areas don't waste time or money on god. Only when they feel incapable, they surrender to faith for mental peace. Scientists should find methods to cure HIV/AIDS. It is not only prevalent among men, it is also found in infants and women, and everybody must be cured.

The government should soon find medicine to cure HIV/AIDS. She says people who are infected with HIV/AIDS

must first be segregated from their family. They should be treated like prisoners, only then the government can stop the spread of the disease.

She says the plain psychological reason is that when men are infected they do not want to use a condom. They don't mind infecting their wife or a CSW.

### 69

We met a 14-year-old boy working as a mechanic. He has heard about HIV/AIDS. He failed in the 6$^{th}$ standard so his parents put him here so that he could learn this trade. He earns Rs.15 a day. We were surprised to hear that he knows about HIV/AIDS. We asked him how he came to know about it. He explained that one of his seniors (aged about 28 years) died of HIV/AIDS last year. He lost weight and had continuous fever. So their employer sent him to a doctor for check-up. Finally, they found that he was suffering from HIV/AIDS. It did not come as a big shock to them because they knew that he often used to visit CSWs.

Whenever their employer went out to make some purchase, he used to slip away telling this boy (interviewee) to look after the shop. Like smokers go out of their office and smoke once every hour, the deceased used to leave the shop and visit CSWs thrice a week. That is why he got the disease. It taught him a lesson that he should be very careful in these 'personal' matters. He said that his friend's father had died eight years ago. His mother worked as a domestic help to make both ends meet.

Even their employer pitied that woman and helped her with money. He always advised that boy to stop his philandering lifestyle.

This 14-year-old boy had often accompanied his friend's mother on her visits to the hospital. Seeing his friend, he learnt that one should never commit any mistakes in life.

The deceased weighed less than him. Day by day, his friend became smaller in size. Only his height didn't change. No one could save him at all.

One day, the doctors asked the boy's mother to take him home because he was so close to dying. He had lost consciousness, and rarely opened his eyes. His body carried an oppressive stench. He could not eat anything; whatever little milk or water his mother poured into his mouth came out. No one ever wished to come near him even when he was little alive. Only his mother was all the time with him. Our owner took all the expenses of cremation. His friend's dead body resembled a skeleton wearing clothes.

Because of this experience, the boy says that he learnt both the horrors and the lessons of life. He has even stopped seeing porn movies.

**70**

We interviewed a 2nd year medical college student. To begin the conversation, we asked her if she could discover some medicine for HIV/AIDS. She said, "Only if there is an incurable disease like HIV/AIDS people will lead a proper life. Otherwise both men and women will go for multiple sex partners." It was surprising for us to hear such things from a girl studying to be a doctor.

She is of the opinion that moral values have deteriorated so badly that men at 70 seek sex from small girls. Some women in their mid-forties have sex with teenage boys. She says the world is going in a topsy-turvy way. Only HIV/AIDS can stop flimsy affairs from taking place. She was very pessimistic and discontinued the interview.

**71**

We interviewed a woman police constable; she was in her early thirties. She is well aware of HIV/AIDS. She says rural men are more prone to HIV/AIDS than city men. We asked why she said so. She said, "Because of AIDS awareness in city, these men use condoms when they visit CSWs. On the other hand, rural men are not so aware, so they make mistakes and end up with HIV." As a policewoman, she has met many CSWs. But she is vexed with these men, because any amount of counseling will not change them.

## 72

We interviewed a physical education teacher working in a private school in Adyar. He is from a village near Cuddalore. We asked him to narrate any unforgettable event in his life related with HIV/AIDS. He told a long story, which is as follows: "One of my rich classmates failed in 12th standard. He left for Bombay. He had no bad habits or association. He used to be a little aloof because he was rich. This boy did not meet me for ten years after school. When he came back from Bombay, he was very ill. His father begged me to take him to a good hospital and give him good treatment so that he can be cured. Accordingly I went with the patient and his father to several hospitals.

To my surprise the doctors disclosed to us that he was fully affected by HIV/AIDS. Several major, reputed hospitals declined to admit him. Several mistook me to be the patient's younger brother. Everyone recommended us to the Tambaram Sanatorium. We met a doctor named Dr. Deivanayagam. He was entirely different from all the doctors whom we met. He asked us to join my friend as an inpatient. The doctor was very affectionate and kind enough to touch the patient.

The doctor asked my friend how he got AIDS. My friend said he did not know anything about how he had got this disease. It was very clear to us that he had visited CSWs from the red light area. Because he had lots of money, he had sex with varieties of CSWs from different states. The father was in tears. In spite of the treatment he died within a year. The government should come forward to honour Dr.Deivanayagam's services as a specialist doctor for HIV/AIDS. Doctors who do service to patients for heart disease, eye disease or cancer are given awards. It is high time government awards doctors for doing selfless service for AIDS."

## 73

We discussed about HIV/AIDS with a pastor from a Chennai city church. He never wished to talk anything other than the religious angle. He said, “For the past ten years several men with HIV/AIDS have visited the church. Some of them have talked heart to heart about their activities to me. Several of the patients came to me for comfort and consolation. Some of them converted to Christianity. Even two Muslim women with HIV/AIDS became Christians.” The pastor said that people were relieved of their tension after they confessed their sins. Patients would ask him to give them a small cross, or photo of the Virgin Mary or Jesus or Bible to keep it with themselves. “When terminally ill cancer patients are given prayers why do they deny such things to HIV/AIDS patients?”

The pastor is of the opinion that they can give up smoking, or alcohol; but CSWs they visit even in their bad state of health, they are not able to give up, that is what some of them confess in their meetings because they feel that the pastor will always maintain their secrets. This is the major reason why they die soon after being infected.

**74**

We interviewed a 28-year-old research scholar who is working for Ph.D. in biomedical sciences. She is a postgraduate in chemistry working on protein levels in cancer patients. We asked her if she would do such research for HIV/AIDS patients. She said she has chosen this problem some seven years ago.

If we had met her some seven years ago her problem for research would have been AIDS.

She said that AIDS is not accepted by society because immorality is not accepted by society. She was also saddened by the fact enough research is not done for HIV/AIDS women patients. She said, “It was reported recently in the dailies that nearly 2 lakh women are infected by their husbands. Unless the traditional set up is changed it would be impossible to stop the number of cases of women being infected by their husbands.”

**75**

We interviewed a 72-year-old uneducated man. He works as a toddy-tapper. In frustration, he said that the current generation lacked character and morality. He says today's youngsters have no character or morality. Their love is only for sex. "Sex is not love and love is not sex. They have not understood this."

He was married when he was eighteen, after 10 years his wife died. They had no children. So his relatives once again married him off to another woman. She lived with him for four years. He came to know that she had sex with another man in his absence. When he questioned her about it she ran away never to return.

Thus he was married for the third time. Now he lives with his third wife and they have two children. He is very proud that his daughter has been married off.

He had never been to a CSW or had sex with women other than his wives. "That is why even at the age of 72, I am able to climb such tall trees and earn my living. I do not depend on the income of my son."

"The only way to eradicate the spread of HIV/AIDS is that every man should be loyal to his wife. If this type of morality is practiced not only can we chase away HIV/AIDS but we can have a healthy and a youthful life."

"All films show half-naked women chased by a lusty man! Is this love?" He was full of praise for MGR and the actresses of those days. Do we have any actresses or actor like those of olden days? Now if they are naked the film becomes a hit and runs for over 100 days."

He was very talkative. We could not counter him because we were afraid that our interruptions might make him hide his true experiences. We asked him how he came to know about HIV/AIDS? "In my village a rich man died of HIV/AIDS. They were 'upper' caste people. They have nearly 200 acres of land and I used to work for them. For enjoyment, the son of the landlord went with his friends for a 30 days tour to Hyderabad, Bangalore and Bombay.

He was just 23 years when he died. His brother works in Dubai. After his return from the tour he was sick. They took him to Chennai for treatment. He had boils all over the body. Pus came out of such boils and no one could ever sit by his side. The villagers said that he died of HIV/AIDS."

## 76

We met a 51-year-old chef who is originally from Tanjore but is settled in Chennai at present. He is a arts graduate. When he became a chef, his father helped in washing vessels and cutting vegetables. After his father was unable to assist him, he employed two women to assist him in the job. Both were from villages near Tanjore.

He was not much aware of the background of these women. Last month, one of them died due to HIV/AIDS. He is yet to recover from that shock. He did not have any non-professional relationship with these women. For three years she had been sick. She tested herself for TB, and he paid her Rs.5000 for the medical expenses.

She never recovered. One day he got a phone call that she was dead. As a courtesy, he went to attend her funeral. There he learnt from her relatives that she had affairs with several auto drivers and she had died of HIV/AIDS. He said this true incident only proved that nobody could guess who has AIDS and who does not have AIDS.

## 77

A 25-year-old diploma holder, who worked as an assistant in a cyber-café (browsing center) in Chennai, spoke to us about his opinion on HIV/AIDS. He is an orphan and he supports his only sister. His owner pays him Rs.3700/- per month, which he feels is a handsome salary. The cafe has 12 computers; the browsing cubicles are partitioned so that one browser cannot observe what his neighbour is viewing. Majority of the customers who visit the café are youth in the age group of 15 to 25 years.

They visit all kinds of porn sites on the Internet. In some cases, couples come to visit the café to see porn. He has himself glanced at these porn sites because some times the customers call to him for help. He feels these cyber cafés are sufficient to promote the spread of HIV/AIDS. These youngsters who frequent the café look like addicts. The order from his boss is that he should never argue or even look at the customers when they are with the computers. He sits in a chair and most of the time closes his eyes so that customers are free to do their work. He feels that his job is worse than the job of a pimp. The unique feature about these customers is that they are tongue-tied. Some of them pay in 500 rupee notes. Some stay in the café for over 4 to 5 hours late in the night!

He said that the spread of HIV/AIDS among youth could be attributed to these Internet cafes. He feels that if he is a policy maker he would ban all Internet cafes and run them like public telephone booths where they can only send messages. He would also order ban on all these porn sites. “One can use Internet to acquire scientific or technical knowledge. Can we use the sites to see porn movies or nude pictures of women? Is it nice?” He did not speak any further about HIV/AIDS.

## 78

We interviewed a housewife in early fifties living in a very posh apartment. She has two children; her husband is a bank officer. She is well aware of HIV/AIDS. Even in her own apartment, a lady in her late thirties died of HIV/AIDS one month ago. “Only in the cremation ground, her close relatives told the news that she had died of HIV/AIDS. Very few people attended her funeral. She was a spinster living with a servant maid. Both men and women visited her. The servants working in her apartment said she was affected with TB. Most of the inmates of the apartment are of the view that she died of TB.” She says that the husbands of the CSWs do not mind their wife earning money by ‘indecent’ means. In some cases they act as pimps. Thus the world has become very corrupt in this ‘Kali Yuga’.

She says CSWs are prevalent in her own apartment. One woman in their apartment is very rich, owns two cars, sends her children to school by car, her husband is unemployed. She is in her late thirties. She is a CSW who caters to very rich clients. We were very surprised at this information. She feels that unless men and women individually choose a moral life no great changes can occur in our society.

## 79

A nearly 27-year-old waiter who works in a five-star hotel was interviewed. He said he is a failed science graduate. Since he knows English and Hindi also, he gets more perks. He is aware of HIV/AIDS. He gave us a customer profile of his restaurant. In his opinion only 10% of customers come to his hotel with their whole family. Some come there to host dinners. Some come to drink and make merry. A few people bring their lovers and spend time together. Some men bring along a CSW and they behave so lewdly that these waiters have to put down their heads in shame.

This man said that he never made any mistake in his life because he was only a runaway from his village. He also had a neighbour in the village who died of HIV/AIDS, so he has learnt a lesson never to make any mistake. In order to deal with his emotional disturbances he has taken up smoking. He said, "I will never take up drinking because one day I might lose my mental balance and end up visiting CSWs. When they are available for Rs,30, which man will not make mistakes?"

## 80

We interviewed a very experienced social worker. He was in his early fifties. Soon after his graduation from the Madras Christian College he joined as a social worker for a paltry stipend of Rs.175/- p.m. in the Tambaram Sanatorium hospital, which was known as TB hospital in those days. His work was to go on rounds to see whether patients were given good food, medicine etc. and also to assist the doctor when the nurses were not available. At that time the hospital was purely run for TB patients.

After several years of work he joined as a social worker with an NGO. Now his work is multifaceted. Sometimes he works with rag pickers, at times with street children, sometimes he was assigned to assist mentally retarded children near Perunkalathur and so on. Since he is paid well, he takes up the assignments. Once, he pretended for a whole week that he was an AIDS patient. Under this guise he moved with the other patients. He observed that young men in the age group 20-30 were very frank about the disease. Whereas middle-aged men said that they were infected by blood transfusion or by needles or some silly story.He also observed that many of them continued to visit CSWs even when they were admitted for treatment. Even after it is confirmed that they have AIDS, they will not wear condoms when they have sex with CSWs. They say in filthy language, "After all I got infected by a CSW, so why should I bother about infecting CSWs?" The social worker is very depressed. He says there is no way in which people can be improved.

## 81

We interviewed a man who runs an NGO involved in condom promotion for the past 10 years. He said, "Men who visited CSWs could not be changed. But I appreciate the fact that I could change the CSWs and can force them and make them use condoms.

The men who visit them never listen to the CSWs advice of using condoms. We have our own limitations. If these men get angry with the CSWs who insists on condom usage, they boycott her. So, out of vexation these women entertain men even if they don't agree to using a condom." He is reminded of a sad story shared by a CSW.

"CSWs have a pimp who regulates the customers and it is an ardent duty for them to pay him at least 60% of the income. It varies from place to place and locality to locality. Unless he gets his due income he would be angry with them and torture them. At times he would physically abuse them. This woman was highly worried about the torture so she used to entertain men without recourse to safe sex. Consequently she was always subconsciously frightened that she would be HIV/AIDS affected. She had no money or time to go for such a test. She had no one to help her in this matter."

He said that the social stigma associated with AIDS was a necessary deterrent which could prevent the immoral attitude of our society.

## 82

We interviewed a doctor who had a private practice. One of his close friends was a pharmacist who used to refer him to patients with STD. Majority of the people who visited him for treatment for VD were poor and lower middle class. They were either daily wagers or contract labourers with no education. These patients usually visited him in the night after 8 p.m. when he was about to shut the clinic. So he could give them treatment without any hurry. He has treated over 90 patients in a year. They were in the age group of 18 to 65. Most of the youngsters who visited him were unmarried. When asked how they got the disease, they used to tell convenient lies: allergy, heat and so on. So the doctor used to tell, “If it is due to heat or allergy, it will be naturally cured. There is no medicine with me for it.” So they would start to disclose details of their visits to CSWs

Only 9 of his patients were found to have HIV/AIDS. He referred them to the TB Sanatorium.

He said, “Men in the age group 35 to 45 never say the truth. They only say that they got the infection from their wives. When I asked them to bring their wives they used to blink and give a sheepish smile. Then they would admit that under the influence of alcohol they had visited CSWs.”

This doctor feels that patients who have no shyness in visiting CSWs have shyness in confiding their activities to a doctor. Sixty percentage of the patients’ sufferings are due to stress and mental worries. Rural patients suffer because they are ostracized by the society.

Many of his doctor friends do not treat HIV/AIDS patients. They don’t have anything against the disease. But they are afraid that they will lose their other clients if it is publicly known that they treat HIV/AIDS.

## 83

We interviewed a young man who has been working as a ward boy in the Tambaram Sanatorium for the past 8 years.

He says that the number of patients visiting this hospital has increased tremendously. “The hospital was not so crowded in the past. The number of beds has also been increased considerably. There has been a considerable increase in the number of patients from Andhra Pradesh." He said, "Male patients do not reciprocate our affection. But the women are much different." We asked him whether he fears that he will get the disease? He was so philosophic; he said, "Everyone must die one day or another, so why should I fear?" Then he began to say his story. "I do not suffer from HIV/AIDS. I have a father and a stepmother. My mother died a long time back. I had an elder brother who was infected with HIV/AIDS. My father took him to the CMC Hospital in Vellore. There we were advised to take treatment in Tambaram. My brother was admitted here as an inpatient. For one whole year I was by his side. When he died I did not go back home because I was scared of my stepmother. So I wept day and night and stayed in the hospital itself. A nurse took pity on me and employed me as a helper. From then onwards I continued to work. Now they have given me the post of ward boy and I earn well enough to maintain myself. So far my father has not come to see me. I too have prestige, so I do not want to visit my father."

He said, "All male HIV/AIDS patients couldn't stay with out 'women' (CSWs). Even in this hospital at least 85% of them visit CSWs in the nights by cheating the hospital authorities. So how can they ever be cured? The treatment is only a waste!" This boy said that these HIV/AIDS affected men must be given proper "coaching' class" to behave well. Might be, by the term "coaching class" he meant, "counseling".

He says that people do not divulge their HIV seropositive status due to the following reasons: fear of losing their job, fear of losing good friends, fear that relatives will isolate them, fear that it will not be possible to marry their children, fear of getting a bad name in society, fear of discrimination, fear of being branded 'immoral.'

## 84

We interviewed a 36-year-old man who left home at the age of 12 and for the past 24 years has been working as a rag picker. He says he earns a decent amount and he spends to his heart's content. He is happy that he is in a position to help his poor and sick friends. Due to the cruelties of his stepmother, he ran away from his home in a village near Dharmapuri. He never felt like returning home. For some time, some missionaries provided him care and shelter. Under their influence he was able to gain moral values and discipline. But he later ran away from them.

He feels that HIV/AIDS patients get the disease because of their bad actions. According to him, HIV/AIDS patients have sinned and hence they should suffer.

We asked him if he knew any rag pickers who had HIV/AIDS. He said that his friend had suffered from constant fever, so they used to buy him food and paracetamol tablets. He never recovered, so they approached a social worker and he helped them to admit the friend in the Tambaram Sanatorium. The friend died a year later. Only then these people came to know that he had died of HIV/AIDS. Since he had no relatives, only his friends had to perform the funeral rites along with the help of a NGO.

He cried as he recollected his friend. His friends (and himself) do not differentiate any disease. He says that this is because they are doing the 'filthy' job of rag picking for a livelihood, so they do not consider any disease to be good or bad.

He said that even the rag pickers were treated like HIV/AIDS patients by the public. People move away from rag-pickers and hold their noses with hand kerchiefs. Working in this profession, they are homeless, they lack public bathrooms, and they don't wear perfumes, and can't afford clean clothes.

He was very philosophic.

When asked if HIV/AIDS was a social stigma he said it depended on a person's perception. Everyone has the right to think whatever he likes. But his right should not hurt the feelings of another!

## 85

We interviewed a porter at the Central Railway Station, Chennai. He claims to be 65 years old, but he has cheated the authorities at the time to registration by saying that he is only 55 years old. He is proud of his young looks. He has worked as a porter for past 46 years. He has no bad habits. If he is very tired he consumes liquor in limited quantity. When we asked him about his opinion about the spread and prevention of HIV/AIDS, he shared his life experiences. He is well aware of HIV/AIDS. Some of his friends have died of it. He talked very degradingly about the easy availability of CSWs. He said they were available for "$3/4^{th}$ of an Anna' "Anna is old system of currency where 16 Annas make one Rupee". He condemned porters who visited CSWs. He thinks that unless there is a total ban on CSWs it is impossible to save India from HIV/AIDS. Every year the number of HIV/AIDS patients only increases. When he recently visited the Tambaram Sanatorium, to see his colleague and friend, he saw that the number of patients has increased.

Thus he says his view is hundred percent correct that HIV/AIDS patients are increasing very steadily. Further, he is of the opinion that free condom distribution program is of no use. It is like advertising, "Smoking is injurious to health" and at the same time selling cigarettes or it is like saying "Drinking will ruin life" and selling liquor. Now the government itself sells liquor under the TASMAC brand. Likewise free distribution of condoms indirectly supports visits to CSWs. Thus he feels that the government does not bother about CSWs. If this is the case, how will the number of HIV/AIDS patients decrease? He says that the poor are the worst affected. Their food habits due to poverty and access to CSWs leave them only as HIV/AIDS affected patients. Majority of the uneducated poor, who come to city in search of jobs are sure victims of HIV/AIDS. He feels only one in 100 migrant uneducated poor labourers can escape without HIV/AIDS. Independent of their marital status, they visit the CSWs when they come to Chennai for jobs.

Their three leisure activities are watching cinemas, visiting CSWs and drinking liquor. Who can save these men? Whatever advice he gives, they retort, "This old man has no work." Only when they are badly affected they come to him for monetary help.

Then he shared his bitter experience by saying that these HIV/AIDS patients are not treated well even after they die. He says only the dead bodies of HIV/AIDS affected persons ooze out some fluid that stinks, so the corpses are put in polythene covers and sealed so that no contagious spread of diseases takes place. In his opinion, free condom distribution may help the urban and semi-urban populace but not rural men, who generally shy away from condom usage. Thus he feels only god can save these rural uneducated men. He is not in a position to find means to help them. A man should lead a healthy life of self-realization. Unless change in attitude comes from within, no external force can change people. Among the porters working with him, 50% of them are originally from rural areas. They have come to the city for one reason or other. He blames the failure of agriculture and the absence of means of sustenance in rural areas are two reasons why they migrate from the village. The most important thing is that these migrant men do not have land of their own in the village. They are landless coolies. Had they owned land, they would not have migrated to the city. The porter is certain that majority of the people affected by HIV/AIDS are very poor. He wants the government to help them. If a person comes to the hospital for HIV/AIDS test, and if he tests positive, the first primary step the hospital should do is to bring all his close relatives, counsel them and test the wife and small children for HIV/AIDS. Unless such stern steps are taken it is impossible to stop the increasing cases of HIV/AIDS. At some point of time India may not find labourers for employment if this situation continues! He concluded by saying that people in white-collared jobs can continue working even if they are infected by HIV/AIDS, but people like him whose profession depends on their physical ability will be worst affected by the disease.

## 86

We discussed about HIV/AIDS with a male nurse from a very reputed hospital in Chennai city. He must be in his late twenties. We asked his opinion about HIV/AIDS patients. He said, “Very rich men suffering from HIV/AIDS come to this hospital. The special thing about our hospital is that we promise 100% confidentiality to the patient. Our hospital has treated even politicians and top businessmen in the city. Some of our doctors even go home and treat these people. Everything can be achieved through money in this world. Some very rich men even go abroad for their treatment. If they come here for treatment in the earlier stages of the disease, they can prolong their life. However, if they come in the later stages, when the disease is highly chronic, there is nothing that can be done. The rich are very careful in maintaining their anonymity. They can protect themselves from suffering. But the case of the poor people is not so. Not only do they lose their health, they also lose their status. It is the responsibility of the poor people to understand that they have great responsibilities towards their wife and children. They should never touch liquor or visit CSWs. They should work to come out of their poverty. Free condom distribution is not going to help the poor rural people. It is the foremost responsibility of politicians to talk openly about HIV/AIDS. They must make it their duty to inform the poor villagers not to lose themselves on such slippery paths." After this speech, he was ready to answer other questions about HIV/AIDS. He feels that the government can order private hospitals to give subsidized treatments and free tests for HIV/AIDS. This will help people a great deal.

## 87

We interviewed an electrician from Saidapet. At the time of the interview, he was accompanying his wife's brother (an HIV/AIDS patient) to the Tambaram Sanatorium. His brother-in-law was in his early thirties. He has been taking treatment for the past four years.

The electrician says, "When they came to know that my wife's brother had AIDS, all his family members asked him to run away from home or to commit suicide. I took pity on him and asked him to come with me. I saved his life. My wife and I admitted him in the hospital and after a week's stay he was discharged. He has been asked to report every month for a check-up.

Now, not only he is healthy, he also found himself a government job as a lower division clerk. He lives with us, and we don't discriminate him.

Not all diseases are curable. There are a few incurable diseases like diabetes and cancer. So AIDS is not a big deal!" He went on with his narration: First I wanted to train him to be an electrician. However he thought that this job was very tough. So I asked him to join a typing, short-hand and computer course and now he is employed."

He is really very proud of his good deeds. If a HIV/AIDS patient leads a healthy life by not getting addicted to CSWs or other bad habits certainly he can lead a normal life. He should eat good food rich in proteins. Then he spoke of his brother-in-law. "He wants to lead a life as a bachelor because he does not want to make one more person a victim of HIV/AIDS. His parents now visit him but he wants to live only with his sister and me. He is very wounded at heart and cannot lead a normal life with his parents. Further they may differentiate or ill-treat him. He said to his parents, "You are not like my brother-in-law, who spontaneously forgave me and made my life meaningful. Except for him I would not be living now."

*Additional Observation:*

*If everyone who was affected by HIV/AIDS had a relative like this electrician certainly it would be a very easy task to rehabilitate the HIV/AIDS patients and wipe off the stigma.*

**88**

A sixty-five year old cobbler who has a pavement shop in Adyar was interviewed about his opinion on HIV/AIDS. He said that due to various problems he has been working in Chennai for the past six years. The ban on cow slaughter left him on the streets. So, he had to come to the city in search of a livelihood. His family had to remain in the platform itself for three months. Now he has built himself a tiny hut in the poromboke land, he has no idea about when he will be chased out of it. He kept lamenting that it is only the rich people who suffer from AIDS. "If the poor make a mistake everybody comes to know about it. On the other hand, if a rich man does any mistake it is immediately covered up. Let there be a compulsory test for everybody. The real truth will be revealed. The world will come to know that only the rich are immoral and not the poor."

## 89

We interviewed a postmaster who is in his early fifties. We first asked him the reason for spread of HIV/AIDS? He said it was because of modernization. For a minute we were little surprised. We asked him why he felt so.

He said, "Modernization makes youth the victims of HIV/AIDS. Nowadays most of the actresses dress in miniskirts, swimsuits etc. All this is very new to the Tamil land. This perversion in dress certainly has a very deep impact on teenagers.

"This sort of modernization or aping the west in dress code has certainly polluted the young minds and has spoilt them. That is why one can see lovers everywhere in the city. Thus media and its modernity has increased the possibility of youth becoming HIV/AIDS victims."

"Secondly, modernization affects poor people also. They do not get enough jobs in the villages, so they come to the city in search of livelihood. They mostly end up working as construction labourers, lorry drivers, bore-pump workers. Their migration makes them susceptible to HIV/AIDS." This postmaster blamed the loss of tradition as the reason for spread of HIV/AIDS. He wouldn't speak about the other manifestations of the disease.

## 90

We interviewed 12$^{th}$ standard boys of a central government run school in the city. A vast majority of them had experimented with smoking and drinking. However, they all maintained that they had not had sex. They were aware of AIDS. Though they knew how the disease spreads, they had many questions that they wanted us to clarify. They said that their parents had never discussed HIV/AIDS with them. These boys also dared not talk such a topic with their parents—because they would be immediately suspected.

Their curiosity to know about the disease and its symptoms was very high. They were crazy to know several things. They said there was no stigma associated with AIDS. It was just another disease, but not so fashionable, that's all. They were giving many suggestions for AIDS treatment such as: replace all the patient's blood with pure blood. They are willing to go and spread awareness to people about it. They say that people collect funds for orphanages from schools but nobody collects donations for AIDS from schools. So in primary classes students don't know anything about AIDS.

These students had utmost sympathy for the patients and wanted to help them. They said they would be happy to share their classroom with a student having HIV/AIDS. They feel that segregating these infected students is a denial of human rights and such evil should never take place. They felt that discrimination is more cruel than the incurable virus itself.

## 91

Consequent to the earlier interview, we interviewed girl students from the same school. When we approached them and wanted to discuss HIV/AIDS they were very reserved and silent. Most of them had learnt about HIV/AIDS in the sex-education classes. They admitted that they came to know more about AIDS through their friends than through the classroom lessons.

They asked us if AIDS was really as dreadful as it was portrayed to be in the media. Unlike other people, these children did not view AIDS with stigma or contempt. The fact that it also spreads through sex was immaterial to these girls. We were happy that at least the younger generation does not discriminate baselessly. They felt that importance needs to be given to the children who were HIV/AIDS infected. These girls also came up in support of a law or rule where both the husband and the wife were required to give a HIV –ve certificate before getting married. They felt that this itself would save a lot of innocent women from becoming victims of this disease.

Some of them agreed that films were a great cause for the spread of immorality in the Indian social setup. But few girls of the same class countered that by saying that India was the same from time immemorial and that even Kamasutra was only an Indian treatise on sex.

**92**

We interviewed 8 male students of classes 10, 11 and 12 from a leading private matriculation school in the city. At the outset they were very unwilling to spend time with us. Only with the persuasion of their teachers they grudgingly spent some time with us. They are highly motivated and their sole aim in life was to score a seat in medicine or engineering.

They are totally unaware of their surroundings and other activities. In fact, their daily routine is school, coaching class, and tuition. They are least bothered about the society. They refused to participate whole-heartedly in the discussion. They said that they suffered only from three diseases: cold, fever and nausea. All these were natural diseases.

There was no point in getting further opinion about AIDS from these boys. So we winded up the interview. We have recorded it here because we wish to say that such an extremely streamlined lifestyle unaware of the social consequences is not going to help.

Out of sheer ignorance these students might contract HIV. On the other hand, being non-sensitized they might discriminate other people ruthlessly.

**93**

We held discussions with the girl students of the same city school. These girls were more polite than the boys and showed greater sincerity than the boys. They were very interested in knowing more about HIV/AIDS. It was very captivating to see their interest in the well being of society.

All the girls said that a stigmatized disease of the present day society was HIV/AIDS. They said earlier it had been leprosy. Some girls felt that only CSWs get AIDS. When we said that even men get AIDS, these girls said that only men without character could get AIDS. They asked us a few questions which we answered. But they held the wrong notion that only poor people and people of Africa suffered from HIV/AIDS.

They denied having seen any advertisement in TV about HIV/AIDS. When we asked them about when AIDS day is celebrated, they did not know. They instead commented that all 365 days in the year are assigned for one cause or the other like Mothers Day, Fathers Day, Martyrs Day and so on. So it was not worth remembering.

The girls could not understand how the child born to a HIV+ mother can be HIV- through medical treatment. They felt that if a growing foetus can be saved, certainly a human being infected with the disease could also be saved.

**94**

As a next part of our interview, we met students from city corporation schools. Most of them were from poor or lower middle class society. The boys said that they knew about HIV/AIDS and visiting CSWs is the main reason for its spread. They recollected the ‘Pulli Raja’ advertisements to us. They repeatedly said that men with bad habits got this disease because of lose character and unwanted association with other women.

These students knew that the disease carried with it a bad social stigma. They said anyone with HIV/AIDS in their slum would be thrown out.

We could see that these boys from the slum knew very well about HIV/AIDS. But unlike the children of other schools they associated with it a very strong social stigma. They were not willing to sympathize that these men with HIV/AIDS they are also human and should be helped. In fact they said these men should be isolated from the society because they are scoundrels.

If they are treated with sympathy they would infect all the women in their colony. We asked them what they would do if one of them came to know that his father had been infected with HIV/AIDS? They said, “We will admit him in the hospital. We will see to it that he never visits our mother and infects her.”

Thus the standard of knowledge about HIV/AIDS among corporation school children of the city was best in comparison with other school children.

Perhaps the transparency of life in their homes and poverty had made them fully aware of the disease. The difference is social conditions was clearly reflected in their world-view.

**95**

These girls from the city corporation schools whom we interviewed can be broadly classified into two types. One, the girls who do not know anything about the disease or the ways in which it spreads. Two, the girls who are fully aware about all aspects of AIDS. The second set of girls came down heavily on CSWs who infect the male society.

They felt that it was a dirty topic to speak about.

None of them had the heart to sympathize or accept HIV/AIDS affected people under any condition . In their social setup it was a matter of prestige. Economically poor condition was acceptable, but not HIV/AIDS; ‘low’ caste was acceptable but not the social stigma due to HIV/AIDS.

**96**

Finally, we give the interviews of girls from convent schools. First we were not allowed to directly go and interact with the students like in other schools. We had to give a letter requesting permission and after a week it was granted by the nuns. The 45-minute discussion was a nice learning experience for us. Girls did not freely discuss the topic with us. They rather contented themselves in merely answering our questions. We kept on talking but they were reserved. Then we could understand that it was actually taught to these girls that to discuss sex-related topics was a taboo, so they had kept quiet.

Despite their seeming reserve, their awareness was really commendable. They had been taught about how it spreads. They had also been sensitized to the plight of the HIV/AIDS patients. They however have been educated into believing that HIV/AIDS patients are sinners. The teaching in the school is that these people suffer because they have sinned a lot. The girls did not mind serving people affected with HIV/AIDS. They even said that they would pray for these people. These girls said that the plight of orphans due to AIDS affected them very much. They felt that couples who don't have children should adopt only orphans due to AIDS. They were quite humane and we were pleased with their outlook.

## 97

We visited a matriculation school run by the Anglo-Indian group. The students there did not seem interested in discussing about HIV/AIDS. They looked at us as if had come to interview them because we did not have anything better to do. They said at the outset that only foolish people get AIDS and that intelligent people wear condoms. We could not refute their logic.

They did not talk about the ethics or discrimination. They said that HIV/AIDS was new to India. It was probably a foreign conspiracy that came into India because of the flesh trade. We talked to them about the rehabilitation of AIDS patients. One boy suggested that they could beg for alms in front of temples, mosques and churches without revealing their HIV+ status. We felt very uneasy to discuss further with these careless students.

Then we talked to them about children who are living with HIV/AIDS. They said, "Anyway these children are going to die. So why waste ones time on them. The sooner they die, the better it is."

However this counter-questioning irked some of their classmates. They started speaking up for the HIV/AIDS patients. It was nice to see controversial viewpoints. At the end of the discussion, the students' attitude towards AIDS was markedly different.

## 98

We interviewed the principal of a reputed central school in the city. He has put up nearly 25 years of service as a teacher. He is vexed with the students. He says that even if there are a few spoilt students, they succeed in ruining the whole class or school. Some boys have the habit of drinking. They introduce other children into these practices. They use many vulgar words frequently.

They are well versed in all bad practices like copying, lying and roaming. The teachers cannot punish these students severely. If they are suspended, they will approach the higher authorities and get the punishment revoked. Movies, magazines and computers help students get free sexual feedback and make them more corrupt. Teachers fear to give any form of sex education in the class because discipline cannot be maintained in a co-educational atmosphere. In the next decade the youth would be the highest population that is HIV/AIDS affected. This headmaster feels that parents are the root cause of such behaviour in children. Because parents have their own way, their children imitate them. No one can ever question them on anything. These boys are not even career oriented. By their natural intelligence they pass the exams. Majority of them don't even aspire for high marks. May be they are frustrated.

He feels the education system does nothing to make a student disciplined or have any morality. At least if they use their intelligence, they can save themselves from HIV/AIDS.

## 99

We interviewed the headmistress of a reputed school. She has 20 years teaching experience. We asked her suggestions to spread awareness about HIV/AIDS. She says her school children do not face any problem of that form. "Right from the day they join my school their only aim is to get very high marks and get a seat in IIT or in a medical college or at least a seat in a good engineering college. So my students have no time to do anything wrong.

Our students know the importance of marks and good education so they dare not waste their time. Majority of the children come from very good families so they have no chance to make mistakes or get HIV/AIDS." We then intervened and said, "we do not say that your school children would get HIV/AIDS we only asked your suggestions of AIDS awareness."

She said, "If I talk these things either to the students or to the parents, I would be first kicked out of the chair. The administration expects only standard, good education. So our students get good marks and are state rank holders. The discipline of the student is never a problem. Only a teacher with no stuff will face problems because she will be harrowed and made to weep by these children. We do not have the work of maintaining discipline among students. There is no time to speak of HIV/AIDS in the school premises. You have to seek counseling and advice only students from ill disciplined schools.

School like ours do not need any such awareness lectures or talks. My students will never waste their time before TV or in a cinema theatre." Thus our interview with her made us think it is inappropriate to talk about HIV/AIDS with her. In fact children from her school were better than the headmistress when we talked about HIV/AIDS with them. This woman behaved as if she belongs to a different planet.

## 100

We interviewed a nun who is the headmistress of a girls' convent school. Her school employs only female teachers. She takes pride in the discipline prevailing in the campus. In her school, HIV/AIDS awareness education is given to high-school children.

They celebrate AIDS days on 1 December of every year. There are two moral education classes every week where topics like sex education and HIV/AIDS are dealt with. This prepares her students not only in the syllabus but also more importantly they are taught overall value education.

She says they have tried to give the students the best of real education. Only 2 to 3 percent of these students go for professional courses. Many of them study up to graduation. A few opt for teacher training and nurse courses. She sadly adds that if she or her team of students visit hospitals to help these HIV/AIDS patients, the orthodox Hindus think that they are up to conversion.

He argues, 'What is the rhyme or reason for us to convert HIV/AIDS patients whose lifespan is short and many of whom are terminally ill?' Such is the thought of some people. After all when we pay visit we give them some consolation, we pray for their wellbeing and help them with some mental strength to bear with the disease.

Some of the people living with HIV/AIDS are orphans or were abandoned on the platforms. A few of them have been joined to a hospital when groups like us found them abandoned in streets. Many of them we picked were women.

She recalls how selflessly Mother Theresa worked for leprosy patients. Earlier, leprosy was a dreaded disease like HIV/AIDS.

Now leprosy patients stay at home and take treatment. It was the work of the missionaries to have removed the stigma of leprosy. In those days, even if a rich man got leprosy, he was chased out of the house. The only means by which he could sustain himself was by begging. This plight has now changed. Likewise, people living with HIV/AIDS must be comforted.

The nation must ensure that they lead a life of comfort and dignity. Everyone must cooperate with this. The social stigma associated with HIV/AIDS is so high that unless some collective steps are taken to help them it is impossible for the affected people to lead a peaceful life.

## 101

We met a driver employed in the state transport corporation. He is from Tiruthani but presently settled in Chennai. He was in his late forties. He has studied up to SSLC. He is married and has two children (one son and one daughter). His marriage was an arranged one. He says he earns Rs.6000 p.m. He lives in a rented house. He claims to know a lot about HIV/AIDS. He says the disease mainly comes from 'ladies' (CSWs). He is against sex with women who have 'loose' character and says that everyone should be very careful. He has no such relationships. He remarks, 'what can one say about CSWs? When people have no other job, they visit CSWs. These poor ladies take up this trade because they have no other option.'

He suggests that people living with HIV/AIDS must not be neglected or ignored. We can advice them even if we are not in a position to help them monetarily. He claims he does not have any bad habits.

He does not see movies. It is 5 or 6 years since he had been to a cinema. He saw a porn movie only once in his life, when he was 20-years-old. He has not read pornographic books since he feels it is wrong to read or see such books. First one would read them, and then perhaps do these things. So he doesn't read them. If he has an opportunity to see a HIV/AIDS patient he would touch them only if they are not visibly infected. He fears the disease and mingling with people who are infected by it.

Chapter Five

# CONCLUSIONS

The study of psychological problems of HIV/AIDS patients is directly related to public opinion. The depression that the people living with HIV/AIDS suffer is due to the social stigma associated with the disease and the concept of immorality. Usually, people living with HIV/AIDS are segregated, ostracized and discriminated. Instead of being treated with understanding and sympathy, they are made to suffer mental agony.

We felt it very important to know the public opinion. The interviewees consist of students, police, teachers, uneducated workers, businessmen, women, and staff from government and private offices. We have selected them arbitrarily. And all the interviewees were people whom we met for the first time.

When we sought the interviews most of the people looked suspiciously at us, some shouted at us, some ran away from us, only a few answered our questions. It is all the more important to note that when certain questions were put forth they refused to answer us and treated us with contempt and anger. The overall view is that even the utterance of the word AIDS seems to be held with disgust.

Almost everybody we interviewed felt shy to talk about AIDS. Everybody spoke to us on the condition of anonymity and confidentiality. We have only mentioned their profession, educational qualification, economic status and age to give an idea of the different viewpoints of different sections of society and some of them did not mind giving their names so in such cases we have also used their names.

Now we give the conclusions based on the interviews, suggestion made by public, the conclusions and suggestion derived from the mathematical analysis.

1. A vast majority declined to treat HIV/AIDS affected people with sympathy or due respect. They argued that patients had to suffer for the crimes that they had committed. The social conscience of our society is so strong. On the other hand, young men and women were not so biased against people living with HIV/AIDS. They expressed their willingness to meet and mingle with them.

2. Everybody pities HIV/AIDS infected children. Some people were willing to help them with money. Many people openly declined to adopt these orphaned HIV/AIDS children but suggested that the government must support these children with free healthcare, education and food. Their pity and understanding originated because the children were seen as 'sinless' victims. The fact that the public is overwhelmingly in support of children shows that there is enormous space for interventions to eradicate the negative opinions held about AIDS.

3. Everybody is sympathetic of women who have lived monogamously and have become patients for no fault of theirs but because of their erring male partners.

4. People felt that initiatives for change must come from government organizations.

5. They urged women's organizations and women self-help groups to take up these issues and help ailing women.

6. Those who were interviewed had varied opinions about CSWs. Many of them did not seek a ban on CSWs. They all wanted CSWs to be rehabilitated and given gainful employment. Many people (especially women) pitied the CSWs who have taken up the trade due to poverty or because of being deserted by their husbands. Some felt that women take up this trade to become rich in a short span of time. The extras and small time actresses do this job to maintain their social status. Commercial sex work is carried out through Internet and cell phones.

7. The students who were interviewed did not shun speaking about HIV/AIDS. They are full of enthusiasm to work for these patients. While other people avoided these topics, students were eager to spread awareness about the disease. They even volunteered to teach children affected with HIV/AIDS.

8. Many students suggested that fund collections for AIDS can be done through installing *hundis* in academic campuses and other public areas. This will spread awareness and also get some money that can be used for the welfare of AIDS patients.

9. Students felt that a single day to observe HIV/AIDS was not sufficient. They wanted at least a week's celebration for spreading AIDS awareness. Student organizations like the NCC and the NSS can be roped in, to spread awareness about AIDS in rural areas through camps organized there. Also, they can be made to visit care homes for HIV/AIDS affected. This would sensitize them to the sufferings of people living with HIV/AIDS.

10. Glamour and vulgarity in films has been attributed as the sole reason for the lowering of public morals. People across different age groups put forth this opinion. Lewd dialogues and obscene song sequences kindle the sexual emotions of people and lead them astray.

11. One of the interviewees suggested the Cuban model where HIV/AIDS patients are quarantined and kept in separate care homes. This has ensured that in Cuba the spread of HIV/AIDS is zero.

12. Cheap availability of alcohol was one of the reasons for accelerating the spread of HIV/AIDS in country. People feel that when men are uninhibited and in a drunken stupor, they invariably follow their friends and end up visiting CSWs.

13. People say that men who know that they are HIV+ keep visiting CSWs. This is one of the reasons for the alarming rise of HIV/AIDS in Tamil Nadu. Our interviews with the ayahs and ward-boys of the Tambaram Sanatorium clearly proved that men who were inpatients and HIV+ kept visiting CSWs in the nights.

14. HIV/AIDS patients in rural areas face social stigma, harassment and discrimination more acutely than in the cities. They are ostracized and in some cases, even chased out of the village. Likewise even awareness about HIV/AIDS is also much lower in rural areas compared to the cities. In a city, a person with HIV/AIDS can survive without identifying himself as a patient. But in a village such things immediately come to light because of the closely-knit social setup.

15. Religious groups have a major role to play in disseminating AIDS awareness. People congregate in large numbers only in the name of religion. If religious leaders take up the responsibility to talk about AIDS, certainly it will influence public opinion.

16. Relatives of the patient living with HIV/AIDS must be counseled to accept him/ her and to give all possible support so that he / she can lead a life of peace. If everybody started tending their relative who had HIV/AIDS, instead of abandoning him/her it would certainly be a positive example in society.

17. Most of the people feel very shy and ashamed to talk about HIV/AIDS in any public place. They seek a secluded, lonely place to discuss with us or to answer the questionnaire about HIV/AIDS. They did not like even to fill up our questionnaire. When we offered to collect the questionnaire later, and asked them to fill it at home they were literally angered to the extent of manhandling us. They preferred us to ask questions which they answered.

18. Some of them held a very strong view that they HIV/AIDS infected men can be isolated or even imprisoned. They can by all means given good food and medical aid but not allowed to visit CSWs or their wife. Unless the government takes strong steps it is impossible to control the spread of HIV/AIDS.

19. A handful of the people feel that proper food habits and protein rich will prolong the lifespan of HIV/AIDS infected person. Government should see to it that rural poor get this type of food.

20. The interviewed people pointed out that since a vast majority cannot read, AIDS awareness advertisements should not be restricted to posters, but must be given through TV and radio.

21. Several people said that CSWs were present all over Chennai (since there was no specific red-light area) and this in turn was dangerous because no awareness program could be carried out successfully. They were also aware of the existence of pimps. They decried the police who are not able to successfully eliminate the flesh trade in the city.

22. Those who had known the behaviour of male HIV/AIDS patients felt that few of them had turned vengeful in life and wanted to infect as many women as possible because they felt betrayed, that only a woman had infected them.

23. Some people said that unless self-realization comes to everybody it is impossible to stop the spread of HIV/AIDS.

24. Students themselves acknowledged that the mass media like TV and films had a major role in youth adopting high-risk

behavior, which in turn could make them HIV/AIDS infected. Everybody shared this view.

25. Most of them welcomed free HIV/AIDS test, with a request that the results of the test must be kept confidential.

26. One interviewee suggested that the government could introduce a miniscule percent of reservation for HIV/AIDS patients in public sector jobs. This will give them an enhanced status in society and bring down the levels of discrimination.

27. Government should arrange for doctors/ sociologists who have expertise in HIV/AIDS to address employees of police departments, government employees, private and quasi government employees so that awareness levels increase and the social stigma associated with AIDS goes away. Similar endeavors have to take place in the private sector too.

28. Adult education has gained popularity in rural areas. It should also include HIV/AIDS awareness education.

29. Some people came down heavily on the concept of free condom distribution. They said that it promotes promiscuous behaviour and visits to CSWs. They compared it to the sale of cigarettes containing the statutory warning "Smoking is injurious to health."

30. Everybody felt that migration was the root cause for rural uneducated people getting infected by HIV/AIDS!

31. Some people wanted to ban CSWs. Almost everybody accepted that these women must be rehabilitated. Only one man said that CSWs were essential for the maintenance of the social structure.

32. Many people wanted scientists to discover a medicine to cure HIV/AIDS. Others differed in their opinions. They felt that a medicine would foster greater immorality in society.

33. People unanimously felt that self-control and good conduct would certainly prevent AIDS.

34. With the exception of five people, all of them felt that prayer couldn't cure this disease. Whereas majority of the HIV/AIDS patients felt prayer can cure the disease.

35. It became clear from our interviews that people discuss HIV/AIDS with their friends, but they never discuss it with their family. Because of the taboo surrounding discussion of sexuality, no discussion of HIV/AIDS takes place in most homes. Efforts must be taken to make the topic discussable in family surroundings, because it is one way in which awareness can spread and discriminatory perspectives can be changed.

36. Some people said that if good movies are made about HIV/AIDS in which the patients are treated very well by the society and with kindness by the doctors certainly it would remove the stigma associated with HIV/AIDS.

37. A few of them said that AIDS like TB and leprosy made the patient lose his appearance because of being covered with scabies, sores and boils.

38. Some people felt that not all cases of HIV/AIDS spread through unsafe sex, but could also occur due to blood transfusion and sharing needles by injecting drug users.

39. The advertisement strategies need to be changed. It has become public knowledge that AIDS has no cure. This is one of the reason why it is dreaded so much. Other diseases like blood pressure and diabetes don't have cure. But people only view them as livable diseases, and not as killer diseases. There are instances of people who have lived for more than two decades with seropositivity. AIDS must be treated on par with these diseases. Only then the stigma can be removed.

40. Most of the interviewed persons felt that if a person suffers from HIV/AIDS he/ she need not give information about his/ her seropositive status to the employer because they could be ill-treated or fired. A few felt that it is the ardent duty of the employee to inform their HIV/AIDS status to their employer.

41. Pension must be paid to CSWs so that this will allow them to lead a life of self-respect.

42. Some people recollected incidents where doctors had denied even first aid to HIV/AIDS patients. Most doctors with private practice do not treat HIV/AIDS patients. They feel that if everybody comes to know that an HIV/AIDS patient is taking treatment in their clinic, other clients will boycott the clinic. On the other hand, there are so many incidents where doctors personally visit the homes of rich people with HIV/AIDS and administer medicine to them and maintain strict confidentiality. This kind of dual behaviour clearly indicates the hypocrisy prevalent in our society.

43. Schools refused to admit children who are affected with HIV/AIDS if this information is revealed to the authorities. Such discrimination must be brought to public attention and there has to be sufficient reaction against it. Only this can improve the state of affairs and reform people who are vindictive of AIDS patients.

44. We came to know that there were several misconceptions among people. These have become myths and popular beliefs. For instance it is a prevalent belief that sex with very young virgins can cure HIV/AIDS. We encountered several such myths during the course of our interview. Some of them had the misconception that HIV/AIDS spreads by air, touch and contact. One man was vehement that it spreads through mosquito bites. In no way could these people be convinced that their view was wrong.

45. It is high time people living with HIV/AIDS rally together and work for their own empowerment. This may give them some power in the nation. Unless they are empowered, nothing can be done to help them from being treated badly. They should emerge with courage to ask for their right. It is pertinent to mention here that even human rights organizations are very quiet over the plight and problems of HIV/AIDS infected patients.

46. People in the age group 20-24 did not show much interest in the HIV/AIDS disease because they had more immediate concerns such as education, or finding employment. This could be seen from the CETD matrix.

47. People in the age group 50-54 showed absolutely no interest. They were engrossed with their own problem of marrying off their children, building homes and making post-retirement plans.

48. The CETD Matrix graph with 6 intervals of age group□ 60, 50-59, 40-49, 30-39, 20-29, □ 19 highlighted a lot of facts. People in the age group □ 19 and 20-29 have an aversion to discuss about HIV/AIDS awareness program. The age group 30 to 49 showed steady interest in the increasing direction and the age group 50-59 showed a sudden decrease i.e., maximum aversion to speak or even discuss about the problem of HIV/AIDS. This is clearly brought out in the graph 3.1.

49. In order to make the study more sensitive we divided the age group from 6 intervals to 11 intervals, □ 65, 60-64, 55-59, …, 20-24, 14-19. The CETD graph reads as follows: Maximum repulsion to discuss such a tabooed topic is shown in the age group 60-64. There is some aversion by the age groups 20-24 and □ 65. Apart from this, the age groups 14-19 and the age group in the range 25 to 59 show steady interest with very small fluctuations. However the age group 50-54 is passive.

50. From the AFTDR matrix we can observe that all the age groups show aversions and some concern but they equally show fear associate a social stigma and are shy about it. However they are sympathetic. Thus the reactions across all age groups is mixed. Other concepts never come to the ON state at all which shows the strong feelings they hold towards the disease.

51. Even after the refinement of the age group they had the same concepts thereby making it very clear that urban people are at large not so averse with patients or AIDS but at the same time they did not show any strong positive feelings.

52. They were all sympathetic and concerned about orphaned children and children infected with HIV/AIDS and women infected by their husbands.

53. Medical and engineering seats should be reserved to children whose parents are infected by HIV/AIDS. This will have two advantages:

    a. The social stigma would be reduced.
    b. People, especially the rich and middle-classes, may have less anxiety to reveal and accept their HIV/AIDS status.

54. Free education should be provided to children whose father/ mother is a HIV/AIDS patient.

Appendix one

# QUESTIONNAIRE

## MANUAL FOR FIELD RESEARCHERS

The points to be had in mind while interviewing the public.

- Names of the Interviewee, will not be disclosed and their identities shall be protected.
- The interviews of the public shall be used by us purely for research purposes.
- The interviewee shall be let to speak out freely and fully. The interviewer shall not intercept.
- The interviewers shall select their questions according to the circumstances/ persons.
- Answers to questions with '*' will be graded.
- When these questions are found not adequate, more questions can be framed as per the requirements or depending upon the situation.
- The interview shall be opened and concluded depending upon how openly the interviewee wishes to speak out.
- When situations arise in which the interviewee conceal certain information / misrepresent, the field researchers shall adopt some psychological techniques and note them down separately.

Interviewee :
Date and place of interview :

**GENERAL PARTICULARS**

1. Name :
2. Age :
3. Sex :
4. Educational qualification :
5. Native place and Address :
6. Present Address :
7. Caste :
8. Religion :

**JOB-RELATED INFORMATION**

9. Profession :
10. The office in which you worked/are working :
11. How many persons are employed in the office? :
12. Name of the profession :
13. At what age did you join the profession? :
14. Whether the job is permanent or temporary ? :
15. What skills do you know? :
16. Do you do this job hereditarily? :
17. Have you any other sources of income? :
18. How much do you earn as monthly income? :
19. For how many years you are living in Chennai? :
20. Have you gone abroad? If yes,the name of country? The purpose of your trip? The income you earned? :
21. Whichever places you went

from your native place
to earn money? :

22. Do you feel the AIDS status must be disclosed to :

- ❑ Friends
- ❑ Siblings (brothers/sisters)
- ❑ Children
- ❑ Other relatives
- ❑ Boss
- ❑ None
- ❑ In-laws
- ❑ Spouse
- ❑ Strangers
- ❑ All

23. Have you talked about AIDS disease with friends, family and fellow employees? (or, it is not necessary?) :

**MARRIAGE**

24. Are you married? :
25. Your age at marriage :
26. Love or arranged marriage :

**HUSBAND/WIFE**

27. Name :
28. Educational qualification :
29. Age :
30. Occupation :
31. Govt./private sector job? :
32. Is husband/wife alive? :
    If no, did you remarry? :

**CHILDREN**

33. How many children? :
34. Up to what standard the children have studied? :
35. What jobs the children are doing? (If over 18) :

36. Have you ever discussed about HIV/AIDS with your children? :

37. How did your children behave when you discuss in general about AIDS disease? :

- ❑ Cold and distant
- ❑ Kind
- ❑ Fear
- ❑ Other reaction
- ❑ No concern

38. Do your children have knowledge about HIV/AIDS?:

39. Do you feel children affected with HIV/AIDS must be helped? :

40. What you expect the Govt. to do for those children infected with HIV/AIDS? :

41. Whoever, and in whichever manner, can AIDS education be imparted to children? :

- ❑ Advertisements in TV
- ❑ Advertisements in Paper
- ❑ Dramas
- ❑ Teaching in School
- ❑ Songs
- ❑ Posters
- ❑ Radio-broadcast
- ❑ Others

42. What do you think about giving education to the AIDS–affected children together with others? :

- ❑ No
- ❑ No reaction
- ❑ Sign on the face
- ❑ Yes
- ❑ Never
- ❑ Other reaction

43. Do you have any close relative affected with HIV/AIDS? :

44. How do you feel for them? :

- ❑ Cold
- ❑ Scornful
- ❑ As usual
- ❑ Dejected
- ❑ Sad
- ❑ No reaction
- ❑ Blaming
- ❑ Depressed
- ❑ Aversion
- ❑ Others

45. Is it own house or rented? :
46. What is the nature of the house you live? :
47. If own house, are there any tenants living in it? :
48. If rented house what rent are you paying? :
49. If apartment whether government housing or private buildings? :
50. Whether nuclear family or a joint family? :
51. Are you living in your ancestral home? :
52. If joint family, with whom are you living? :
53. What type of locality is your house situated? :
54. Whether any of your friends have this disease? :
55. How many of your friends have contracted AIDS? :
56. What is their condition? :
57. Have you ever spoken with your friends / colleagues about AIDS? :

- ❑ Never
- ❑ Fully
- ❑ Aggressive
- ❑ Despair
- ❑ Dejection
- ❑ Fear
- ❑ Sometimes
- ❑ Silent
- ❑ Pathetic
- ❑ Anger
- ❑ Others

**COMMERCIAL SEX WORKERS**

58. Can you tell the reasons for seeking improper sex partner? :

- ❑ Urge for sex
- ❑ Not satisfied
- ❑ Male Ego
- ❑ Silent
- ❑ Doesn't like to answer
- ❑ No reaction
- ❑ Looks sad
- ❑ Others
- ❑ Quarrel at home
- ❑ Others

59. What do you think about the living standard of the CSW? :

- ❑ Forced into the trade
- ❑ Destitute
- ❑ Poverty
- ❑ Family burden
- ❑ As per liking
- ❑ For money
- ❑ Others
- ❑ No answer

**GENERAL QUESTIONS**

60. Have you ever talked with your family members regarding AIDS? :
61. Are there AIDS patients in your family/relations? :
62. Do you know anybody, known to you who has HIV?:
63. How your family members/ relations/friends should behave in taking care of the AIDS patients? :

- ❑ Kind
- ❑ Harsh
- ❑ Open
- ❑ Considerate

64.* What is your view on the opinion that if medicine is discovered to fully cure

AIDS, social deterioration will further be worsened? :

65.* What is your opinion on the concept that the government should call out all the AIDS patients for free medical examination/checkup? :

66.* What do you think about those who say: AIDS will not touch me. Don't talk about it to me? :

67.* What is your view on the statement that "it is the social duty of everyone to help AIDS patients?" :

68. What do you think about giving sex education in school? :

69. What do you think about yellow journalism and pornographic/blue films? Does it have any relevance with the spread of HIV/AIDS? :

70. Different views prevail in different countries on the use of contraceptives. What is your opinion? Can condom help in the prevention of HIV/AIDS? :

71. Whether one has or has not the right to know whether his colleagues in his workplace has AIDS? :

72.* It is natural to have a fear while speaking, reading and commenting about AIDS. Is there a way to change this attitude? :

73.* Have you had a discussion with your friends on AIDS like cancer or TB? :

74. About what you would talk during such discussion? :

75.* What is your mental state on seeing an AIDS patient? :

76. Will you explain to the public about AIDS if you are given an opportunity to do so? :

77. Many people speak differently about AIDS.You might have heard/spoken to some of them; Can you recall any one unforgettable incident related with HIV/ AIDS patients or knowing more about it or any other situation? :

78. Is the AIDS patient alone responsible for contracting the disease? :

- ❑ Society
- ❑ Wife / Husband
- ❑ Government
- ❑ CSWs
- ❑ Frustration / depression

79. “When anyone sees/hears an AIDS patient, s/he is emotionally surcharged. S/he is subjected to mental affliction/grief/ anger/insult/frustration and behaves so. Can you tell the reason for such behaviour? :

80. “The foremost duty of the government is to ascertain how many citizens are suffering from AIDS. But no such precise assessment has been made in India.” What is your opinion on this statement? :

- ❑ Can the data from the rich and upper middle class be collected?
- ❑ Is the data pertaining only to the poor and rural middle class?

81. Why no free medical test is not being conducted on all citizens to find out whether anyone has AIDS like for smallpox/ diabetes/leprosy? :

82.* We talk about diabetes, leprosy and cancer in all places. But, we do not talk about AIDS. Why? :

83. What do you think about those who say "Increasing of AIDS disease is not my problem. It is the problem of concerned persons"? :

84.* "Suppose when you come to know your friend has AIDS how will you move with him"? Why? :

85. What expenses you want the government to meet for AIDS awareness? :

86. What should be banned to prevent spread of HIV/AIDS?

| | |
|---|---|
| ❑ CSWs | ❑ Bars |
| ❑ Internet-cafe | ❑ Cinema |
| ❑ Liquor shops | ❑ Porn Magazines |
| ❑ TV serials | ❑ Porn movies |
| ❑ Smoking | |

87. What further measures do you expect the government to take to control AIDS? :

| | |
|---|---|
| ❑ Advertisement | ❑ Condom promotion |
| ❑ Awareness education | |

88. Which of the measures the government is undertaking to prevent AIDS is most attractive to you? :

89. In which medical system you have faith that can cure HIV/AIDS : Allopathy/Siddha/ Ayurveda/Unani?
The reason for it? :

90. Do you fear that AIDS will spread through mosquito biting/kissing on hand/lips and eating in same plate? :

91. Which advertisements for AIDS awareness is mostly liked / disliked by you? :

- ❑ Pulliraja
- ❑ When you have wife why go for CSWs?

92. What do you think when you see in private the advertisements on AIDS? :

| | |
|---|---|
| ❑ Fear | ❑ Guilt |
| ❑ Unconcerned | ❑ Sadness |

93. What do you think when you see advertisements about HIV/AIDS in the company of your family members? :

| | |
|---|---|
| ❑ Fear | ❑ Guilt |
| ❑ Unconcerned | ❑ Sadness |

94. What do you think of AIDS awareness advertisements?:

| | |
|---|---|
| ❑ not sufficient | ❑ Sufficient |
| ❑ not of standard | ❑ not of good taste |

95. Can you suggest good advertisement techniques for AIDS awareness/ publicity? :

- ❑ Drama
- ❑ Skit/ Street thereat
- ❑ Cinema
- ❑ Shows by top actress
- ❑ Address by political leaders

96. What further steps the central/state governments must take to eradicate AIDS?:

97. Should all people agree to undergo AIDS test if the government comes forward to do the test freely? :

- ❑ Yes
- ❑ No
- ❑ No response

98. What kind of expense you want the government to bear for the AIDS control? :

- ❑ Medical
- ❑ Money
- ❑ Supply of good food
- ❑ Free treatment

99. Do you think that the AIDS patients should/shouldn't be employed in govt./pvt. services? :

100. "Any disease can be contracted in the world. But the AIDS should never come" what do you think about this? :

101.* "AIDS patients should be given top priority in medical treatment above all other patients". What is your view on this opinion? :

102. Do you think that the AIDS patient should also be present when the doctors are discussing about the nature of AIDS? :

103.*What do you think about interviewing you by people like us about AIDS awareness? :

104. What other questions you want us to ask you, other than those we asked you? :

105. Please tell what part of our questions, you like/don't like? :

**ADVERTISEMENTS**

106.*Through which source you came to know about AIDS? :

107. If you came to know it through advertisements, what type was the media? :

108. Our government has provided telephone facilities in several places throughout the nation to tell confidentially doubts questions therein. Have you used this facility? Do you know about this facility? :

109. By whatever further steps the government can take the advts. so that it reaches to the masses? :

110.* What is your view about fully eradicating prostitution from our country? :

111. How can the NGO/social organizations make the advertisements on AIDS reach the public? :

112. In what ways can your religion publicize AIDS? :

113. Do you feel school and colleges do enough to spread awareness about HIV/AIDS? :

114. Do you accept the participation of students in spread of awareness about HIV/AIDS? :

115. Does Hindu religious leaders ever take steps to spread about awareness about HIV/AIDS? :

116. Can you give any good awareness program about HIV/AIDS? :

**RELIGION**

117. Do you profess the same religion hereditarily? If not, in which generation you converted to this religion? :

118. Do you think conversion can help HIV/AIDS patients? :

119. What reasons you can give for your conversion about HIV/AIDS patients? :

120. Do you think HIV/AIDS patients can visit religious institutions? :

121. Can faith in god cure HIV/AIDS? :

122. Is the traditional set-up of India a cause for the social stigma the disease holds? :

123. Is it ever possible to wipe out the social stigma associated with the disease? Your suggestions! :

Appendix Two

# TABLE OF STATISTICS

| NO | SEX | AGE | PROFESSION | EDUCATIONAL QUALIFICATION | SOCIAL STATUS | MARITAL STATUS |
|---|---|---|---|---|---|---|
| 1 | M | 53 | Professor (Arts College) | Doctorate | Rich | Yes |
| 2 | F | 49 | House wife | Graduate | Rich | Yes |
| 3 | M | 35 | Asst. Professor in Engg College | (Doctorate) | Middle class | Yes |
| 4 | M | 68 | Provision shop owner | 9th std | Upper middle class | Yes |
| 5 | M | Around 40 | LDC in Secretariat | Arts Graduate B.A. | Middle class | Yes |
| 6 | F | 50 | House wife, earns by tailoring | Uneducated | Middle class | Yes |
| 7 | M | 19 | Auto driver | 9th std failed | Poor | No |
| 8 | F | Late forties | Vegetable vendor | Uneducated | Middle class | Yes |
| 9 | M | 21 | Post graduate student studying | MA studying | Upper middle class | No |
| 10 | M | 45 | Barber | School dropout | Upper middle class | Yes |
| 11 | M | Late thirties | Police constable | Art graduate, BA | Middle Class | Yes |
| 12 | F | 22 | Housewife | 6th std | Poor | Yes |
| 13 | F | 22 | Unemployed | Science graduate | Middle class | No |
| 14 | F | 21 | Studying | BSMS (studying) | Rich | No |

| NO | SEX | AGE | PROFESSION | EDUCATIONAL QUALIFICATION | SOCIAL STATUS | MARITAL STATUS |
|---|---|---|---|---|---|---|
| 15 | F | 26 | Preparing for IAS | Post graduate | Middle class | No |
| 16 | M | 46 | Vegetable vendor | 10th std | Poor | Yes |
| 17 | M | 43 | Petti shop owner | Graduate | Poor | Yes |
| 18 | F | 22 | Working in a shoe Company | 7th std | Poor | No |
| 19 | M | Late forties | Mobile tailor | Did not disclose | Middle class | Yes |
| 20 | M | 27 | Commando squad | Graduate | Middle class | |
| 21 | M | 38 | Upholstery company | 10th std | Poor | Yes |
| 22 | M | 22 | Unemployed | 10th std | Poor | No |
| 23 | M | 32 | Police department | 12th std | Middle class | Yes |
| 24 | M | 68 | Security working in a private company | 12th std | Middle class | Yes |
| 25 | M | 52 | Porter | 5th std | Poor | Yes |
| 26 | M | 32 | Working in a fancy store (before a driver) | 10th std | Middle class | Yes |
| 27 | M | 45 | Self employed owns a electronic repair shop | 10th std | Middle class | Yes |
| 28 | M | 19 | Unemployed | 12th std | Poor | No |
| 29 | M | 26 | Computer Artist | 10th std | Middle class | No |
| 30 | M | 24 | Watch repair - shop | 10th std | Middle class | Yes |
| 31 | M | 18 | Working in a fancy shop. | 8th std | Poor | No |

| NO | SEX | AGE | PROFESSION | EDUCATIONAL QUALIFICATION | SOCIAL STATUS | MARITAL STATUS |
|---|---|---|---|---|---|---|
| 32 | M | 59 | Works in Madras Corporation | $10^{th}$ std | Middle class | Yes |
| 33 | M | 24 | Lottery shop & STD booth owner | Did not disclose | Middle class | No |
| 34 | M | 35 | Cassette shop owner | | Middle class | Yes |
| 35 | M | 24 | Business | Plus two | Middle class | No |
| 36 | M | 26 | Cloth merchant | - | Rich | Yes |
| 37 | M | 23 | Fish market | $7^{th}$ std | Poor | No. |
| 38 | M | 56 | Motor mechanic | Diploma Holder | Middle class | Yes |
| 39 | M | 29 | Works in teashop | Plus two | Poor | Yes |
| 40 | M | 30 | Platform shopkeeper | $9^{th}$ std | Poor | Yes |
| 41 | M | 38 | Government Employee | Graduate | Middle class | Yes |
| 42 | M | 25 | Workers in an Auto company | ITI Diploma | Middle class | Yes |
| 43 | M | 30 | Motor mechanic | $7^{th}$ std | Middle class | Yes |
| 44 | M | Late forties | Rice merchant | $5^{th}$ std | Middle class | Yes |
| 45 | M | Late twenties | Centering | – | Middle class | Yes |
| 46 | M | Late thirties | Salesman | – | Poor | No |
| 47 | M | 29 | Transport PT dept. | Graduate | Middle class | Yes |
| 48 | M | 22 | Daily wager in a teashop | $10^{th}$ std. | Poor | No |
| 49 | M | 23 | Watch repairer | $9^{th}$ std | Middle class | No |

| NO | SEX | AGE | PROFESSION | EDUCATIONAL QUALIFICATION | SOCIAL STATUS | MARITAL STATUS |
|---|---|---|---|---|---|---|
| 50 | M | 56 | Construction labourer | $4^{th}$ std | Middle class | Yes |
| 51 | M | Late thirties | Works in a company | ITI | Middle class | Yes |
| 52 | M | 33 | – | Science graduate | Middle class | Yes |
| 53 | M | 28 | Carpenter | $10^{th}$ std | Middle class | Yes |
| 54 | M | 49 | Barber | $8^{th}$ std | Middle class | Yes |
| 55 | M | Early thirties | Selling petty articles | Uneducated | Middle class | Yes |
| 56 | M | 29 | Working in Navy | $10^{th}$ std | Middle class | No |
| 57 | F | 22 | – | Graduate | poor | No |
| 58 | F | 26 | Works in a private firm | MBA | Middle in come group | No |
| 59 | F | 29 | House wife | $6^{th}$ std | Middle income group | No |
| 60 | M | 32 | Accountant / assistant in a Jain shop | Plus two | Middle class | Unmarried<br>Unmarried |
| 61 | M | 26 | Unemployed | Art graduate | Middle class | – |
| 62 | M | 22 | Works in the liquor shop owned by father | B.com., | Rich | Unmarried |
| 63 | M | 23 | Unemployed | $8^{th}$ std | Middle class | Unmarried |
| 64 | M | 47 | Owns a petty shop | Not known | Middle class | Not known |
| 65 | M | 25 | Police constable | Plus two | Middle class | Unmarried |
| 66 | F | 19 | Studying | B.Sc., studying | Middle class | Unmarried |
| 67 | F | 22 | Self employed making wire bags | $8^{th}$ std | Poor | Unmarried |

| NO | SEX | AGE | PROFESSION | EDUCATIONAL QUALIFICATION | SOCIAL STATUS | MARITAL STATUS |
|---|---|---|---|---|---|---|
| 68 | F | 24 | Unemployed | MCA | Middle class | Unmarried |
| 69 | M | 14 | Works in a mechanic shop | 6th std | Poor | Unmarried |
| 70 | F | 19 | Studying 2nd year MBBS | 2nd year MBBS | Middle class | Unmarried |
| 71 | F | Early thirtees | Women police Constable | Don't know | Middle class | Unmarried |
| 72 | M | Early forties | P.T. master in Adyar | P.T trained | Middle class | Not known |
| 73 | M | Early forties | Pastor | Not known | Middle class | Not married |
| 74 | F | 28 | Working for PhD in Biomedical sciences | PhD studying | Middle class | Unmarried |
| 75 | M | 72 | Sannar one who climbs the coconut & palm trees & cleans it | Uneducated | Middle class | Married |
| 76 | M | 51 | Head cook | Art graduate | Middle class | Not known |
| 77 | M | 25 | Works in are internet café | Diploma in computer Applications | Poor | Unmarried |
| 78 | F | Early fifties | House wife | Graduate | Upper middle class | Married |
| 79 | M | Late twenties | Waiter in a city hotel | Science graduate | Middle class | Unmarried |
| 80 | M | Early fifties | Social worker | Graduate from MCC | Middle class | Yes |
| 81 | M | Late fifties | NGO | Unknown | Middle class | Yes |
| 82 | M | Late fifties | Doctor | MBBS | Rich | Yes |

| NO | SEX | AGE | PROFESSION | EDUCATIONAL QUALIFICATION | SOCIAL STATUS | MARITAL STATUS |
|---|---|---|---|---|---|---|
| 83 | M | Early twenties | Ward boy in a Tambaram Hospital | Uneducated | Poor | No |
| 84 | M | 36 | Rag Picker | Primary school dropout | Poor | Yes |
| 85 | M | 65 | Porter in the central Railway station | - | Middle class | Yes |
| 86 | M | Late twenties | Male Nurse | Nursing | Middle class | Yes |
| 87 | M | Late forties | Electrician | Diploma holder | Middle class | Yes |
| 88 | M | Late sixties | Cobbler | Uneducated | Poor | Yes |
| 89 | M | Early fifties | Post master | - | Middle class | Yes |
| 90 | M | Late teens | High school students | Studying 11th or 12th | Middle class | No |
| 91 | F | Late teens | High school students | Studying 11th or 12th | Middle class | No |
| 92 | M | Late teens | Very good school | Studying 11th or 12th | Rich | No |
| 93 | F | Late teens | Very good school | Studying 11th or 12th | Rich | No |
| 94 | M | Late teens | Corporation school | Studying 11th or 12th | Poor | No |
| 95 | F | Late teens | Corporation school | Studying 11th or 12th | Middle class | No |
| 96 | F | In late teens | Convent | Studying | Rich | No |
| 97 | M | Late teens | Anglo Indian school | Studying 11th or 12th std | Rich | No |

| NO | SEX | AGE | PROFESSION | EDUCATIONAL QUALIFICATION | SOCIAL STATUS | MARITAL STATUS |
|---|---|---|---|---|---|---|
| 98 | M | Late fifties | Principal central school | Post graduate | Rich | Yes |
| 99 | F | Early forties | HM | Post graduate | Rich | Yes |
| 100 | F | Forties | HM | Post graduate | Rich | Yes |
| 101 | M | Late forties | Working as driver in state govt. transport | SSLC | Middle class | Yes |

Appendix Three

# C-PROGRAM FOR CETD AND RTD MATRIX

**PROGRAM IN C LANGUAGE TO FIND THE ROW SUMS OF THE CETD MATRIX AND REFINED TIME DEPENDENT MATRIX**

```
#include<stdio.h>
#include<math.h>
void main()
{
float a[8][8],x[8][8],me[8],std[8];
float t[8],c1=0.0,c2=0.0,c=0.0,max[8];
float ce[8][8],row[8],c3=0.0,min,r,s;
  int i,j,m,n,e[8][8],p,q,u,v;
  clrscr();
  printf("enter the no. of rows and columns\n");
  scanf("%d %d",&m,&n);
  printf("enter the initial matrix\n");
  for(i=0;i<m;i++)
  {
    for(j=0;j<n;j++)
      scanf("%f",&a[i][j]);
    printf("\n");
  }
  printf("enter the time interval in each row\n");
  for(i=0;i<m;i++)
    scanf("%f",&t[i]);
  printf("enter the epsilon\n");
  scanf("%f",&s);
  for(i=0;i<m;i++)
  {for(j=0;j<n;j++)
```

```
    ce[i][j]=0.0;}
for(i=0;i<m;i++)
{for(j=0;j<n;j++)
    a[i][j]=a[i][j]/t[i];}
for(j=0;j<n;j++)
{for(i=0;i<m;i++)
    c1=c1+a[i][j];
  me[j]=c1/m;
  c1=0.0;}
for(j=0;j<n;j++)
{for(i=0;i<m;i++)
    c2=c2+pow((a[i][j]-me[j]),2);
  std[j]=sqrt(c2/m);
  c2=0.0;}
for(i=0;i<m;i++)
{
  for(j=0;j<n;j++)
    x[i][j]=(a[i][j]-me[j])/std[j];
}
for(i=0;i<m;i++)
{
  c=fabs(x[i][0]);
  for(j=0;j<n;j++)
  {
    if(fabs(x[i][j])>c)
    c=fabs(x[i][j]);
  }
  max[i]=c;
}
min=max[0];
for(i=1;i<m;i++)
{
  if(max[i]<min)
  min=max[i];
}
printf("ATD matrix:\n");
for(i=0;i<m;i++)
{
```

```
  for(j=0;j<n;j++)
    printf("%f\t",a[i][j]);
  printf("\n");
 }
 getch();
 printf("\n%f",min);
 printf("\n");
 for(p=0;p<m;p++)
 {
  for(q=0;q<n;q++)
  {
    printf("\n\n");
    if(fabs((fabs(x[p][q])-min))<s)
    {
      r=fabs(x[p][q]);
      printf("RTD matrix for alpha=
%f\n\n",r);
      for(i=0;i<m;i++)
      {for(j=0;j<n;j++)
        { if(x[i][j]>=r)
            e[i][j]=1;
          if(fabs(x[i][j])<r)
            e[i][j]=0;
          if(x[i][j]<=(-r))
            e[i][j]=-1;
          printf("%d\t",e[i][j]);
        }
        printf("\n");
      }
      for(u=0;u<m;u++)
      {for(v=0;v<n;v++)
          ce[u][v]=ce[u][v]+e[u][v];}
    }getch();
  }
 }
 printf("\nCETD matrix is:\n\n");
 for(i=0;i<m;i++)
 {
```

```
  for(j=0;j<n;j++)
    printf("%f\t",ce[i][j]);
  printf("\n");
 }
 getch();
 for(i=0;i<m;i++)
 {
  c3=0.0;
  for(j=0;j<n;j++)
    c3=c3+ce[i][j];
  row[i]=c3;
 }
 printf("\nRow sums of the CETD matrix;\n");
 for(i=0;i<m;i++)
  printf("row %d= %f\n",i+1,row[i]);
 getch();
}
```

Appendix Four

# C-PROGRAM FOR FRM

**PROGRAM IN C LANGUAGE TO FIND THE FIXED POINT FOR A GIVEN INPUT VECTOR FOR FRM**

```
#include<stdio.h>
#include<math.h>
void main()
{
 int a[8][8],b[8][8],c[8][8],y[8],y1[8],e[8][8];
 int x[8],x1[8],i,j,k,t=1,t1=1,d,c1,c2,s=0,u,v,u1,v1,s1=0,p,q,p1;
 clrscr();
 printf("enter the number of rows and columns of the matrix: ");
 scanf("%d %d",&p,&q);
 printf("enter the initial matrix:\n\n");
 for(i=0;i<p;i++)
 {
  for(j=0;j<q;j++)
    scanf("%d",&a[i][j]);
  printf("\n");
 }
 printf("\nThe given matrix is:\n");
 for(i=0;i<p;i++)
 {
  for(j=0;j<q;j++)
    printf("%d\t",a[i][j]);
  printf("\n");
 }
 getch();
 printf("enter the dimension of the input vector: ");
 scanf("%d",&p1);
 if(p1==p)
 {
  for(i=0;i<p;i++)
  {
    for(j=0;j<q;j++)
      e[i][j]=a[i][j];
  }
```

```
}
else
{
  for(i=0;i<p;i++)
   {
    for(j=0;j<q;j++)
      e[j][i]=a[i][j];
   }
q=p;
p=p1;
}
printf("\nenter the input vector:\n");
for(i=0;i<p;i++)
  scanf("%d",&x[i]);
for(i=0;i<p;i++)
{
  x1[i]=x[i];
}
for(j=0;j<q;j++)
{
  c2=0;
  for(i=0;i<p;i++)
    c2=c2+(x1[i]*e[i][j]);
  y[j]=c2;
}
for(i=0;i<q;i++)
{
  if(y[i]>0)
    y[i]=1;
  if(y[i]<0)
    y[i]=-1;
  if(y[i]==0)
    y[i]=0;
}
for(i=0;i<q;i++)
  y1[i]=y[i];
for(k=1;(t==1)&&(t1==1);k++)
{
  for(i=0;i<p;i++)
    b[s1][i]=x[i];
  s1++;
  for(i=0;i<p;i++)
   {
    c1=0;
    for(j=0;j<q;j++)
      c1=c1+(y[j]*e[i][j]);
    x[i]=c1;
   }
```

```
for(i=0;i<p;i++)
{
 if(x[i]>0)
   x[i]=1;
 if(x[i]<0)
   x[i]=-1;
 if(x[i]==0)
   x[i]=0;
}
printf("\n\n");
for(i=0;i<p;i++)
{
 if(abs(x[i]+x1[i])<=1)
   x[i]=x[i]+x1[i];
}
for(j=0;j<s1;j++)
{
 t1=0;
 for(i=0;i<p;i++)
 {
   if(b[j][i]!=x[i])
     t1=1;
 }
 if(t1==0)
 break;
}
for(i=0;i<q;i++)
 c[s][i]=y[i];
s++;
for(j=0;j<q;j++)
{
 d=0;
 for(i=0;i<p;i++)
   d=d+(x[i]*e[i][j]);
 y[j]=d;
}
for(i=0;i<q;i++)
{
 if(y[i]>0)
   y[i]=1;
 if(y[i]<0)
   y[i]=-1;
 if(y[i]==0)
   y[i]=0;
}
printf("\n\n");
for(i=0;i<q;i++)
{
```

```
    if(abs(y[i]+y1[i])<=1)
      y[i]=y[i]+y1[i];
  }
  for(j=0;j<s;j++)
  {
    t=0;
    for(i=0;i<q;i++)
    {
      if(c[j][i]!=y[i])
        t=1;
    }
    if(t==0)
    break;
  }
}
printf("\nthe input vector is:\n");
for(i=0;i<p;i++)
  printf("%d\t",x1[i]);
for(i=0;i<p;i++)
  b[s1][i]=x[i];
for(i=0;i<q;i++)
  c[s][i]=y[i];
if(t1==0)
{
  for(j=0;j<=s1;j++)
  {
    v1=0;
    for(i=0;i<p;i++)
    {
      if(b[j][i]!=x[i])
        v1=1;
    }
    if(v1==0)
    break;
  }
  u1=j;
  printf("\n");
  for(j=u1;j<=s1;j++)
  {
    printf("\nthe input vector for the %d iteration:",j+1);
    for(i=0;i<p;i++)
      printf("%d\t",b[j][i]);
    printf("\n");
  }
}
else
{
  for(j=0;j<=s;j++)
```

```
  {
   v=0;
   for(i=0;i<q;i++)
   {
    if(c[j][i]!=y[i])
      v=1;
   }
   if(v==0)
   break;
  }
  u=j;
  printf("\n");
  for(j=u;j<=s;j++)
  {
   printf("\nthe input vector for the %d iteration:",j+1);
   for(i=0;i<q;i++)
    printf("%d\t",c[j][i]);
   printf("\n");
  }
 }
}
```

Appendix Five

# C-PROGRAM FOR CFRM

**PROGRAM IN C LANGUAGE TO FIND THE FIXED POINT FOR A GIVEN INPUT VECTOR FOR CFRM**

```
#include<stdio.h>
#include<math.h>
void main()
{
 int a[8][8],b[8][8],c[8][8],y[8],y1[8],e[8][8],m[8][8];
 int x[8],x1[8],i,j,k,t=1,t1=1,d,c1,c2,s=0,u,v,u1,v1,s1=0;
 int p,q,p1,l;
 clrscr();
 printf("enter the number of rows and columns of the matrix: ");
 scanf("%d %d",&p,&q);
 printf("enter the number of experts: ");
 scanf("%d",&l);
 for(i=0;i<p;i++)
 {
  for(j=0;j<q;j++)
   a[i][j]=0;
 }
 for(k=0;k<l;k++)
 {
  printf("enter the %d expert's matrix\n",k+1);
  for(i=0;i<p;i++)
  {
   for(j=0;j<q;j++)
   {
     scanf("%d",&m[i][j]);
     a[i][j]=a[i][j]+m[i][j];
   }
  }
 }
 printf("\nThe given matrix is:\n");
 for(i=0;i<p;i++)
 {
  for(j=0;j<q;j++)
   printf("%d\t",a[i][j]);
  printf("\n");
```

```
}
getch();
printf("enter the dimension of the input vector: ");
scanf("%d",&p1);
if(p1==p)
{
  for(i=0;i<p;i++)
  {
    for(j=0;j<q;j++)
      e[i][j]=a[i][j];
  }
}
else
{
  for(i=0;i<p;i++)
  {
    for(j=0;j<q;j++)
      e[j][i]=a[i][j];
  }
q=p;
p=p1;
}
printf("\nenter the input vector:\n");
for(i=0;i<p;i++)
  scanf("%d",&x[i]);
for(i=0;i<p;i++)
{
  x1[i]=x[i];
}
for(j=0;j<q;j++)
{
  c2=0;
  for(i=0;i<p;i++)
    c2=c2+(x1[i]*e[i][j]);
  y[j]=c2;
}
for(i=0;i<q;i++)
{
  if(y[i]>0)
    y[i]=1;
  if(y[i]<0)
    y[i]=-1;
  if(y[i]==0)
    y[i]=0;
}
for(i=0;i<q;i++)
  y1[i]=y[i];
for(k=1;(t==1)&&(t1==1);k++)
```

```
{
 for(i=0;i<p;i++)
   b[s1][i]=x[i];
 s1++;
 for(i=0;i<p;i++)
 {
   c1=0;
   for(j=0;j<q;j++)
     c1=c1+(y[j]*e[i][j]);
   x[i]=c1;
 }
 for(i=0;i<p;i++)
 {
   if(x[i]>0)
     x[i]=1;
   if(x[i]<0)
     x[i]=-1;
   if(x[i]==0)
     x[i]=0;
 }
 printf("\n\n");
 for(i=0;i<p;i++)
 {
   if(abs(x[i]+x1[i])<=1)
     x[i]=x[i]+x1[i];
 }
 for(j=0;j<s1;j++)
 {
   t1=0;
   for(i=0;i<p;i++)
   {
     if(b[j][i]!=x[i])
       t1=1;
   }
   if(t1==0)
   break;
 }
 for(i=0;i<q;i++)
   c[s][i]=y[i];
 s++;
 for(j=0;j<q;j++)
 {
   d=0;
   for(i=0;i<p;i++)
     d=d+(x[i]*e[i][j]);
   y[j]=d;
 }
 for(i=0;i<q;i++)
```

```
  {
   if(y[i]>0)
     y[i]=1;
   if(y[i]<0)
     y[i]=-1;
   if(y[i]==0)
     y[i]=0;
  }
  printf("\n\n");
  for(i=0;i<q;i++)
  {
   if(abs(y[i]+y1[i])<=1)
     y[i]=y[i]+y1[i];
  }
  for(j=0;j<s;j++)
  {
   t=0;
   for(i=0;i<q;i++)
   {
    if(c[j][i]!=y[i])
      t=1;
   }
   if(t==0)
   break;
  }
}
printf("\nthe input vector is:\n");
for(i=0;i<p;i++)
  printf("%d\t",x1[i]);
for(i=0;i<p;i++)
  b[s1][i]=x[i];
for(i=0;i<q;i++)
  c[s][i]=y[i];
if(t1==0)
{
  for(j=0;j<=s1;j++)
  {
   v1=0;
   for(i=0;i<p;i++)
   {
    if(b[j][i]!=x[i])
      v1=1;
   }
   if(v1==0)
   break;
  }
  u1=j;
  printf("\n");
```

```
  for(j=u1;j<=s1;j++)
  {printf("\nthe input vector for the %d iteration:",j+1);
   for(i=0;i<p;i++)
     printf("%d\t",b[j][i]);
   printf("\n");
  }
 }
 else
 {for(j=0;j<=s;j++)
  {v=0;
   for(i=0;i<q;i++)
   {
    if(c[j][i]!=y[i])
      v=1;
   }
   if(v==0)
   break;
  }
  u=j;
  printf("\n");
  for(j=u;j<=s;j++)
  {
   printf("\nthe input vector for the %d iteration:",j+1);
   for(i=0;i<q;i++)
     printf("%d\t",c[j][i]);
   printf("\n");
  }}
}
```

Appendix Six

# C-PROGRAM FOR BAM

**PROGRAM IN C LANGUAGE TO FIND THE FIXED POINT OF THE GIVEN INPUT VECTOR IN BIDIRECTIONAL ASSOCIATIVE MEMORIES MATRIX**

```
#include<stdio.h>
void main()
{
  int m[8][8],x[8],sx[8],y[8],sy[8];
int x1[8],y1[8],sx1[8],sy1[8],i,j,a,b,t;
  int p,q,u=0,v=0;
  clrscr();
  printf("enter the no. of rows\n");
  scanf("%d",&p);
  printf("enter the no. of columns\n");
  scanf("%d",&q);
  printf("enter the matrix\n");
  for(i=0;i<p;i++)
  {
    for(j=0;j<q;j++)
      scanf("%d",&m[i][j]);
    printf("\n");
  }
  printf("enter the input vector\n");
  for(i=0;i<p;i++)
    scanf("%d",&x[i]);
  printf("\n");
  printf("\nenter the initial signal fn\n");
  for(i=0;i<q;i++)
    scanf("%d",&sy[i]);
  for(i=0;i<p;i++)
  {
    if(x[i]<=0)
      sx[i]=0;
    if(x[i]>0)
      sx[i]=1;
  }
  for(t=0;(u==0)||(v==0);t++)
```

```
{
 u=1;
 v=1;
 printf("\n\n");
 for(j=0;j<q;j++)
 {
   a=0;
   for(i=0;i<p;i++)
     a=a+sx[i]*m[i][j];
   y1[j]=a;
 }
 for(i=0;i<q;i++)
 {
   if(y1[i]<0)
     sy1[i]=0;
   if(y1[i]==0)
     sy1[i]=sy[i];
   if(y1[i]>0)
     sy1[i]=1;
 }
 for(i=0;i<p;i++)
 {
   b=0;
   for(j=0;j<q;j++)
     b=b+sy1[j]*m[i][j];
   x1[i]=b;
 }
 for(i=0;i<p;i++)
 {
   if(x1[i]<0)
     sx1[i]=0;
   if(x1[i]==0)
     sx1[i]=sx[i];
   if(x1[i]>0)
     sx1[i]=1;
 }
 for(i=0;i<p;i++)
 {
   if(sx[i]!=sx1[i])
     u=0;
   sx[i]=sx1[i];
 }
 for(i=0;i<q;i++)
 {
   if(sy[i]!=sy1[i])
     v=0;
   sy[i]=sy1[i];
 }
```

```
 }
 printf("\ntotal no of iterations\n");
 printf("%d",2*t);
 printf("\n\nthe fixed points are:\n");
 for(i=0;i<p;i++)
   printf("%d\t",sx1[i]);
 printf("\n");
 for(i=0;i<q;i++)
   printf("%d\t",sy1[i]);
}
```

# FURTHER READING

# INDEX

**A**

**B**

**C**

**D**

**E**

**F**

**H**

**I**

**L**

**M**

**N**

**R**

**S**

**U**

# About the Authors

**Dr.W.B.Vasantha Kandasamy** is an Associate Professor in the Department of Mathematics, Indian Institute of Technology Madras, Chennai, where she lives with her husband Dr.K.Kandasamy and daughters Meena and Kama. Her current interests include Smarandache algebraic structures, fuzzy theory, coding/ communication theory. In the past decade she has guided 11 Ph.D. scholars in the different fields of non-associative algebras, algebraic coding theory, transportation theory, fuzzy groups, and applications of fuzzy theory of the problems faced in chemical industries and cement industries. Currently, four Ph.D. scholars are working under her guidance.

She has to her credit 612 research papers of which 209 are individually authored. Apart from this, she and her students have presented around 329 papers in national and international conferences. She teaches both undergraduate and post-graduate students and has guided over 45 M.Sc. and M.Tech. projects. She has worked in collaboration projects with the Indian Space Research Organization and with the Tamil Nadu State AIDS Control Society. This is her 23rd book.

She can be contacted at vasantha@iitm.ac.in
You can visit her work on the web at: http://mat.iitm.ac.in/~wbv

---

**Dr.Florentin Smarandache** is an Associate Professor of Mathematics at the University of New Mexico in USA. He published over 75 books and 100 articles and notes in mathematics, physics, philosophy, psychology, literature, rebus. In mathematics his research is in number theory, non-Euclidean geometry, synthetic geometry, algebraic structures, statistics, neutrosophic logic and set (generalizations of fuzzy logic and set respectively), neutrosophic probability (generalization of classical and imprecise probability). Also, small contributions to nuclear and particle physics, information fusion, neutrosophy (a generalization of dialectics), law of sensations and stimuli, etc.).

He can be contacted at smarand@unm.edu